\documentclass{article}

\usepackage{a4wide,cite,latexsym,amsfonts,amssymb,exscale,enumerate}
\usepackage[centertags,sumlimits,intlimits,namelimits,reqno]{amsmath}
\usepackage{amsthm}
\usepackage{hyperref}

\usepackage[all]{xy}
\usepackage{fancyheadings}
\def\theauthor{}
\def\empty{}
\def\theaffiliation{}
\def\preprint#1{
  \thispagestyle{plain}
  \def\theauthor{#1}
  \ifx\theauthor\empty
  \else
    \begin{flushright}{\small #1\par}\end{flushright}
  \fi
  \begin{center}}
\def\title#1{
  {\LARGE #1\par}\vskip 1em}
\def\author#1{
  \ifx\theaffiliation\empty
  \else
    \par
  \fi
  \def\theauthor{#1}\def\theaffiliation{}}
\def\email#1{
  \vskip 1em{\large\theauthor\footnote{\small email: {\tt #1}}\par}\vskip .5em}
\def\affiliation#1{
  \ifx\theaffiliation\empty
    \def\theaffiliation{second}
  \else
    \par and\par
  \fi
  {\small\sl #1}}
\def\date#1{
  \vskip 1em{(#1)\par}\end{center}\vskip 2em}

\newcommand{\BOX}{\hbox {$\sqcap$ \kern -1em $\sqcup$}}
\newcommand{\Hom}{{\rm Hom}}
\renewcommand{\to}{\rightarrow}
\newcommand{\To}{\Rightarrow}
\newcommand{\TO}{\Rrightarrow}
\newcommand{\maps}{\colon}
\newcommand{\op}{{\rm op}}
\newcommand{\T}{{\mathbb T}}
\newcommand{\G}{{\mathbb G}}
\newcommand{\ten}{\otimes }

\newcommand{\cat}[1]{\ensuremath{\mbox{\bfseries {\upshape {#1}}}}}

\newtheorem{thm}{Theorem}    
\newtheorem{cor}[thm]{Corollary}

\newtheorem{lem}[thm]{Lemma}
\newtheorem{prop}[thm]{Proposition}
\newtheorem{defn}[thm]{Definition}
\newcommand{\Proof}{\noindent {\bf Proof. }}

        \newcommand{\be}{\begin{equation}}
        \newcommand{\ee}{\end{equation}}
        \newcommand{\ba}{\begin{eqnarray}}
        \newcommand{\ea}{\end{eqnarray}}
        \newcommand{\ban}{\begin{eqnarray*}}
        \newcommand{\ean}{\end{eqnarray*}}

\newcommand{\adjunction}[4]{%
  \ensuremath{\xymatrix{%
    {#1} \ar@<4pt>[r]^{{#3}} \ar@{}[r]|{\scriptscriptstyle{\bot}} &%
    {#2} \ar@<4pt>[l]^{{#4}}%
  }}%
}

\newcommand{\radjunction}[4]{%
  \ensuremath{\xymatrix{%
    {#1} \ar@<4pt>[r]^{{#3}} \ar@{}[r]|{\scriptscriptstyle{\top}} &%
    {#2} \ar@<4pt>[l]^{{#4}}%
  }}%
}

\newcommand{\ambijunction}[4]{%
  \ensuremath{\xymatrix{%
    {#1} \ar@<4pt>[r]^{{#3}} \ar@{}[r]|{\scriptscriptstyle{\top \;\bot}} &%
    {#2} \ar@<4pt>[l]^{{#4}}%
  }}%
}

\newcommand{\padjunction}[4]{%
  \ensuremath{\xymatrix{%
    {#1} \ar@<4pt>[r]^{{#3}} \ar@{}[r]|{\scriptscriptstyle{\bot_p}} &%
    {#2} \ar@<4pt>[l]^{{#4}}%
  }}%
}

\newcommand{\pradjunction}[4]{%
  \ensuremath{\xymatrix{%
    {#1} \ar@<4pt>[r]^{{#3}} \ar@{}[r]|{\scriptscriptstyle{\top_p}} &%
    {#2} \ar@<4pt>[l]^{{#4}}%
  }}%
}

\newcommand{\pseudoambijunction}[4]{%
  \ensuremath{\xymatrix{%
    {#1} \ar@<4pt>[r]^{{#3}} \ar@{}[r]|{\scriptscriptstyle{\top_p \;\bot_p}} &%
    {#2} \ar@<4pt>[l]^{{#4}}%
  }}%
}

\SelectTips{cm}{}

\begin{document}
\preprint{\medskip}

\title{Frobenius algebras and ambidextrous adjunctions}

\author{Aaron D.\ Lauda}
\email{A.Lauda@dpmms.cam.ac.uk}
\affiliation{Department of Pure Mathematics and Mathematical Statistics,\\
  University of Cambridge, Cambridge CB3 0WB, United Kingdom}


\date{January 19, 2006}

\begin{abstract}
In this paper we explain the relationship between Frobenius
objects in monoidal categories and adjunctions in 2-categories.
Specifically, we show that every Frobenius object in a monoidal
category $M$ arises from an ambijunction (simultaneous left and
right adjoints) in some 2-category $\mathcal{D}$ into which $M$
fully and faithfully embeds. Since a 2D topological quantum field
theory is equivalent to a commutative Frobenius algebra, this
result also shows that every 2D TQFT is obtained from an
ambijunction in some 2-category. Our theorem is proved by
extending the theory of adjoint monads to the context of an
arbitrary 2-category and utilizing the free completion under
Eilenberg-Moore objects.  We then categorify this theorem by
replacing the monoidal category $M$ with a semistrict monoidal
2-category $M$, and replacing the 2-category $\mathcal{D}$ into
which it embeds by a semistrict 3-category.  To state this more
powerful result, we must first define the notion of a `Frobenius
pseudomonoid', which categorifies that of a Frobenius object. We
then define the notion of a `pseudo ambijunction', categorifying
that of an ambijunction. In each case, the idea is that all the
usual axioms now hold only up to coherent isomorphism.  Finally,
we show that every Frobenius pseudomonoid in a semistrict monoidal
2-category arises from a pseudo ambijunction in some semistrict
3-category.
\end{abstract}

\section{Introduction}

In this paper we aim to illuminate the relationship between
Frobenius objects in monoidal categories and adjunctions in
2-categories.  One of the results we prove is that:
\begin{quote}
{\it Every Frobenius object in \emph{any} monoidal category $M$
arises from simultaneous left and right adjoints in some
2-category into which $M$ fully and faithfully embeds.}
\end{quote}
To indicate the \textit{two-handedness} of these simultaneous left
and right adjoints we refer to them as ambidextrous adjunctions
following Baez~\cite{ThisWeek}. We sometimes refer to an
ambidextrous adjunction as an {\em ambijunction} for short.

Intuitively, the relationship between Frobenius objects and
adjunctions can best be understood geometrically. This geometry
arises naturally from the language of 2-categorical string
diagrams~\cite{js2,Street}. In string diagram notation, objects
$A$ and $B$ of the 2-category $\mathcal{D}$ are depicted as
2-dimensional regions which we sometimes shade to differentiate
between different objects:
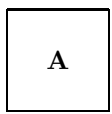
\[
\xy  0;/r.16pc/: 
    (0,0)*++{\textbf{A}};
    (-10,10)*{}="x1";
      (-10,-10)*{}="x2";
      (10,10)*{}="x4";
      (10,-10)*{}="x3";
        "x1";"x2"; **\dir{-};
        "x2";"x3"; **\dir{-};
        "x3";"x4"; **\dir{-};
        "x4";"x1"; **\dir{-};
\endxy
\qquad \qquad \qquad
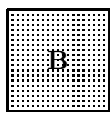
\xy 0;/r.16pc/: 
    (0,0)*++{\textbf{B}};
    (-10,10)*{}="x1";
      (-10,-10)*{}="x2";
      (10,10)*{}="x4";
      (10,-10)*{}="x3";
        "x1";"x2"; **\dir{-};
        "x2";"x3"; **\dir{-};
        "x3";"x4"; **\dir{-};
        "x4";"x1"; **\dir{-};
        (10,9)*{};  (-10,9)*{} **\dir{.};
        (10,8)*{};(-10,8)    **\dir{.};
        (10,7)*{};(-10,7)    **\dir{.};
        (10,6)*{};(-10,6)    **\dir{.};
        (10,5)*{};(-10,5)    **\dir{.};
        (10,4)*{};(-10,4)    **\dir{.};
        (10,3)*{};(-10,3)    **\dir{.};
        (10,2)*{};(-10,2)    **\dir{.};
        (10,1)*{};(-10,1)    **\dir{.};
        (10,0)*{};(-10,0)    **\dir{.};
        (10,-1)*{};(-10,-1)  **\dir{.};
        (10,-2)*{};(-10,-2)  **\dir{.};
        (10,-3)*{};(-10,-3)  **\dir{.};
        (10,-4)*{};(-10,-4)   **\dir{.};
        (10,-5)*{};(-10,-5)   **\dir{.};
        (10,-6)*{};(-10,-6)   **\dir{.};
        (10,-7)*{};(-10,-7)   **\dir{.};
        (10,-8)*{};(-10,-8)   **\dir{.};
        (10,-9)*{};(-10,-9)   **\dir{.};
\endxy
 \]
The morphisms of $\mathcal{D}$ are depicted as one dimensional
edges. Thus, if $L\maps A \to B$ and $R\maps B \to A$ are
morphisms in $\mathcal{D}$, then they are depicted as follows:
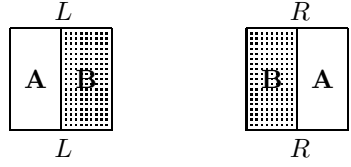
\[
\xy 0;/r.16pc/:
    (-5,0)*++{\textbf{A}};
    (5,0)*++{\textbf{B}};
    (0,10)*{ }="Top";
    (0,-10)*{ }="Bot";;
    "Top";"Bot" **\dir{-};
    (-10,10)*{}="x1";
      (-10,-10)*{}="x2";
      (10,10)*{}="x4";
      (10,-10)*{}="x3";
        "x1";"x2"; **\dir{-};
        "x2";"x3"; **\dir{-};
        "x3";"x4"; **\dir{-};
        "x4";"x1"; **\dir{-};
        (.5,13.5)*{L};
        (.5,-13.5)*{L};
        (10,9)*{};  (0,9)*{} **\dir{.};
        (10,8)*{};(0,8)    **\dir{.};
        (10,7)*{};(0,7)    **\dir{.};
        (10,6)*{};(0,6)    **\dir{.};
        (10,5)*{};(0,5)    **\dir{.};
        (10,4)*{};(0,4)    **\dir{.};
        (10,3)*{};(0,3)    **\dir{.};
        (10,2)*{};(0,2)    **\dir{.};
        (10,1)*{};(0,1)    **\dir{.};
        (10,0)*{};(0,0)    **\dir{.};
        (10,-1)*{};(0,-1)  **\dir{.};
        (10,-2)*{};(0,-2)  **\dir{.};
        (10,-3)*{};(0,-3)  **\dir{.};
        (10,-4)*{};(0,-4)   **\dir{.};
        (10,-5)*{};(0,-5)   **\dir{.};
        (10,-6)*{};(0,-6)   **\dir{.};
        (10,-7)*{};(0,-7)   **\dir{.};
        (10,-8)*{};(0,-8)   **\dir{.};
        (10,-9)*{};(0,-9)   **\dir{.};
\endxy
\qquad \qquad \quad
\xy 0;/r.16pc/:
    (-5,0)*++{\textbf{B}};
    (5,0)*++{\textbf{A}};
    (0,10)*{ }="Top";
    (0,-10)*{ }="Bot";;
    "Top";"Bot" **\dir{-};
    (-10,10)*{}="x1";
      (-10,-10)*{}="x2";
      (10,10)*{}="x4";
      (10,-10)*{}="x3";
        "x1";"x2"; **\dir{-};
        "x2";"x3"; **\dir{-};
        "x3";"x4"; **\dir{-};
        "x4";"x1"; **\dir{-};
        (.5,13.5)*{R};
        (.5,-13.5)*{R};
        (-10,9)*{};  (0,9)*{} **\dir{.};
        (-10,8)*{};(0,8)    **\dir{.};
        (-10,7)*{};(0,7)    **\dir{.};
        (-10,6)*{};(0,6)    **\dir{.};
        (-10,5)*{};(0,5)    **\dir{.};
        (-10,4)*{};(0,4)    **\dir{.};
        (-10,3)*{};(0,3)    **\dir{.};
        (-10,2)*{};(0,2)    **\dir{.};
        (-10,1)*{};(0,1)    **\dir{.};
        (-10,0)*{};(0,0)    **\dir{.};
        (-10,-1)*{};(0,-1)  **\dir{.};
        (-10,-2)*{};(0,-2)  **\dir{.};
        (-10,-3)*{};(0,-3)  **\dir{.};
        (-10,-4)*{};(0,-4)   **\dir{.};
        (-10,-5)*{};(0,-5)   **\dir{.};
        (-10,-6)*{};(0,-6)   **\dir{.};
        (-10,-7)*{};(0,-7)   **\dir{.};
        (-10,-8)*{};(0,-8)   **\dir{.};
        (-10,-9)*{};(0,-9)   **\dir{.};
\endxy
 \]
and their composite $RL\maps A \to A$ as:
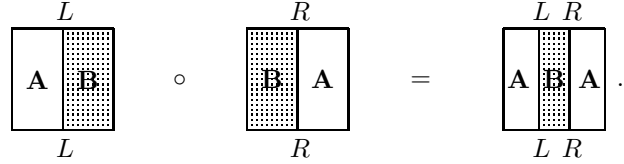
\[
\xy 0;/r.16pc/:
    (-5,0)*++{\textbf{A}};
    (5,0)*++{\textbf{B}};
    (0,10)*{ }="Top";
    (0,-10)*{ }="Bot";;
    "Top";"Bot" **\dir{-};
    (-10,10)*{}="x1";
      (-10,-10)*{}="x2";
      (10,10)*{}="x4";
      (10,-10)*{}="x3";
        "x1";"x2"; **\dir{-};
        "x2";"x3"; **\dir{-};
        "x3";"x4"; **\dir{-};
        "x4";"x1"; **\dir{-};
        (.5,13.5)*{L};
        (.5,-13.5)*{L};
        (10,9)*{};  (0,9)*{} **\dir{.};
        (10,8)*{};(0,8)    **\dir{.};
        (10,7)*{};(0,7)    **\dir{.};
        (10,6)*{};(0,6)    **\dir{.};
        (10,5)*{};(0,5)    **\dir{.};
        (10,4)*{};(0,4)    **\dir{.};
        (10,3)*{};(0,3)    **\dir{.};
        (10,2)*{};(0,2)    **\dir{.};
        (10,1)*{};(0,1)    **\dir{.};
        (10,0)*{};(0,0)    **\dir{.};
        (10,-1)*{};(0,-1)  **\dir{.};
        (10,-2)*{};(0,-2)  **\dir{.};
        (10,-3)*{};(0,-3)  **\dir{.};
        (10,-4)*{};(0,-4)   **\dir{.};
        (10,-5)*{};(0,-5)   **\dir{.};
        (10,-6)*{};(0,-6)   **\dir{.};
        (10,-7)*{};(0,-7)   **\dir{.};
        (10,-8)*{};(0,-8)   **\dir{.};
        (10,-9)*{};(0,-9)   **\dir{.};
\endxy
\qquad \circ \qquad
\xy 0;/r.16pc/:
    (-5,0)*++{\textbf{B}};
    (5,0)*++{\textbf{A}};
    (0,10)*{ }="Top";
    (0,-10)*{ }="Bot";;
    "Top";"Bot" **\dir{-};
    (-10,10)*{}="x1";
      (-10,-10)*{}="x2";
      (10,10)*{}="x4";
      (10,-10)*{}="x3";
        "x1";"x2"; **\dir{-};
        "x2";"x3"; **\dir{-};
        "x3";"x4"; **\dir{-};
        "x4";"x1"; **\dir{-};
        (.5,13.5)*{R};
        (.5,-13.5)*{R};
        (-10,9)*{};  (0,9)*{} **\dir{.};
        (-10,8)*{};(0,8)    **\dir{.};
        (-10,7)*{};(0,7)    **\dir{.};
        (-10,6)*{};(0,6)    **\dir{.};
        (-10,5)*{};(0,5)    **\dir{.};
        (-10,4)*{};(0,4)    **\dir{.};
        (-10,3)*{};(0,3)    **\dir{.};
        (-10,2)*{};(0,2)    **\dir{.};
        (-10,1)*{};(0,1)    **\dir{.};
        (-10,0)*{};(0,0)    **\dir{.};
        (-10,-1)*{};(0,-1)  **\dir{.};
        (-10,-2)*{};(0,-2)  **\dir{.};
        (-10,-3)*{};(0,-3)  **\dir{.};
        (-10,-4)*{};(0,-4)   **\dir{.};
        (-10,-5)*{};(0,-5)   **\dir{.};
        (-10,-6)*{};(0,-6)   **\dir{.};
        (-10,-7)*{};(0,-7)   **\dir{.};
        (-10,-8)*{};(0,-8)   **\dir{.};
        (-10,-9)*{};(0,-9)   **\dir{.};
\endxy
\qquad = \qquad
\xy 0;/r.16pc/:
    (-7,0)*++{\textbf{A}};
    (7,0)*++{\textbf{A}};
    (0,0)*++{\textbf{B}};
    (-3,10)*{ }="Top";
    (-3,-10)*{ }="Bot";
    "Top";"Bot" **\dir{-};
    (3,10)*{ }="Top";
    (3,-10)*{ }="Bot";
    "Top";"Bot" **\dir{-};
    (-10,10)*{}="x1";
      (-10,-10)*{}="x2";
      (10,10)*{}="x4";
      (10,-10)*{}="x3";
        "x1";"x2"; **\dir{-};
        "x2";"x3"; **\dir{-};
        "x3";"x4"; **\dir{-};
        "x4";"x1"; **\dir{-};
        (3.5,13.5)*{R};
        (3.5,-13.5)*{R};
        (-2.5,13.5)*{L};
        (-2.5,-13.5)*{L};
        (-3,9)*{};  (3,9)*{} **\dir{.};
        (-3,8)*{};(3,8)    **\dir{.};
        (-3,7)*{};(3,7)    **\dir{.};
        (-3,6)*{};(3,6)    **\dir{.};
        (-3,5)*{};(3,5)    **\dir{.};
        (-3,4)*{};(3,4)    **\dir{.};
        (-3,3)*{};(3,3)    **\dir{.};
        (-3,2)*{};(3,2)    **\dir{.};
        (-3,1)*{};(3,1)    **\dir{.};
        (-3,0)*{};(3,0)    **\dir{.};
        (-3,-1)*{};(3,-1)  **\dir{.};
        (-3,-2)*{};(3,-2)  **\dir{.};
        (-3,-3)*{};(3,-3)  **\dir{.};
        (-3,-4)*{};(3,-4)   **\dir{.};
        (-3,-5)*{};(3,-5)   **\dir{.};
        (-3,-6)*{};(3,-6)   **\dir{.};
        (-3,-7)*{};(3,-7)   **\dir{.};
        (-3,-8)*{};(3,-8)   **\dir{.};
        (-3,-9)*{};(3,-9)   **\dir{.};
\endxy .
 \]
As a convenient convention, the identity morphisms of objects in
$\mathcal{D}$ are not drawn.  This convention allows for the
identification:
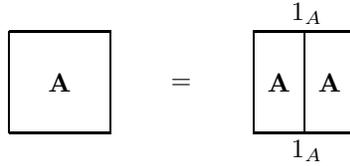
\[ 
\xy 0;/r.16pc/:
    (0,0)*++{\textbf{A}};
      (-10,10)*{}="x1";
      (-10,-10)*{}="x2";
      (10,10)*{}="x4";
      (10,-10)*{}="x3";
        "x1";"x2"; **\dir{-};
        "x2";"x3"; **\dir{-};
        "x3";"x4"; **\dir{-};
        "x4";"x1"; **\dir{-};
 \endxy
 \qquad = \qquad
 \xy 0;/r.16pc/:
    (-5,0)*++{\textbf{A}};
    (5,0)*++{\textbf{A}};
    (0,10)*{ }="Top";
    (0,-10)*{ }="Bot";;
    "Top";"Bot" **\dir{-};
    (-10,10)*{}="x1";
      (-10,-10)*{}="x2";
      (10,10)*{}="x4";
      (10,-10)*{}="x3";
        "x1";"x2"; **\dir{-};
        "x2";"x3"; **\dir{-};
        "x3";"x4"; **\dir{-};
        "x4";"x1"; **\dir{-};
        (.5,13.5)*{1_A};
        (.5,-13.5)*{1_A};
\endxy
 \]
of string diagrams.

The 2-morphisms of $\mathcal{D}$ are drawn as 0-dimensional
vertices or as small discs if we want to label them.  Hence, the
unit $i\maps 1 \To RL$ and counit $e \maps LR \To 1$ of an
adjunction $\adjunction{A}{B}{L}{R}$ are depicted as:
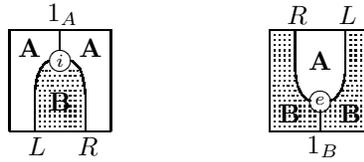
\[
\xy 0;/r.16pc/: 
    (-6,6.5)*+{\textbf{A}};
    (6,6.5)*+{\textbf{A}};
    (0,-4)*+{\textbf{B}};
    (-5,-10)*{}="TL";
    (5,-10)*{}="TR";
    (0,10)*{}="B";
    (-4.5,-13)*{L};
    (5.5,-13)*{R};
    (.5,13)*{1_A};
    (0,4)*+{\scriptstyle i}*\cir{}="m";
    (-5,-4)*{}="l";
    (5,-4)*{}="r";
    "l";"m" **\crv{(-5,5)};
    "r";"m" **\crv{(5,5)};
    "TL";"l" **\dir{-};
    "m";"B" **\dir{-};
    "r";"TR" **\dir{-};
       (-10,10)*{}="x1";
       (-10,-10)*{}="x2";
       (10,10)*{}="x4";
       (10,-10)*{}="x3";
        "x1";"x2"; **\dir{-};
        "x2";"x3"; **\dir{-};
        "x3";"x4"; **\dir{-};
        "x4";"x1"; **\dir{-};
        (-5,-9)*{};(5,-9)    **\dir{.};
        (-5,-8)*{};(5,-8)    **\dir{.};
        (-5,-7)*{};(5,-7)    **\dir{.};
        (-5,-6)*{};(5,-6)    **\dir{.};
        (-5,-5)*{};(5,-5)    **\dir{.};
        (-5,-4)*{};(5,-4)    **\dir{.};
        (-5,-3)*{};(5,-3)    **\dir{.};
        (-5,-2)*{};(5,-2)    **\dir{.};
        (-4.5,-1)*{};(4.5,-1)    **\dir{.};
        (-4.5,0)*{};(4.5,0)    **\dir{.};
        (-4,1)*{};(4,1)    **\dir{.};
        (2,2)*{};(3.5,2)    **\dir{.};
        (-3.5,2)*{};(-2,2)    **\dir{.};
        (-3,3)*{};(-2.5,3)    **\dir{.};
        (3,3)*{};(2.5,3)    **\dir{.};
        \endxy
\qquad \qquad \qquad
 \xy 0;/r.16pc/: 
    (-6,-6.5)*+{\textbf{B}};
    (6,-6.5)*+{\textbf{B}};
    (0,4)*+{\textbf{A}};
    (-5,10)*{}="TL";
    (5,10)*{}="TR";
    (0,-10)*{}="B";
    (-4.5,13)*{R};
    (5.5,13)*{L};
    (.5,-13)*{1_B};
    (0,-4)*+{\scriptstyle e}*\cir{}="m";
    (-5,4)*{}="l";
    (5,4)*{}="r";
    "l";"m" **\crv{(-5,-5)};
    "r";"m" **\crv{(5,-5)};
    "TL";"l" **\dir{-};
    "m";"B" **\dir{-};
    "r";"TR" **\dir{-};
       (-10,10)*{}="x1";
       (-10,-10)*{}="x2";
       (10,10)*{}="x4";
       (10,-10)*{}="x3";
        "x1";"x2"; **\dir{-};
        "x2";"x3"; **\dir{-};
        "x3";"x4"; **\dir{-};
        "x4";"x1"; **\dir{-};
        (-10,9)*{};(-5,9)    **\dir{.};
        (-10,8)*{};(-5,8)    **\dir{.};
        (-10,7)*{};(-5,7)    **\dir{.};
        (-10,6)*{};(-5,6)    **\dir{.};
        (-10,5)*{};(-5,5)    **\dir{.};
        (-10,4)*{};(-5,4)    **\dir{.};
        (-10,3)*{};(-5,3)    **\dir{.};
        (-10,2)*{};(-5,2)    **\dir{.};
        (-10,1)*{};(-5,1)    **\dir{.};
        (-10,0)*{};(-5,0)    **\dir{.};
        (-10,-1)*{};(-4.5,-1)    **\dir{.};
        (-10,-2)*{};(-4,-2)    **\dir{.};
        (-10,-3)*{};(-3.5,-3)    **\dir{.};
        (-10,-4)*{};(-2,-4)    **\dir{.};
        (-10,-5)*{};(-2,-5)    **\dir{.};
        (-10,-6)*{};(0,-6)    **\dir{.};
        (-10,-7)*{};(0,-7)    **\dir{.};
        (-10,-8)*{};(0,-8)    **\dir{.};
        (-10,-9)*{};(0,-9)    **\dir{.};
          (10,9)*{};(5,9)    **\dir{.};
        (10,8)*{};(5,8)    **\dir{.};
        (10,7)*{};(5,7)    **\dir{.};
        (10,6)*{};(5,6)    **\dir{.};
        (10,5)*{};(5,5)    **\dir{.};
        (10,4)*{};(5,4)    **\dir{.};
        (10,3)*{};(5,3)    **\dir{.};
        (10,2)*{};(5,2)    **\dir{.};
        (10,1)*{};(5,1)    **\dir{.};
        (10,0)*{};(5,0)    **\dir{.};
        (10,-1)*{};(4.5,-1)    **\dir{.};
         (10,-2)*{};(4.5,-2)    **\dir{.};
        (10,-3)*{};(3.5,-3)    **\dir{.};
        (10,-4)*{};(2,-4)    **\dir{.};
        (10,-5)*{};(2,-5)    **\dir{.};
        (10,-6)*{};(0,-6)    **\dir{.};
        (10,-7)*{};(0,-7)    **\dir{.};
        (10,-8)*{};(0,-8)    **\dir{.};
        (10,-9)*{};(0,-9)    **\dir{.};
        \endxy
\]
However, using the convention for the identity morphisms mentioned
above and omitting the labels we can simplify these string
diagrams as follows:
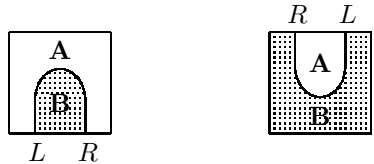
\[ 
\xy 0;/r.16pc/:
    (0,6.5)*+{\textbf{A}};
    (0,-4)*+{\textbf{B}};
    (-5,-10)*{ }="TL";
    (5,-10)*{ }="TR";
    (-4.5,-13.5)*{L};
    (5.5,-13.5)*{R};
    (-5,-4)*{}="l";
    (5,-4)*{}="r";
    "l";"r" **\crv{(-5,5)&(5,5)};
    "TL";"l" **\dir{-};
    "r";"TR" **\dir{-};
       (-10,10)*{}="x1";
       (-10,-10)*{}="x2";
       (10,10)*{}="x4";
       (10,-10)*{}="x3";
        "x1";"x2"; **\dir{-};
        "x2";"x3"; **\dir{-};
        "x3";"x4"; **\dir{-};
        "x4";"x1"; **\dir{-};
        (-5,-9)*{};(5,-9)    **\dir{.};
        (-5,-8)*{};(5,-8)    **\dir{.};
        (-5,-7)*{};(5,-7)    **\dir{.};
        (-5,-6)*{};(5,-6)    **\dir{.};
        (-5,-5)*{};(5,-5)    **\dir{.};
        (-5,-4)*{};(5,-4)    **\dir{.};
        (-5,-3)*{};(5,-3)    **\dir{.};
        (-4.5,-2)*{};(4.5,-2)    **\dir{.};
        (-4,-1)*{};(4,-1)    **\dir{.};
        (-4,0)*{};(4,0)    **\dir{.};
        (-3,1)*{};(3,1)    **\dir{.};
        (-2,2)*{};(2,2)    **\dir{.};
        \endxy
\qquad \qquad \qquad
\xy 0;/r.16pc/: 
    (0,-6.5)*+{\textbf{B}};
    (0,4)*+{\textbf{A}};
    (-5,10)*{ }="TL";
    (5,10)*{ }="TR";
    (-4.5,13.5)*{R};
    (5.5,13.5)*{L};
    (-5,4)*{}="l";
    (5,4)*{}="r";
    "l";"r" **\crv{(-5,-5)&(5,-5)};
    "TL";"l" **\dir{-};
    "r";"TR" **\dir{-};
       (-10,10)*{}="x1";
       (-10,-10)*{}="x2";
       (10,10)*{}="x4";
       (10,-10)*{}="x3";
        "x1";"x2"; **\dir{-};
        "x2";"x3"; **\dir{-};
        "x3";"x4"; **\dir{-};
        "x4";"x1"; **\dir{-};
        (-10,9)*{};(-5,9)    **\dir{.};
        (-10,8)*{};(-5,8)    **\dir{.};
        (-10,7)*{};(-5,7)    **\dir{.};
        (-10,6)*{};(-5,6)    **\dir{.};
        (-10,5)*{};(-5,5)    **\dir{.};
        (-10,4)*{};(-5,4)    **\dir{.};
        (-10,3)*{};(-5,3)    **\dir{.};
        (-10,2)*{};(-5,2)    **\dir{.};
        (-10,1)*{};(-5,1)    **\dir{.};
        (-10,0)*{};(-5,0)    **\dir{.};
          (-10,-1)*{};(-4,-1)    **\dir{.};
          (-10,-2)*{};(-3,-2)    **\dir{.};
        (10,9)*{};(5,9)    **\dir{.};
        (10,8)*{};(5,8)    **\dir{.};
        (10,7)*{};(5,7)    **\dir{.};
        (10,6)*{};(5,6)    **\dir{.};
        (10,5)*{};(5,5)    **\dir{.};
        (10,4)*{};(5,4)    **\dir{.};
        (10,3)*{};(5,3)    **\dir{.};
        (10,2)*{};(5,2)    **\dir{.};
        (10,1)*{};(5,1)    **\dir{.};
        (10,0)*{};(5,0)    **\dir{.};
          (10,-1)*{};(4,-1)    **\dir{.};
          (10,-2)*{};(3,-2)    **\dir{.};
            (-10,-9)*{};(10,-9)    **\dir{.};
            (-10,-8)*{};(10,-8)    **\dir{.};
            (-10,-7)*{};(10,-7)    **\dir{.};
            (-10,-6)*{};(10,-6)    **\dir{.};
            (-10,-5)*{};(10,-5)    **\dir{.};
            (-10,-4)*{};(10,-4)    **\dir{.};
            (-10,-3)*{};(10,-3)    **\dir{.};
        \endxy
\]
We can also express the axioms for an adjunction, often referred
to as the triangle identities or zig-zag identities, by the
following equations of string diagrams:
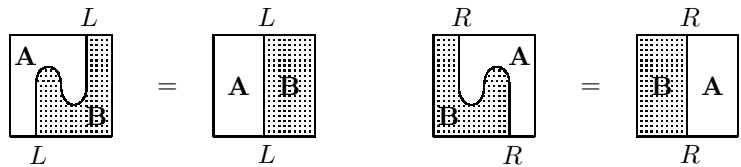
\[
 \xy 0;/r.16pc/: 
    (-7,6)*[o]++{\textbf{A}};
    (7,-6)*[o]++{\textbf{B}};
    (-5,-10)*{ }="Top";
    (5,10)*{ }="Bot";
    (-4.5,-13.5)*{L};
    (5.5,13.5)*{L};
    (-5,0)*{}="l";
    (-0,0)*{}="m";
    (5,0)*{}="r";
    "l";"m" **\crv{(-5,5)&(0,5)};
    "m";"r" **\crv{(0,-5)&(5,-5)};
    "Top";"l" **\dir{-};
    "r";"Bot" **\dir{-};
    (-10,10)*{}="x1";
      (-10,-10)*{}="x2";
      (10,10)*{}="x4";
      (10,-10)*{}="x3";
        "x1";"x2"; **\dir{-};
        "x2";"x3"; **\dir{-};
        "x3";"x4"; **\dir{-};
        "x4";"x1"; **\dir{-};
        (10,9)*{};(5,9)    **\dir{.};
        (10,8)*{};(5,8)    **\dir{.};
        (10,7)*{};(5,7)    **\dir{.};
        (10,6)*{};(5,6)    **\dir{.};
        (10,5)*{};(5,5)    **\dir{.};
        (10,4)*{};(5,4)    **\dir{.};
        (10,3)*{};(5,3)    **\dir{.};
        (10,2)*{};(5,2)    **\dir{.};
        (10,1)*{};(5,1)    **\dir{.};
        (10,0)*{};(5,0)    **\dir{.};
        (10,-1)*{};(5,-1)  **\dir{.};
        (10,-2)*{};(5,-2)  **\dir{.};
        (10,-3)*{};(5,-3)  **\dir{.};
        (-1,3)*{} ;(-4,3)      **\dir{.};
        (-1,2)*{} ;(-4,2)      **\dir{.};
        (-.5,1)*{} ;(-4.5,1)      **\dir{.};
        (0,0)*{} ;(-5,0)      **\dir{.};
        (0,-1)*{};(-5,-1)   **\dir{.};
        (0,-2)*{};(-5,-2)     **\dir{.};
        (1,-3)*{};(-5,-3)   **\dir{.};
        (10,-4)*{};(-5,-4)   **\dir{.};
        (10,-5)*{};(-5,-5)   **\dir{.};
        (10,-6)*{};(-5,-6)   **\dir{.};
        (10,-7)*{};(-5,-7)   **\dir{.};
        (10,-8)*{};(-5,-8)   **\dir{.};
        (10,-9)*{};(-5,-9)   **\dir{.};
\endxy
\quad = \quad
\xy 0;/r.16pc/: 
    (-5,0)*++{\textbf{A}};
    (5,0)*++{\textbf{B}};
    (0,10)*{ }="Top";
    (0,-10)*{ }="Bot";;
    "Top";"Bot" **\dir{-};
    (-10,10)*{}="x1";
      (-10,-10)*{}="x2";
      (10,10)*{}="x4";
      (10,-10)*{}="x3";
        "x1";"x2"; **\dir{-};
        "x2";"x3"; **\dir{-};
        "x3";"x4"; **\dir{-};
        "x4";"x1"; **\dir{-};
        (.5,13.5)*{L};
        (.5,-13.5)*{L};
        (10,9)*{};  (0,9)*{} **\dir{.};
        (10,8)*{};(0,8)    **\dir{.};
        (10,7)*{};(0,7)    **\dir{.};
        (10,6)*{};(0,6)    **\dir{.};
        (10,5)*{};(0,5)    **\dir{.};
        (10,4)*{};(0,4)    **\dir{.};
        (10,3)*{};(0,3)    **\dir{.};
        (10,2)*{};(0,2)    **\dir{.};
        (10,1)*{};(0,1)    **\dir{.};
        (10,0)*{};(0,0)    **\dir{.};
        (10,-1)*{};(0,-1)  **\dir{.};
        (10,-2)*{};(0,-2)  **\dir{.};
        (10,-3)*{};(0,-3)  **\dir{.};
        (10,-4)*{};(0,-4)   **\dir{.};
        (10,-5)*{};(0,-5)   **\dir{.};
        (10,-6)*{};(0,-6)   **\dir{.};
        (10,-7)*{};(0,-7)   **\dir{.};
        (10,-8)*{};(0,-8)   **\dir{.};
        (10,-9)*{};(0,-9)   **\dir{.};
\endxy
 \qquad  \qquad
\xy 0;/r.16pc/: 
    (-7,-6)*[o]++{\textbf{B}};
    (7,6)*[o]++{\textbf{A}};
    (-4.5,13.5)*{R};
    (5.5,-13.5)*{R};
    (-5,10)*{ }="Top";
    (5,-10)*{ }="Bot";
    (-5,0)*{}="l";
    (-0,0)*{}="m";
    (5,0)*{}="r";
    "l";"m" **\crv{(-5,-5)&(0,-5)};
    "m";"r" **\crv{(0,5)&(5,5)};
    "Top";"l" **\dir{-};
    "r";"Bot" **\dir{-};
    (-10,10)*{}="x1";
      (-10,-10)*{}="x2";
      (10,10)*{}="x4";
      (10,-10)*{}="x3";
        "x1";"x2"; **\dir{-};
        "x2";"x3"; **\dir{-};
        "x3";"x4"; **\dir{-};
        "x4";"x1"; **\dir{-};
        (-10,9)*{};(-5,9)    **\dir{.};
        (-10,8)*{};(-5,8)    **\dir{.};
        (-10,7)*{};(-5,7)    **\dir{.};
        (-10,6)*{};(-5,6)    **\dir{.};
        (-10,5)*{};(-5,5)    **\dir{.};
        (-10,4)*{};(-5,4)    **\dir{.};
        (-10,3)*{};(-5,3)    **\dir{.};
        (-10,2)*{};(-5,2)    **\dir{.};
        (-10,1)*{};(-5,1)    **\dir{.};
        (-10,0)*{};(-5,0)    **\dir{.};
        (-10,-1)*{};(-5,-1)  **\dir{.};
        (-10,-2)*{};(-5,-2)  **\dir{.};
        (-10,-3)*{};(-5,-3)  **\dir{.};
        (1,3)*{} ;(4,3)      **\dir{.};
        (1,2)*{} ;(4,2)      **\dir{.};
        (.5,1)*{} ;(4.5,1)      **\dir{.};
        (0,0)*{} ;(5,0)      **\dir{.};
        (0,-1)*{};(5,-1)   **\dir{.};
        (0,-2)*{};(5,-2)     **\dir{.};
        (-1,-3)*{};(5,-3)   **\dir{.};
        (-10,-4)*{};(5,-4)   **\dir{.};
        (-10,-5)*{};(5,-5)   **\dir{.};
        (-10,-6)*{};(5,-6)   **\dir{.};
        (-10,-7)*{};(5,-7)   **\dir{.};
        (-10,-8)*{};(5,-8)   **\dir{.};
        (-10,-9)*{};(5,-9)   **\dir{.};
\endxy
\quad = \quad
\xy 0;/r.16pc/: 
    (-5,0)*++{\textbf{B}};
    (5,0)*++{\textbf{A}};
    (0,10)*{ }="Top";
    (0,-10)*{ }="Bot";;
    "Top";"Bot" **\dir{-};
    (-10,10)*{}="x1";
      (-10,-10)*{}="x2";
      (10,10)*{}="x4";
      (10,-10)*{}="x3";
        "x1";"x2"; **\dir{-};
        "x2";"x3"; **\dir{-};
        "x3";"x4"; **\dir{-};
        "x4";"x1"; **\dir{-};
        (.5,13.5)*{R};
        (.5,-13.5)*{R};
        (-10,9)*{};  (0,9)*{} **\dir{.};
        (-10,8)*{};(0,8)    **\dir{.};
        (-10,7)*{};(0,7)    **\dir{.};
        (-10,6)*{};(0,6)    **\dir{.};
        (-10,5)*{};(0,5)    **\dir{.};
        (-10,4)*{};(0,4)    **\dir{.};
        (-10,3)*{};(0,3)    **\dir{.};
        (-10,2)*{};(0,2)    **\dir{.};
        (-10,1)*{};(0,1)    **\dir{.};
        (-10,0)*{};(0,0)    **\dir{.};
        (-10,-1)*{};(0,-1)  **\dir{.};
        (-10,-2)*{};(0,-2)  **\dir{.};
        (-10,-3)*{};(0,-3)  **\dir{.};
        (-10,-4)*{};(0,-4)   **\dir{.};
        (-10,-5)*{};(0,-5)   **\dir{.};
        (-10,-6)*{};(0,-6)   **\dir{.};
        (-10,-7)*{};(0,-7)   **\dir{.};
        (-10,-8)*{};(0,-8)   **\dir{.};
        (-10,-9)*{};(0,-9)   **\dir{.};
\endxy
 \]

Early work on homological algebra~\cite{god,mac2} and monad theory
\cite{EM,Kleisli} showed that an adjunction
$\adjunction{A}{B}{L}{R}$ endows the composite morphism $RL$ with
a monoid structure in the monoidal category $\Hom(A,A)$. This
monoid in $\Hom(A,A)$ can be vividly seen using the language of
2-categorical string diagrams.  The multiplication on $RL$ is
defined using the unit for the adjunction as seen below:
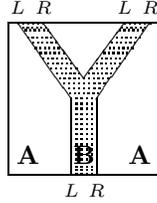
\[
\xy 0;/r.16pc/: 
    (-11,-11)*++{\textbf{A}};
    (11,-11)*++{\textbf{A}};
    (0,-11)*++{\textbf{B}};
    (-13,15)*{ }="1";
    (-8,15)*{ }="2";
    (8,15)*{ }="3";
    (13,15)*{ }="4";
    (-2.5,0)*{}="2M";
    (2.5,0)*{}="3M";
    (0,4)*{}="2X";
    (-2.5,-15)*{ }="2B";
    (2.5,-15)*{ }="3B";
    "2M";"2B" **\dir{-};
    "3M";"3B" **\dir{-};
    "2X"; "2" **\dir{-};
    "3"; "2X" **\dir{-};
    "1";"2M" **\dir{-};
    "4";"3M" **\dir{-};
    (-15,15)*{}="x1"; 
      (-15,-15)*{}="x2";
      (15,15)*{}="x4";
      (15,-15)*{}="x3";
        "x1";"x2"; **\dir{-};
        "x2";"x3"; **\dir{-};
        "x3";"x4"; **\dir{-};
        "x4";"x1"; **\dir{-}; 
        (-13,18)*{\scriptstyle L};
        (-8, 18)*{\scriptstyle R};
        (13,18)*{\scriptstyle R};
        (8,18)*{\scriptstyle L};
        (-2.5,-18)*{\scriptstyle L};
        (2.5,-18)*{\scriptstyle R};
            (-13,14.9)*{};  (-8,14.9)*{} **\dir{.};
            (-12,14)*{};  (-7.5,14)*{} **\dir{.};
            (-11.5,13.5)*{};  (-7,13.5)*{} **\dir{.};
            (-11,12)*{};  (-6.5,12)*{} **\dir{.};
            (-10,11)*{};  (-5.5,11)*{} **\dir{.};
            (-9.5,10)*{};  (-5,10)*{} **\dir{.};
            (-8.5,9)*{};  (-4,9)*{} **\dir{.};
            (-8,8)*{};  (-3.5,8)*{} **\dir{.};
            (-7,7)*{};  (-2.5,7)*{} **\dir{.};
            (-6.5,6)*{};  (-1.5,6)*{} **\dir{.};
            (-5.25,5)*{};  (-1,5)*{} **\dir{.};
            (-4.75,4)*{};  (0,4)*{} **\dir{.};
            (-4,3)*{};  (4,3)*{} **\dir{.};
            (-3,2)*{};  (3,2)*{} **\dir{.};
            (13,14.9)*{};  (8,14.9)*{} **\dir{.};
            (12,14)*{};  (7.5,14)*{} **\dir{.};
            (11.5,13.5)*{};  (7,13.5)*{} **\dir{.};
            (11,12)*{};  (6.5,12)*{} **\dir{.};
            (10,11)*{};  (5.5,11)*{} **\dir{.};
            (9.5,10)*{};  (5,10)*{} **\dir{.};
            (8.5,9)*{};  (4,9)*{} **\dir{.};
            (8,8)*{};  (3.5,8)*{} **\dir{.};
            (7,7)*{};  (2.5,7)*{} **\dir{.};
            (6.5,6)*{};  (1.5,6)*{} **\dir{.};
            (5.25,5)*{};  (1,5)*{} **\dir{.};
            (4.75,4)*{};  (0,4)*{} **\dir{.};
            (-4,3)*{};  (4,3)*{} **\dir{.};
            (-3,2)*{};  (3,2)*{} **\dir{.};
            (-2.75,1)*{};  (2.75,1)*{} **\dir{.};
            (-2.5,0)*{};  (2.5,0)*{} **\dir{.};
            (-2.5,-1)*{};  (2.5,-1)*{} **\dir{.};
            (-2.5,-2)*{};  (2.5,-2)*{} **\dir{.};
            (-2.5,-3)*{};  (2.5,-3)*{} **\dir{.};
            (-2.5,-4)*{};  (2.5,-4)*{} **\dir{.};
            (-2.5,-5)*{};  (2.5,-5)*{} **\dir{.};
            (-2.5,-6)*{};  (2.5,-6)*{} **\dir{.};
            (-2.5,-7)*{};  (2.5,-7)*{} **\dir{.};
            (-2.5,-8)*{};  (2.5,-8)*{} **\dir{.};
            (-2.5,-9)*{};  (2.5,-9)*{} **\dir{.};
            (-2.5,-10)*{};  (2.5,-10)*{} **\dir{.};
            (-2.5,-11)*{};  (2.5,-11)*{} **\dir{.};
            (-2.5,-12)*{};  (2.5,-12)*{} **\dir{.};
            (-2.5,-13.5)*{};  (2.5,-13.5)*{} **\dir{.};
            (-2.5,-14)*{};  (2.5,-14)*{} **\dir{.};
            (-2.5,-14.9)*{};  (2.5,-14.9)*{} **\dir{.};
\endxy
 \]
and the unit for multiplication is
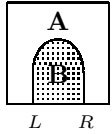
\[
\xy 0;/r.16pc/:
    (0,6.5)*+{\textbf{A}};
    (0,-4)*+{\textbf{B}};
    (-5,-10)*{ }="TL";
    (5,-10)*{ }="TR";
    (-4.5,-13.5)*{\scriptstyle L};
    (5.5,-13.5)*{\scriptstyle R};
    (-5,-4)*{}="l";
    (5,-4)*{}="r";
    "l";"r" **\crv{(-5,5)&(5,5)};
    "TL";"l" **\dir{-};
    "r";"TR" **\dir{-};
       (-10,10)*{}="x1";
       (-10,-10)*{}="x2";
       (10,10)*{}="x4";
       (10,-10)*{}="x3";
        "x1";"x2"; **\dir{-};
        "x2";"x3"; **\dir{-};
        "x3";"x4"; **\dir{-};
        "x4";"x1"; **\dir{-};
        (-5,-9)*{};(5,-9)    **\dir{.};
        (-5,-8)*{};(5,-8)    **\dir{.};
        (-5,-7)*{};(5,-7)    **\dir{.};
        (-5,-6)*{};(5,-6)    **\dir{.};
        (-5,-5)*{};(5,-5)    **\dir{.};
        (-5,-4)*{};(5,-4)    **\dir{.};
        (-5,-3)*{};(5,-3)    **\dir{.};
        (-4.5,-2)*{};(4.5,-2)    **\dir{.};
        (-4,-1)*{};(4,-1)    **\dir{.};
        (-4,0)*{};(4,0)    **\dir{.};
        (-3,1)*{};(3,1)    **\dir{.};
        (-2,2)*{};(2,2)    **\dir{.};
        \endxy
\]
the unit of the adjunction.  The associativity axiom:
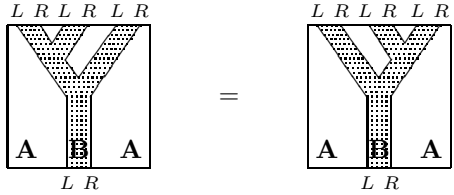
\[
\xy 0;/r.15pc/:
    (-11,-11)*++{\textbf{A}};
    (11,-11)*++{\textbf{A}};
    (0,-11)*++{\textbf{B}};
    (-13,15)*{ }="1";
    (-8,15)*{ }="2";
    (-3,15)*{ }="V1";
    (2,15)*{ }="V2";
    (8,15)*{ }="3";
    (13,15)*{ }="4";
    (-2.75,7.75)*{}="V2B";
    (-5.5,11.25)*{}="V1B";
    (-2.5,0)*{}="2M";
    (2.5,0)*{}="3M";
    (0,4)*{}="2X";
    (-2.5,-15)*{ }="2B";
    (2.5,-15)*{ }="3B";
    "V2";"V2B" **\dir{-};
    "2M";"2B" **\dir{-};
    "3M";"3B" **\dir{-};
    "V1"; "V1B" **\dir{-};
    "V1B"; "2" **\dir{-};
    "V2B"; "2X" **\dir{-};
    "3"; "2X" **\dir{-};
    "1";"2M" **\dir{-};
    "4";"3M" **\dir{-};
    (-15,15)*{}="x1"; 
      (-15,-15)*{}="x2";
      (15,15)*{}="x4";
      (15,-15)*{}="x3";
        "x1";"x2"; **\dir{-};
        "x2";"x3"; **\dir{-};
        "x3";"x4"; **\dir{-};
        "x4";"x1"; **\dir{-}; 
        (-13,18)*{\scriptstyle L};
        (-8, 18)*{\scriptstyle R};
        (-3,18)*{\scriptstyle L};
        (2, 18)*{\scriptstyle R};
        (13,18)*{\scriptstyle R};
        (8,18)*{\scriptstyle L};
        (-2.5,-18)*{\scriptstyle L};
        (2.5,-18)*{\scriptstyle R};
            (-13,14.9)*{};  (-8,14.9)*{} **\dir{.};
            (-12,14)*{};  (-7.5,14)*{} **\dir{.};
            (-11.5,13)*{};  (-7,13)*{} **\dir{.};
            (-11,12)*{};  (-6.5,12)*{} **\dir{.};
            (-3,14.9)*{};  (2,14.9)*{} **\dir{.};
            (-3.5,14)*{};  (1,14)*{} **\dir{.};
            (-4,13)*{};  (.5,13)*{} **\dir{.};
            (-4.5,12)*{};  (0,12)*{} **\dir{.};
            (-10,11)*{};  (-1,11)*{} **\dir{.};
            (-9.5,10)*{};  (-1.5,10)*{} **\dir{.};
            (-8.5,9)*{};  (-2,9)*{} **\dir{.};
            (-8,8)*{};  (-3.5,8)*{} **\dir{.};
            (-7,7)*{};  (-2.5,7)*{} **\dir{.};
            (-6.5,6)*{};  (-1.5,6)*{} **\dir{.};
            (-5.25,5)*{};  (-1,5)*{} **\dir{.};
            (-4.75,4)*{};  (0,4)*{} **\dir{.};
            (-4,3)*{};  (4,3)*{} **\dir{.};
            (-3,2)*{};  (3,2)*{} **\dir{.};
            (13,14.9)*{};  (8,14.9)*{} **\dir{.};
            (12,14)*{};  (7.5,14)*{} **\dir{.};
            (11.5,13)*{};  (7,13)*{} **\dir{.};
            (11,12)*{};  (6.5,12)*{} **\dir{.};
            (10,11)*{};  (5.5,11)*{} **\dir{.};
            (9.5,10)*{};  (5,10)*{} **\dir{.};
            (8.5,9)*{};  (4,9)*{} **\dir{.};
            (8,8)*{};  (3.5,8)*{} **\dir{.};
            (7,7)*{};  (2.5,7)*{} **\dir{.};
            (6.5,6)*{};  (1.5,6)*{} **\dir{.};
            (5.25,5)*{};  (1,5)*{} **\dir{.};
            (4.75,4)*{};  (0,4)*{} **\dir{.};
            (-4,3)*{};  (4,3)*{} **\dir{.};
            (-3,2)*{};  (3,2)*{} **\dir{.};
            (-2.75,1)*{};  (2.75,1)*{} **\dir{.};
            (-2.5,0)*{};  (2.5,0)*{} **\dir{.};
            (-2.5,-1)*{};  (2.5,-1)*{} **\dir{.};
            (-2.5,-2)*{};  (2.5,-2)*{} **\dir{.};
            (-2.5,-3)*{};  (2.5,-3)*{} **\dir{.};
            (-2.5,-4)*{};  (2.5,-4)*{} **\dir{.};
            (-2.5,-5)*{};  (2.5,-5)*{} **\dir{.};
            (-2.5,-6)*{};  (2.5,-6)*{} **\dir{.};
            (-2.5,-7)*{};  (2.5,-7)*{} **\dir{.};
            (-2.5,-8)*{};  (2.5,-8)*{} **\dir{.};
            (-2.5,-9)*{};  (2.5,-9)*{} **\dir{.};
            (-2.5,-10)*{};  (2.5,-10)*{} **\dir{.};
            (-2.5,-11)*{};  (2.5,-11)*{} **\dir{.};
            (-2.5,-12)*{};  (2.5,-12)*{} **\dir{.};
            (-2.5,-13)*{};  (2.5,-13)*{} **\dir{.};
            (-2.5,-14)*{};  (2.5,-14)*{} **\dir{.};
            (-2.5,-14.9)*{};  (2.5,-14.9)*{} **\dir{.};
\endxy
\qquad = \qquad
\xy 0;/r.15pc/:
    (-11,-11)*++{\textbf{A}};
    (11,-11)*++{\textbf{A}};
    (0,-11)*++{\textbf{B}};
    (-13,15)*{ }="1";
    (-8,15)*{ }="2";
    (3,15)*{ }="V1";
    (-2,15)*{ }="V2";
    (8,15)*{ }="3";
    (13,15)*{ }="4";
    (2.75,7.75)*{}="V2B";
    (5.5,11.25)*{}="V1B";
    (-2.5,0)*{}="2M";
    (2.5,0)*{}="3M";
    (0,4)*{}="2X";
    (-2.5,-15)*{ }="2B";
    (2.5,-15)*{ }="3B";
    "V2";"V2B" **\dir{-};
    "2M";"2B" **\dir{-};
    "3M";"3B" **\dir{-};
    "V1"; "V1B" **\dir{-};
    "V1B"; "3" **\dir{-};
    "V2B"; "2X" **\dir{-};
    "2"; "2X" **\dir{-};
    "1";"2M" **\dir{-};
    "4";"3M" **\dir{-};
    (-15,15)*{}="x1"; 
      (-15,-15)*{}="x2";
      (15,15)*{}="x4";
      (15,-15)*{}="x3";
        "x1";"x2"; **\dir{-};
        "x2";"x3"; **\dir{-};
        "x3";"x4"; **\dir{-};
        "x4";"x1"; **\dir{-}; 
        (-13,18)*{\scriptstyle L};
        (-8, 18)*{\scriptstyle R};
        (-3,18)*{\scriptstyle L};
        (2, 18)*{\scriptstyle R};
        (13,18)*{\scriptstyle R};
        (8,18)*{\scriptstyle L};
        (-2.5,-18)*{\scriptstyle L};
        (2.5,-18)*{\scriptstyle R};
            (-13,14.9)*{};  (-8,14.9)*{} **\dir{.};
            (-12,14)*{};  (-7.5,14)*{} **\dir{.};
            (-11.5,13)*{};  (-7,13)*{} **\dir{.};
            (-11,12)*{};  (-6.5,12)*{} **\dir{.};
            (3,14.9)*{};  (-2,14.9)*{} **\dir{.};
            (3.5,14)*{};  (-1,14)*{} **\dir{.};
            (4,13)*{};  (-.5,13)*{} **\dir{.};
            (4.5,12)*{};  (0,12)*{} **\dir{.};
            (10,11)*{};  (1,11)*{} **\dir{.};
            (9.5,10)*{};  (1.5,10)*{} **\dir{.};
            (8.5,9)*{};  (2,9)*{} **\dir{.};
            (8,8)*{};  (3.5,8)*{} **\dir{.};
            (-7,7)*{};  (-2.5,7)*{} **\dir{.};
            (-6.5,6)*{};  (-1.5,6)*{} **\dir{.};
            (-5.25,5)*{};  (-1,5)*{} **\dir{.};
            (-4.75,4)*{};  (0,4)*{} **\dir{.};
            (-4,3)*{};  (4,3)*{} **\dir{.};
            (-3,2)*{};  (3,2)*{} **\dir{.};
            (13,14.9)*{};  (8,14.9)*{} **\dir{.};
            (12,14)*{};  (7.5,14)*{} **\dir{.};
            (11.5,13)*{};  (7,13)*{} **\dir{.};
            (11,12)*{};  (6.5,12)*{} **\dir{.};
            (-10,11)*{};  (-5.5,11)*{} **\dir{.};
            (-9.5,10)*{};  (-5,10)*{} **\dir{.};
            (-8.5,9)*{};  (-4,9)*{} **\dir{.};
            (-8,8)*{};  (-3.5,8)*{} **\dir{.};
            (7,7)*{};  (2.5,7)*{} **\dir{.};
            (6.5,6)*{};  (1.5,6)*{} **\dir{.};
            (5.25,5)*{};  (1,5)*{} **\dir{.};
            (4.75,4)*{};  (0,4)*{} **\dir{.};
            (-4,3)*{};  (4,3)*{} **\dir{.};
            (-3,2)*{};  (3,2)*{} **\dir{.};
            (-2.75,1)*{};  (2.75,1)*{} **\dir{.};
            (-2.5,0)*{};  (2.5,0)*{} **\dir{.};
            (-2.5,-1)*{};  (2.5,-1)*{} **\dir{.};
            (-2.5,-2)*{};  (2.5,-2)*{} **\dir{.};
            (-2.5,-3)*{};  (2.5,-3)*{} **\dir{.};
            (-2.5,-4)*{};  (2.5,-4)*{} **\dir{.};
            (-2.5,-5)*{};  (2.5,-5)*{} **\dir{.};
            (-2.5,-6)*{};  (2.5,-6)*{} **\dir{.};
            (-2.5,-7)*{};  (2.5,-7)*{} **\dir{.};
            (-2.5,-8)*{};  (2.5,-8)*{} **\dir{.};
            (-2.5,-9)*{};  (2.5,-9)*{} **\dir{.};
            (-2.5,-10)*{};  (2.5,-10)*{} **\dir{.};
            (-2.5,-11)*{};  (2.5,-11)*{} **\dir{.};
            (-2.5,-12)*{};  (2.5,-12)*{} **\dir{.};
            (-2.5,-13)*{};  (2.5,-13)*{} **\dir{.};
            (-2.5,-14)*{};  (2.5,-14)*{} **\dir{.};
            (-2.5,-14.9)*{};  (2.5,-14.9)*{} **\dir{.};
\endxy
 \]
follows from the interchange law in the 2-category $\mathcal{D}$,
and the unit laws:
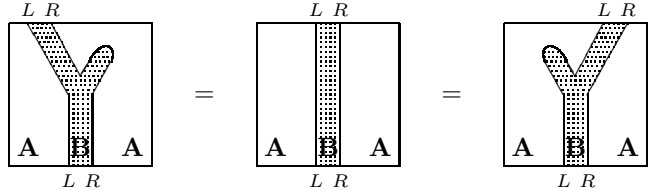
\[
\xy 0;/r.15pc/:
    (-11,-11)*++{\textbf{A}};
    (11,-11)*++{\textbf{A}};
    (0,-11)*++{\textbf{B}};
    (-6,15)*{ }="3";
    (-11,15)*{ }="4";
    (2.5,0)*{}="2M";
    (-2.5,0)*{}="3M";
    (0,4)*{}="2X";
    (2.5,-15)*{ }="2B";
    (-2.5,-15)*{ }="3B";
         (6.5,7)*{}="A";
        (1.85,7)*{}="B";
        "A";"B" **\crv{(8,11) & (3.75, 11.5)};
    "2M";"2B" **\dir{-};
    "3M";"3B" **\dir{-};
    "B"; "2X" **\dir{-};
    "3"; "2X" **\dir{-};
    "A";"2M" **\dir{-};
    "4";"3M" **\dir{-};
    (-15,15)*{}="x1"; 
      (-15,-15)*{}="x2";
      (15,15)*{}="x4";
      (15,-15)*{}="x3";
        "x1";"x2"; **\dir{-};
        "x2";"x3"; **\dir{-};
        "x3";"x4"; **\dir{-};
        "x4";"x1"; **\dir{-}; 
        (-11,18)*{\scriptstyle L};
        (-6,18)*{\scriptstyle R};
        (-2.5,-18)*{\scriptstyle L};
        (2.5,-18)*{\scriptstyle R};
            (6,10)*{};  (5,10)*{} **\dir{.};
            (6.5,9)*{};  (3.75,9)*{} **\dir{.};
            (6.5,8)*{};  (3,8)*{} **\dir{.};
            (6,7)*{};  (2.25,7)*{} **\dir{.};
            (-5.5,6)*{};  (-1.5,6)*{} **\dir{.};
            (-5,5)*{};  (-1,5)*{} **\dir{.};
            (-4.5,4)*{};  (0,4)*{} **\dir{.};
            (-4,3)*{};  (4,3)*{} **\dir{.};
            (-3,2)*{};  (3,2)*{} **\dir{.};
            (-11,14.9)*{};  (-6,14.9)*{} **\dir{.};
            (-10,14)*{};  (-5.5,14)*{} **\dir{.};
            (-9.5,13)*{};  (-5,13)*{} **\dir{.};
            (-9,12)*{};  (-4.5,12)*{} **\dir{.};
            (-8.5,11)*{};  (-4,11)*{} **\dir{.};
            (-8,10)*{};  (-3.5,10)*{} **\dir{.};
            (-7.5,9)*{};  (-3,9)*{} **\dir{.};
            (-7,8)*{};  (-2.5,8)*{} **\dir{.};
            (-6,7)*{};  (-2.25,7)*{} **\dir{.};
            (6,7)*{};  (2.25,7)*{} **\dir{.};
            (5.5,6)*{};  (1.5,6)*{} **\dir{.};
            (5,5)*{};  (1,5)*{} **\dir{.};
            (4.5,4)*{};  (0,4)*{} **\dir{.};
            (-2.75,1)*{};  (2.75,1)*{} **\dir{.};
            (-2.5,0)*{};  (2.5,0)*{} **\dir{.};
            (-2.5,-1)*{};  (2.5,-1)*{} **\dir{.};
            (-2.5,-2)*{};  (2.5,-2)*{} **\dir{.};
            (-2.5,-3)*{};  (2.5,-3)*{} **\dir{.};
            (-2.5,-4)*{};  (2.5,-4)*{} **\dir{.};
            (-2.5,-5)*{};  (2.5,-5)*{} **\dir{.};
            (-2.5,-6)*{};  (2.5,-6)*{} **\dir{.};
            (-2.5,-7)*{};  (2.5,-7)*{} **\dir{.};
            (-2.5,-8)*{};  (2.5,-8)*{} **\dir{.};
            (-2.5,-9)*{};  (2.5,-9)*{} **\dir{.};
            (-2.5,-10)*{};  (2.5,-10)*{} **\dir{.};
            (-2.5,-11)*{};  (2.5,-11)*{} **\dir{.};
            (-2.5,-12)*{};  (2.5,-12)*{} **\dir{.};
            (-2.5,-13)*{};  (2.5,-13)*{} **\dir{.};
            (-2.5,-14)*{};  (2.5,-14)*{} **\dir{.};
            (-2.5,-14.9)*{};  (2.5,-14.9)*{} **\dir{.};
\endxy
\quad = \quad
\xy 0;/r.15pc/:
    (-11,-11)*++{\textbf{A}};
    (11,-11)*++{\textbf{A}};
    (0,-11)*++{\textbf{B}};
    (-2.5,15)*{ }="2";
    (2.5,15)*{ }="3";
    (-2.5,-15)*{ }="2B";
    (2.5,-15)*{ }="3B";
    "2";"2B" **\dir{-};
    "3";"3B" **\dir{-};
    (-15,15)*{}="x1"; 
      (-15,-15)*{}="x2";
      (15,15)*{}="x4";
      (15,-15)*{}="x3";
        "x1";"x2"; **\dir{-};
        "x2";"x3"; **\dir{-};
        "x3";"x4"; **\dir{-};
        "x4";"x1"; **\dir{-}; 
        (-2.5,18)*{\scriptstyle L};
        (2.5, 18)*{\scriptstyle R};
        (-2.5,-18)*{\scriptstyle L};
        (2.5,-18)*{\scriptstyle R};
            (-2.5,2)*{};  (2.5,2)*{} **\dir{.};
            (-2.5,3)*{};  (2.5,3)*{} **\dir{.};
            (-2.5,4)*{};  (2.5,4)*{} **\dir{.};
            (-2.5,5)*{};  (2.5,5)*{} **\dir{.};
            (-2.5,6)*{};  (2.5,6)*{} **\dir{.};
            (-2.5,7)*{};  (2.5,7)*{} **\dir{.};
            (-2.5,8)*{};  (2.5,8)*{} **\dir{.};
            (-2.5,9)*{};  (2.5,9)*{} **\dir{.};
            (-2.5,10)*{};  (2.5,10)*{} **\dir{.};
            (-2.5,11)*{};  (2.5,11)*{} **\dir{.};
            (-2.5,12)*{};  (2.5,12)*{} **\dir{.};
            (-2.5,13)*{};  (2.5,13)*{} **\dir{.};
            (-2.5,14)*{};  (2.5,14)*{} **\dir{.};
            (-2.5,14.9)*{};  (2.5,14.9)*{} **\dir{.};
            (-2.5,1)*{};  (2.5,1)*{} **\dir{.};
            (-2.5,0)*{};  (2.5,0)*{} **\dir{.};
            (-2.5,-1)*{};  (2.5,-1)*{} **\dir{.};
            (-2.5,-2)*{};  (2.5,-2)*{} **\dir{.};
            (-2.5,-3)*{};  (2.5,-3)*{} **\dir{.};
            (-2.5,-4)*{};  (2.5,-4)*{} **\dir{.};
            (-2.5,-5)*{};  (2.5,-5)*{} **\dir{.};
            (-2.5,-6)*{};  (2.5,-6)*{} **\dir{.};
            (-2.5,-7)*{};  (2.5,-7)*{} **\dir{.};
            (-2.5,-8)*{};  (2.5,-8)*{} **\dir{.};
            (-2.5,-9)*{};  (2.5,-9)*{} **\dir{.};
            (-2.5,-10)*{};  (2.5,-10)*{} **\dir{.};
            (-2.5,-11)*{};  (2.5,-11)*{} **\dir{.};
            (-2.5,-12)*{};  (2.5,-12)*{} **\dir{.};
            (-2.5,-13)*{};  (2.5,-13)*{} **\dir{.};
            (-2.5,-14)*{};  (2.5,-14)*{} **\dir{.};
            (-2.5,-14.9)*{};  (2.5,-14.9)*{} **\dir{.};
\endxy
\quad = \quad
\xy 0;/r.15pc/: 
    (-11,-11)*++{\textbf{A}};
    (11,-11)*++{\textbf{A}};
    (0,-11)*++{\textbf{B}};
    (6,15)*{ }="3";
    (11,15)*{ }="4";
    (-2.5,0)*{}="2M";
    (2.5,0)*{}="3M";
    (0,4)*{}="2X";
    (-2.5,-15)*{ }="2B";
    (2.5,-15)*{ }="3B";
         (-6.5,7)*{}="A";
        (-1.85,7)*{}="B";
        "A";"B" **\crv{(-8,11) & (-3.75, 11.5)};
    "2M";"2B" **\dir{-};
    "3M";"3B" **\dir{-};
    "B"; "2X" **\dir{-};
    "3"; "2X" **\dir{-};
    "A";"2M" **\dir{-};
    "4";"3M" **\dir{-};
    (-15,15)*{}="x1"; 
      (-15,-15)*{}="x2";
      (15,15)*{}="x4";
      (15,-15)*{}="x3";
        "x1";"x2"; **\dir{-};
        "x2";"x3"; **\dir{-};
        "x3";"x4"; **\dir{-};
        "x4";"x1"; **\dir{-}; 
        (11,18)*{\scriptstyle R};
        (6,18)*{\scriptstyle L};
        (-2.5,-18)*{\scriptstyle L};
        (2.5,-18)*{\scriptstyle R};
            (-6,10)*{};  (-5,10)*{} **\dir{.};
            (-6.5,9)*{};  (-3.75,9)*{} **\dir{.};
            (-6.5,8)*{};  (-3,8)*{} **\dir{.};
            (-6,7)*{};  (-2.25,7)*{} **\dir{.};
            (-5.5,6)*{};  (-1.5,6)*{} **\dir{.};
            (-5,5)*{};  (-1,5)*{} **\dir{.};
            (-4.5,4)*{};  (0,4)*{} **\dir{.};
            (-4,3)*{};  (4,3)*{} **\dir{.};
            (-3,2)*{};  (3,2)*{} **\dir{.};
            (11,14.9)*{};  (6,14.9)*{} **\dir{.};
            (10,14)*{};  (5.5,14)*{} **\dir{.};
            (9.5,13)*{};  (5,13)*{} **\dir{.};
            (9,12)*{};  (4.5,12)*{} **\dir{.};
            (8.5,11)*{};  (4,11)*{} **\dir{.};
            (8,10)*{};  (3.5,10)*{} **\dir{.};
            (7.5,9)*{};  (3,9)*{} **\dir{.};
            (7,8)*{};  (2.5,8)*{} **\dir{.};
            (6,7)*{};  (2.25,7)*{} **\dir{.};
            (5.5,6)*{};  (1.5,6)*{} **\dir{.};
            (5,5)*{};  (1,5)*{} **\dir{.};
            (4.5,4)*{};  (0,4)*{} **\dir{.};
            (-2.75,1)*{};  (2.75,1)*{} **\dir{.};
            (-2.5,0)*{};  (2.5,0)*{} **\dir{.};
            (-2.5,-1)*{};  (2.5,-1)*{} **\dir{.};
            (-2.5,-2)*{};  (2.5,-2)*{} **\dir{.};
            (-2.5,-3)*{};  (2.5,-3)*{} **\dir{.};
            (-2.5,-4)*{};  (2.5,-4)*{} **\dir{.};
            (-2.5,-5)*{};  (2.5,-5)*{} **\dir{.};
            (-2.5,-6)*{};  (2.5,-6)*{} **\dir{.};
            (-2.5,-7)*{};  (2.5,-7)*{} **\dir{.};
            (-2.5,-8)*{};  (2.5,-8)*{} **\dir{.};
            (-2.5,-9)*{};  (2.5,-9)*{} **\dir{.};
            (-2.5,-10)*{};  (2.5,-10)*{} **\dir{.};
            (-2.5,-11)*{};  (2.5,-11)*{} **\dir{.};
            (-2.5,-12)*{};  (2.5,-12)*{} **\dir{.};
            (-2.5,-13)*{};  (2.5,-13)*{} **\dir{.};
            (-2.5,-14)*{};  (2.5,-14)*{} **\dir{.};
            (-2.5,-14.9)*{};  (2.5,-14.9)*{} **\dir{.};
\endxy
 \]
follow from the triangle identities in the definition of an
adjunction.

Starting with an adjunction $\radjunction{A}{B}{L}{R}$ where $L$
is the right adjoint produces the color inverted versions of the
diagrams above. The unit $j \maps 1_B \To LR$ and counit $k \maps
RL \To 1_A$ would appear as:
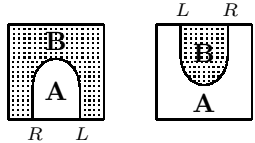
\[ \xy 0;/r.15pc/:  
    (0,6.5)*+{\textbf{B}};
    (0,-4)*+{\textbf{A}};
    (-5,-10)*{ }="TL";
    (5,-10)*{ }="TR";
    (-4.5,-13)*{\scriptstyle R};
    (5.5,-13)*{\scriptstyle L};
    (-5,-4)*{}="l";
    (5,-4)*{}="r";
    "l";"r" **\crv{(-5,5)&(5,5)};
    "TL";"l" **\dir{-};
    "r";"TR" **\dir{-};
       (-10,10)*{}="x1";
       (-10,-10)*{}="x2";
       (10,10)*{}="x4";
       (10,-10)*{}="x3";
        "x1";"x2"; **\dir{-};
        "x2";"x3"; **\dir{-};
        "x3";"x4"; **\dir{-};
        "x4";"x1"; **\dir{-};
        (-10,-9)*{};(-5,-9)    **\dir{.};
        (-10,-8)*{};(-5,-8)    **\dir{.};
        (-10,-7)*{};(-5,-7)    **\dir{.};
        (-10,-6)*{};(-5,-6)    **\dir{.};
        (-10,-5)*{};(-5,-5)    **\dir{.};
        (-10,-4)*{};(-5,-4)    **\dir{.};
        (-10,-3)*{};(-5,-3)    **\dir{.};
        (-10,-2)*{};(-5,-2)    **\dir{.};
        (-10,-1)*{};(-5,-1)    **\dir{.};
        (-10,0)*{};(-5,0)    **\dir{.};
        (-10,1)*{};(-4,1)    **\dir{.};
        (-10,2)*{};(-3,2)    **\dir{.};
        (10,-9)*{};(5,-9)    **\dir{.};
        (10,-8)*{};(5,-8)    **\dir{.};
        (10,-7)*{};(5,-7)    **\dir{.};
        (10,-6)*{};(5,-6)    **\dir{.};
        (10,-5)*{};(5,-5)    **\dir{.};
        (10,-4)*{};(5,-4)    **\dir{.};
        (10,-3)*{};(5,-3)    **\dir{.};
        (10,-2)*{};(5,-2)    **\dir{.};
        (10,-1)*{};(5,-1)    **\dir{.};
        (10,0)*{};(5,0)    **\dir{.};
          (10,1)*{};(4,1)    **\dir{.};
          (10,2)*{};(3,2)    **\dir{.};
            (-10,9)*{};(10,9)    **\dir{.};
            (-10,8)*{};(10,8)    **\dir{.};
            (-10,7)*{};(10,7)    **\dir{.};
            (-10,6)*{};(10,6)    **\dir{.};
            (-10,5)*{};(10,5)    **\dir{.};
            (-10,4)*{};(10,4)    **\dir{.};
            (-10,3)*{};(10,3)    **\dir{.};
        \endxy
\qquad
\xy 0;/r.15pc/:  
    (0,-6.5)*+{\textbf{A}};
    (0,4)*+{\textbf{B}};
    (-5,10)*{ }="TL";
    (5,10)*{ }="TR";
    (-4.5,13)*{\scriptstyle L};
    (5.5,13)*{\scriptstyle R};
    (-5,4)*{}="l";
    (5,4)*{}="r";
    "l";"r" **\crv{(-5,-5)&(5,-5)};
    "TL";"l" **\dir{-};
    "r";"TR" **\dir{-};
       (-10,10)*{}="x1";
       (-10,-10)*{}="x2";
       (10,10)*{}="x4";
       (10,-10)*{}="x3";
        "x1";"x2"; **\dir{-};
        "x2";"x3"; **\dir{-};
        "x3";"x4"; **\dir{-};
        "x4";"x1"; **\dir{-};
        (-5,9)*{};(5,9)    **\dir{.};
        (-5,8)*{};(5,8)    **\dir{.};
        (-5,7)*{};(5,7)    **\dir{.};
        (-5,6)*{};(5,6)    **\dir{.};
        (-5,5)*{};(5,5)    **\dir{.};
        (-5,4)*{};(5,4)    **\dir{.};
        (-5,3)*{};(5,3)    **\dir{.};
        (-4.5,2)*{};(4.5,2)    **\dir{.};
        (-4,1)*{};(4,1)    **\dir{.};
        (-4,0)*{};(4,0)    **\dir{.};
        (-3,-1)*{};(3,-1)    **\dir{.};
        (-2,-2)*{};(2,-2)    **\dir{.};
        \endxy
\]
In this case $RL$ becomes a comonoid in $\Hom(A,A)$ whose
comultiplication is given by the diagram:
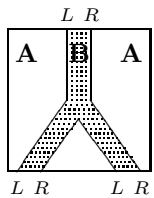
\[
\xy 0;/r.15pc/: 
    (-11,10)*++{\textbf{A}};
    (11,10)*++{\textbf{A}};
    (0,10)*++{\textbf{B}};
    (-13,-15)*{ }="1";
    (-8,-15)*{ }="2";
    (8,-15)*{ }="3";
    (13,-15)*{ }="4";
    (-2.5,0)*{}="2M";
    (2.5,0)*{}="3M";
    (0,-4)*{}="2X";
    (-2.5,15)*{ }="2B";
    (2.5,15)*{ }="3B";
    "2M";"2B" **\dir{-};
    "3M";"3B" **\dir{-};
    "2X"; "2" **\dir{-};
    "3"; "2X" **\dir{-};
    "1";"2M" **\dir{-};
    "4";"3M" **\dir{-};
    (-15,-15)*{}="x1"; 
      (-15,15)*{}="x2";
      (15,-15)*{}="x4";
      (15,15)*{}="x3";
        "x1";"x2"; **\dir{-};
        "x2";"x3"; **\dir{-};
        "x3";"x4"; **\dir{-};
        "x4";"x1"; **\dir{-}; 
        (-13,-18.5)*{\scriptstyle L};
        (-8, -18.5)*{\scriptstyle R};
        (13,-18.5)*{\scriptstyle R};
        (8,-18.5)*{\scriptstyle L};
        (-2.5,18)*{\scriptstyle L};
        (2.5,18)*{\scriptstyle R};
            (-13,-14.9)*{};  (-8,-14.9)*{} **\dir{.};
            (-12,-14)*{};  (-7.5,-14)*{} **\dir{.};
            (-11.5,-13)*{};  (-7,-13)*{} **\dir{.};
            (-11,-12)*{};  (-6.5,-12)*{} **\dir{.};
            (-10,-11)*{};  (-5.5,-11)*{} **\dir{.};
            (-9.5,-10)*{};  (-5,-10)*{} **\dir{.};
            (-8.5,-9)*{};  (-4,-9)*{} **\dir{.};
            (-8,-8)*{};  (-3.5,-8)*{} **\dir{.};
            (-7,-7)*{};  (-2.5,-7)*{} **\dir{.};
            (-6.5,-6)*{};  (-1.5,-6)*{} **\dir{.};
            (-5.25,-5)*{};  (-1,-5)*{} **\dir{.};
            (-4.75,-4)*{};  (0,-4)*{} **\dir{.};
            (-4,-3)*{};  (4,-3)*{} **\dir{.};
            (-3,-2)*{};  (3,-2)*{} **\dir{.};
            (13,-14.9)*{};  (8,-14.9)*{} **\dir{.};
            (12,-14)*{};  (7.5,-14)*{} **\dir{.};
            (11.5,-13)*{};  (7,-13)*{} **\dir{.};
            (11,-12)*{};  (6.5,-12)*{} **\dir{.};
            (10,-11)*{};  (5.5,-11)*{} **\dir{.};
            (9.5,-10)*{};  (5,-10)*{} **\dir{.};
            (8.5,-9)*{};  (4,-9)*{} **\dir{.};
            (8,-8)*{};  (3.5,-8)*{} **\dir{.};
            (7,-7)*{};  (2.5,-7)*{} **\dir{.};
            (6.5,-6)*{};  (1.5,-6)*{} **\dir{.};
            (5.25,-5)*{};  (1,-5)*{} **\dir{.};
            (4.75,-4)*{};  (0,-4)*{} **\dir{.};
            (-4,-3)*{};  (4,-3)*{} **\dir{.};
            (-3,-2)*{};  (3,-2)*{} **\dir{.};
            (-2.75,-1)*{};  (2.75,-1)*{} **\dir{.};
            (-2.5,0)*{};  (2.5,0)*{} **\dir{.};
            (-2.5,1)*{};  (2.5,1)*{} **\dir{.};
            (-2.5,2)*{};  (2.5,2)*{} **\dir{.};
            (-2.5,3)*{};  (2.5,3)*{} **\dir{.};
            (-2.5,4)*{};  (2.5,4)*{} **\dir{.};
            (-2.5,5)*{};  (2.5,5)*{} **\dir{.};
            (-2.5,6)*{};  (2.5,6)*{} **\dir{.};
            (-2.5,7)*{};  (2.5,7)*{} **\dir{.};
            (-2.5,8)*{};  (2.5,8)*{} **\dir{.};
            (-2.5,9)*{};  (2.5,9)*{} **\dir{.};
            (-2.5,10)*{};  (2.5,10)*{} **\dir{.};
            (-2.5,11)*{};  (2.5,11)*{} **\dir{.};
            (-2.5,12)*{};  (2.5,12)*{} **\dir{.};
            (-2.5,13)*{};  (2.5,13)*{} **\dir{.};
            (-2.5,14)*{};  (2.5,14)*{} **\dir{.};
            (-2.5,14.9)*{};  (2.5,14.9)*{} **\dir{.};
\endxy
 \]
 and whose counit is:
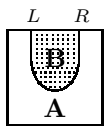
\[ \xy 0;/r.15pc/:  
    (0,-6.5)*+{\textbf{A}};
    (0,4)*+{\textbf{B}};
    (-5,10)*{ }="TL";
    (5,10)*{ }="TR";
    (-4.5,13)*{\scriptstyle L};
    (5.5,13)*{\scriptstyle R};
    (-5,4)*{}="l";
    (5,4)*{}="r";
    "l";"r" **\crv{(-5,-5)&(5,-5)};
    "TL";"l" **\dir{-};
    "r";"TR" **\dir{-};
       (-10,10)*{}="x1";
       (-10,-10)*{}="x2";
       (10,10)*{}="x4";
       (10,-10)*{}="x3";
        "x1";"x2"; **\dir{-};
        "x2";"x3"; **\dir{-};
        "x3";"x4"; **\dir{-};
        "x4";"x1"; **\dir{-};
        (-5,9)*{};(5,9)    **\dir{.};
        (-5,8)*{};(5,8)    **\dir{.};
        (-5,7)*{};(5,7)    **\dir{.};
        (-5,6)*{};(5,6)    **\dir{.};
        (-5,5)*{};(5,5)    **\dir{.};
        (-5,4)*{};(5,4)    **\dir{.};
        (-5,3)*{};(5,3)    **\dir{.};
        (-4.5,2)*{};(4.5,2)    **\dir{.};
        (-4,1)*{};(4,1)    **\dir{.};
        (-4,0)*{};(4,0)    **\dir{.};
        (-3,-1)*{};(3,-1)    **\dir{.};
        (-2,-2)*{};(2,-2)    **\dir{.};
        \endxy
\]
the counit of the right adjunction. By similar diagrams as those
above, the coassociativity and counit axioms follow from the
axioms of a 2-category and the axioms of an adjunction.

When the morphisms $L$ is both left and right adjoint to $R$ the
object $RL$ of $\Hom(A,A)$ is both a monoid and a comonoid. These
structures satisfy compatibility conditions, known as the
Frobenius identities, making $RL$ into a Frobenius object. Indeed,
the Frobenius identities:
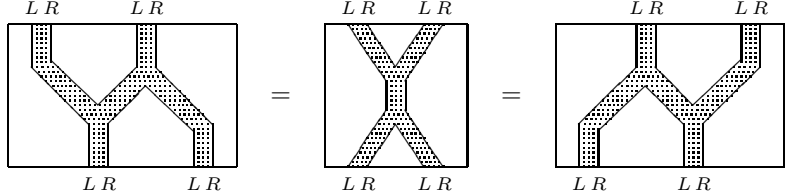
\[
\xy  0;/r.15pc/: 
    ( 3 , 15)*{ }="TR";
    ( 7, 15)*{ }="TRR";
    ( 3 , 6)*{}="TR2";
    ( 7, 6)*{}="TRR2";
    (-19, 15)*{ }="TLL";
    (-15, 15)*{ }="TL";
    (-19, 6)*{}="TLL2";
    (-15, 7.5)*{}="TL2";
    ( -3,-15)*{ }="BL";
    (-7,-15)*{ }="BLL";
    ( -3,-6)*{}="BL2";
    (-7,-6)*{}="BLL2";
    ( 15,-15)*{ }="BR";
    ( 19,-15)*{ }="BRR";
    ( 15,-7.5)*{}="BR2";
    ( 19,-6)*{}="BRR2";
    ( -5, -2)*{}="M1";
    ( 5 , 2)*{}="M2";
        "TR" ;"TR2" **\dir{-};
        "TRR";"TRR2" **\dir{-};
        "TL" ;"TL2" **\dir{-};
        "TLL";"TLL2" **\dir{-};
        "BR" ;"BR2" **\dir{-};
        "BRR";"BRR2" **\dir{-};
        "BL" ;"BL2" **\dir{-};
        "BLL";"BLL2" **\dir{-};
        "TL2" ;"M1" **\dir{-};
        "TR2" ;"M1" **\dir{-};
        "BL2" ;"M2" **\dir{-};
        "BR2" ;"M2" **\dir{-};
        "TLL2" ;"BLL2" **\dir{-};
        "TRR2" ;"BRR2" **\dir{-};
    (-24,-15)*{}="x1"; 
      (-24,15)*{}="x2";
      (24,-15)*{}="x4";
      (24,15)*{}="x3";
        "x1";"x2"; **\dir{-};
        "x2";"x3"; **\dir{-};
        "x3";"x4"; **\dir{-};
        "x4";"x1"; **\dir{-}; 
        (3,18.5)*{\scriptstyle L};
        (7, 18.5)*{\scriptstyle R};
        (-19,18.5)*{\scriptstyle L};
        (-15, 18.5)*{\scriptstyle R};
        (-7,-18.5)*{\scriptstyle L};
        (-3,-18.5)*{\scriptstyle R};
        (19,-18.5)*{\scriptstyle R};
        (15, -18.5)*{\scriptstyle L};
            (-19,14.9)*{}; (-15,14.9)*{} **\dir{.};
            (-19,14)*{}; (-15,14)*{} **\dir{.};
            (-19,13)*{}; (-15,13)*{} **\dir{.};
            (-19,12)*{}; (-15,12)*{} **\dir{.};
            (-19,11)*{}; (-15,11)*{} **\dir{.};
            (-19,10)*{}; (-15,10)*{} **\dir{.};
            (-19,9)*{}; (-15,9)*{} **\dir{.};
            (-19,8)*{}; (-15,8)*{} **\dir{.};
            (-19,7)*{}; (-15,7)*{} **\dir{.};
            (-19,6)*{}; (-13.5,6)*{} **\dir{.};
            (-17.5,5)*{}; (-12.5,5)*{} **\dir{.};
            (-17,4)*{}; (-12,4)*{} **\dir{.};
            (-16,3)*{}; (-11,3)*{} **\dir{.};
            (-14.5,2)*{}; (-10,2)*{} **\dir{.};
            (-14,1)*{}; (-9,1)*{} **\dir{.};
            (-13,0)*{}; (-8,0)*{} **\dir{.};
            (-12,-1)*{}; (-7,-1)*{} **\dir{.};
            (3,14.9)*{}; (7,14.9)*{} **\dir{.};
            (3,14)*{}; (7,14)*{} **\dir{.};
            (3,13)*{}; (7,13)*{} **\dir{.};
            (3,12)*{}; (7,12)*{} **\dir{.};
            (3,11)*{}; (7,11)*{} **\dir{.};
            (3,10)*{}; (7,10)*{} **\dir{.};
            (3,9)*{}; (7,9)*{} **\dir{.};
            (3,8)*{}; (7,8)*{} **\dir{.};
            (3,7)*{}; (7,7)*{} **\dir{.};
            (3,6)*{}; (7,6)*{} **\dir{.};
            (2,5)*{}; (8,5)*{} **\dir{.};
            (1.5,4)*{}; (8.5,4)*{} **\dir{.};
            (.5,3)*{}; (9,3)*{} **\dir{.};
            (-.5,2)*{}; (10,2)*{} **\dir{.};
            (-3,-14.9)*{}; (-7,-14.9)*{} **\dir{.};
            (-3,-14)*{}; (-7,-14)*{} **\dir{.};
            (-3,-13)*{}; (-7,-13)*{} **\dir{.};
            (-3,-12)*{}; (-7,-12)*{} **\dir{.};
            (-3,-11)*{}; (-7,-11)*{} **\dir{.};
            (-3,-10)*{}; (-7,-10)*{} **\dir{.};
            (-3,-9)*{}; (-7,-9)*{} **\dir{.};
            (-3,-8)*{}; (-7,-8)*{} **\dir{.};
            (-3,-7)*{}; (-7,-7)*{} **\dir{.};
            (-3,-6)*{}; (-7,-6)*{} **\dir{.};
            (-2,-5)*{}; (-8,-5)*{} **\dir{.};
            (-1.5,-4)*{}; (-8.5,-4)*{} **\dir{.};
            (-.5,-3)*{}; (-9,-3)*{} **\dir{.};
            (.5,-2)*{}; (-10,-2)*{} **\dir{.};
            (19,-14.9)*{}; (15,-14.9)*{} **\dir{.};
            (19,-14)*{}; (15,-14)*{} **\dir{.};
            (19,-13)*{}; (15,-13)*{} **\dir{.};
            (19,-12)*{}; (15,-12)*{} **\dir{.};
            (19,-11)*{}; (15,-11)*{} **\dir{.};
            (19,-10)*{}; (15,-10)*{} **\dir{.};
            (19,-9)*{}; (15,-9)*{} **\dir{.};
            (19,-8)*{}; (15,-8)*{} **\dir{.};
            (19,-7)*{}; (15,-7)*{} **\dir{.};
            (19,-6)*{}; (13.5,-6)*{} **\dir{.};
            (17.5,-5)*{}; (12.5,-5)*{} **\dir{.};
            (17,-4)*{}; (12,-4)*{} **\dir{.};
            (16,-3)*{}; (11,-3)*{} **\dir{.};
            (14.5,-2)*{}; (10,-2)*{} **\dir{.};
            (14,-1)*{}; (9,-1)*{} **\dir{.};
            (13,0)*{}; (8,0)*{} **\dir{.};
            (12,1)*{}; (7,1)*{} **\dir{.};
            (-1.5,1)*{}; (3.5,1)*{} **\dir{.};
            (-2.5,0)*{}; (2.5,0)*{} **\dir{.};
            (-3.5,-1)*{}; (1.5,-1)*{} **\dir{.};
\endxy
\quad = \quad
\xy  0;/r.15pc/: 
    ( 6 , 15)*{ }="TR";
    ( 10, 15)*{ }="TRR";
    ( -2 , 3)*{}="TL2";
    ( 2, 3)*{}="TR2";
    ( -2 , -3)*{}="BL2";
    ( 2, -3)*{}="BR2";
    (-10, 15)*{ }="TLL";
    (-6, 15)*{ }="TL";
    ( -6,-15)*{ }="BL";
    (-10,-15)*{ }="BLL";
    ( 6,-15)*{ }="BR";
    ( 10,-15)*{ }="BRR";
    ( 0, 6)*{}="M1";
    ( 0 , -6)*{}="M2";
        "TRR" ;"TR2" **\dir{-};
        "TLL" ;"TL2" **\dir{-};
        "BRR" ;"BR2" **\dir{-};
        "BLL" ;"BL2" **\dir{-};
        "TL2" ;"BL2" **\dir{-};
        "TR2" ;"BR2" **\dir{-};
        "TL" ;"M1" **\dir{-};
        "TR" ;"M1" **\dir{-};
        "BL" ;"M2" **\dir{-};
        "BR" ;"M2" **\dir{-};
    (-15,-15)*{}="x1"; 
      (-15,15)*{}="x2";
      (15,-15)*{}="x4";
      (15,15)*{}="x3";
        "x1";"x2"; **\dir{-};
        "x2";"x3"; **\dir{-};
        "x3";"x4"; **\dir{-};
        "x4";"x1"; **\dir{-}; 
        (-10,18.5)*{\scriptstyle L};
        (-6, 18.5)*{\scriptstyle R};
        (10,18.5)*{\scriptstyle R};
        (6, 18.5)*{\scriptstyle L};
        (-10,-18.5)*{\scriptstyle L};
        (-6, -18.5)*{\scriptstyle R};
        (10,-18.5)*{\scriptstyle R};
        (6, -18.5)*{\scriptstyle L};
            (-2,3)*{}; (2,3)*{} **\dir{.};
            (-2,2)*{}; (2,2)*{} **\dir{.};
            (-2,1)*{}; (2,1)*{} **\dir{.};
            (-2,0)*{}; (2,0)*{} **\dir{.};
            (-2,-1)*{}; (2,-1)*{} **\dir{.};
            (-2,-2)*{}; (2,-2)*{} **\dir{.};
            (-2,-3)*{}; (2,-3)*{} **\dir{.};
            (-2.5,4)*{}; (2.5,4)*{} **\dir{.};
            (-3,5)*{}; (3,5)*{} **\dir{.};
            (-3.5,6)*{}; (3.5,6)*{} **\dir{.};
            (-2.5,-4)*{}; (2.5,-4)*{} **\dir{.};
            (-3,-5)*{}; (3,-5)*{} **\dir{.};
            (-3.5,-6)*{}; (3.5,-6)*{} **\dir{.};
            (-10,14.9)*{}; (-6,14.9)*{} **\dir{.};
            (-9,14)*{}; (-5.5,14)*{} **\dir{.};
            (-8,13)*{}; (-5,13)*{} **\dir{.};
            (-7.5,12)*{}; (-4.5,12)*{} **\dir{.};
            (-7,11)*{}; (-3.5,11)*{} **\dir{.};
            (-6.5,10)*{}; (-3,10)*{} **\dir{.};
            (-6,9)*{}; (-2.5,9)*{} **\dir{.};
            (-5,8)*{}; (-1.5,8)*{} **\dir{.};
            (-4.5,7)*{}; (-1,7)*{} **\dir{.};
            (10,14.9)*{}; (6,14.9)*{} **\dir{.};
            (9,14)*{}; (5.5,14)*{} **\dir{.};
            (8,13)*{}; (5,13)*{} **\dir{.};
            (7.5,12)*{}; (4.5,12)*{} **\dir{.};
            (7,11)*{}; (3.5,11)*{} **\dir{.};
            (6.5,10)*{}; (3,10)*{} **\dir{.};
            (6,9)*{}; (2.5,9)*{} **\dir{.};
            (5,8)*{}; (1.5,8)*{} **\dir{.};
            (4.5,7)*{}; (1,7)*{} **\dir{.};
            (-10,-14.9)*{}; (-6,-14.9)*{} **\dir{.};
            (-9,-14)*{}; (-5.5,-14)*{} **\dir{.};
            (-8,-13)*{}; (-5,-13)*{} **\dir{.};
            (-7.5,-12)*{}; (-4.5,-12)*{} **\dir{.};
            (-7,-11)*{}; (-3.5,-11)*{} **\dir{.};
            (-6.5,-10)*{}; (-3,-10)*{} **\dir{.};
            (-6,-9)*{}; (-2.5,-9)*{} **\dir{.};
            (-5,-8)*{}; (-1.5,-8)*{} **\dir{.};
            (-4.5,-7)*{}; (-1,-7)*{} **\dir{.};
            (10,-14.9)*{}; (6,-14.9)*{} **\dir{.};
            (9,-14)*{}; (5.5,-14)*{} **\dir{.};
            (8,-13)*{}; (5,-13)*{} **\dir{.};
            (7.5,-12)*{}; (4.5,-12)*{} **\dir{.};
            (7,-11)*{}; (3.5,-11)*{} **\dir{.};
            (6.5,-10)*{}; (3,-10)*{} **\dir{.};
            (6,-9)*{}; (2.5,-9)*{} **\dir{.};
            (5,-8)*{}; (1.5,-8)*{} **\dir{.};
            (4.5,-7)*{}; (1,-7)*{} **\dir{.};
\endxy
\quad = \quad
\xy  0;/r.15pc/: 
    ( -3 , 15)*{ }="TR";
    ( -7, 15)*{ }="TRR";
    ( -3 , 6)*{}="TR2";
    ( -7, 6)*{}="TRR2";
    (19, 15)*{ }="TLL";
    (15, 15)*{ }="TL";
    (19, 6)*{}="TLL2";
    (15, 7.5)*{}="TL2";
    ( 3,-15)*{ }="BL";
    (7,-15)*{ }="BLL";
    ( 3,-6)*{}="BL2";
    (7,-6)*{}="BLL2";
    ( -15,-15)*{ }="BR";
    ( -19,-15)*{ }="BRR";
    ( -15,-7.5)*{}="BR2";
    ( -19,-6)*{}="BRR2";
    ( 5, -2)*{}="M1";
    ( -5 , 2)*{}="M2";
        "TR" ;"TR2" **\dir{-};
        "TRR";"TRR2" **\dir{-};
        "TL" ;"TL2" **\dir{-};
        "TLL";"TLL2" **\dir{-};
        "BR" ;"BR2" **\dir{-};
        "BRR";"BRR2" **\dir{-};
        "BL" ;"BL2" **\dir{-};
        "BLL";"BLL2" **\dir{-};
        "TL2" ;"M1" **\dir{-};
        "TR2" ;"M1" **\dir{-};
        "BL2" ;"M2" **\dir{-};
        "BR2" ;"M2" **\dir{-};
        "TLL2" ;"BLL2" **\dir{-};
        "TRR2" ;"BRR2" **\dir{-};
    (-24,-15)*{}="x1"; 
      (-24,15)*{}="x2";
      (24,-15)*{}="x4";
      (24,15)*{}="x3";
        "x1";"x2"; **\dir{-};
        "x2";"x3"; **\dir{-};
        "x3";"x4"; **\dir{-};
        "x4";"x1"; **\dir{-}; 
        (-3,18.5)*{\scriptstyle R};
        (-7, 18.5)*{\scriptstyle L};
        (19,18.5)*{\scriptstyle R};
        (15, 18.5)*{\scriptstyle L};
        (7,-18.5)*{\scriptstyle R};
        (3,-18.5)*{\scriptstyle L};
        (-19,-18.5)*{\scriptstyle L};
        (-15, -18.5)*{\scriptstyle R};
            (19,14.9)*{}; (15,14.9)*{} **\dir{.};
            (19,14)*{}; (15,14)*{} **\dir{.};
            (19,13)*{}; (15,13)*{} **\dir{.};
            (19,12)*{}; (15,12)*{} **\dir{.};
            (19,11)*{}; (15,11)*{} **\dir{.};
            (19,10)*{}; (15,10)*{} **\dir{.};
            (19,9)*{}; (15,9)*{} **\dir{.};
            (19,8)*{}; (15,8)*{} **\dir{.};
            (19,7)*{}; (15,7)*{} **\dir{.};
            (19,6)*{}; (13.5,6)*{} **\dir{.};
            (17.5,5)*{}; (12.5,5)*{} **\dir{.};
            (17,4)*{}; (12,4)*{} **\dir{.};
            (16,3)*{}; (11,3)*{} **\dir{.};
            (14.5,2)*{}; (10,2)*{} **\dir{.};
            (14,1)*{}; (9,1)*{} **\dir{.};
            (13,0)*{}; (8,0)*{} **\dir{.};
            (12,-1)*{}; (7,-1)*{} **\dir{.};
            (-3,14.9)*{}; (-7,14.9)*{} **\dir{.};
            (-3,14)*{}; (-7,14)*{} **\dir{.};
            (-3,13)*{}; (-7,13)*{} **\dir{.};
            (-3,12)*{}; (-7,12)*{} **\dir{.};
            (-3,11)*{}; (-7,11)*{} **\dir{.};
            (-3,10)*{}; (-7,10)*{} **\dir{.};
            (-3,9)*{}; (-7,9)*{} **\dir{.};
            (-3,8)*{}; (-7,8)*{} **\dir{.};
            (-3,7)*{}; (-7,7)*{} **\dir{.};
            (-3,6)*{}; (-7,6)*{} **\dir{.};
            (-2,5)*{}; (-8,5)*{} **\dir{.};
            (-1.5,4)*{}; (-8.5,4)*{} **\dir{.};
            (-.5,3)*{}; (-9,3)*{} **\dir{.};
            (.5,2)*{}; (-10,2)*{} **\dir{.};
            (3,-14.9)*{}; (7,-14.9)*{} **\dir{.};
            (3,-14)*{}; (7,-14)*{} **\dir{.};
            (3,-13)*{}; (7,-13)*{} **\dir{.};
            (3,-12)*{}; (7,-12)*{} **\dir{.};
            (3,-11)*{}; (7,-11)*{} **\dir{.};
            (3,-10)*{}; (7,-10)*{} **\dir{.};
            (3,-9)*{}; (7,-9)*{} **\dir{.};
            (3,-8)*{}; (7,-8)*{} **\dir{.};
            (3,-7)*{}; (7,-7)*{} **\dir{.};
            (3,-6)*{}; (7,-6)*{} **\dir{.};
            (2,-5)*{}; (8,-5)*{} **\dir{.};
            (1.5,-4)*{}; (8.5,-4)*{} **\dir{.};
            (.5,-3)*{}; (9,-3)*{} **\dir{.};
            (-.5,-2)*{}; (10,-2)*{} **\dir{.};
            (-19,-14.9)*{}; (-15,-14.9)*{} **\dir{.};
            (-19,-14)*{}; (-15,-14)*{} **\dir{.};
            (-19,-13)*{}; (-15,-13)*{} **\dir{.};
            (-19,-12)*{}; (-15,-12)*{} **\dir{.};
            (-19,-11)*{}; (-15,-11)*{} **\dir{.};
            (-19,-10)*{}; (-15,-10)*{} **\dir{.};
            (-19,-9)*{}; (-15,-9)*{} **\dir{.};
            (-19,-8)*{}; (-15,-8)*{} **\dir{.};
            (-19,-7)*{}; (-15,-7)*{} **\dir{.};
            (-19,-6)*{}; (-13.5,-6)*{} **\dir{.};
            (-17.5,-5)*{}; (-12.5,-5)*{} **\dir{.};
            (-17,-4)*{}; (-12,-4)*{} **\dir{.};
            (-16,-3)*{}; (-11,-3)*{} **\dir{.};
            (-14.5,-2)*{}; (-10,-2)*{} **\dir{.};
            (-14,-1)*{}; (-9,-1)*{} **\dir{.};
            (-13,0)*{}; (-8,0)*{} **\dir{.};
            (-12,1)*{}; (-7,1)*{} **\dir{.};
            (1.5,1)*{}; (-3.5,1)*{} **\dir{.};
            (2.5,0)*{}; (-2.5,0)*{} **\dir{.};
            (3.5,-1)*{}; (-1.5,-1)*{} **\dir{.};
\endxy
\]
follow from the interchange axiom of the 2-category $\mathcal{D}$.
Thus we have shown that the axioms of an ambijunction in a
2-category beautifully imply all the axioms of a Frobenius object;
by drawing string diagrams their relationship becomes much more
transparent.

The converse, that every Frobenius object in the monoidal category
$M$ actually arises in this way from an ambidextrous adjunction in
some 2-category $\mathcal{D}$ into which $M$ fully and faithfully
embeds, has thus far not been proven in a completely general
context.  An attempt to prove this result was made by
M\"{u}ger~\cite{Muger} who showed that, with certain extra
assumptions about the monoidal category $M$, a 2-category
$\mathcal{E}$ into which $M$ fully and faithfully embeds can be
constructed, and Frobenius objects in $M$ correspond precisely to
ambijunctions in this 2-category $\mathcal{E}$.  However, we will
see that the converse can be proven quite naturally using the
language of monad theory where questions like this have already
been resolved. Extending the work of Lawvere~\cite{law},
Street~\cite{StreetFrob} has made substantial progress by proving
this result in the context of the 2-category \cat{Cat}. Street's
approach suggests a natural framework for proving the completely
general result and experts in 2-categorical monad theory and
enriched category theory will find that our proof is a straight
forward extension of Street's work. Using a wealth of results from
category theory, especially the formal theory of
monads~\cite{LS1}, we extend Street's work and prove this result
for Frobenius objects in a generic monoidal category $M$.

To understand the relevance of monad theory a bit of background is
in order.  A monad on a category $\mathcal{A}$ can be defined as a
monoid in the functor category $[\mathcal{A},\mathcal{A}]$. The
theory of monads has been well developed and, in particular, it is
well known that every monad $\T$ on a category $\mathcal{A}$
arises from a pair of adjoint functors
$\adjunction{\mathcal{A}}{\mathcal{B}}{}{}$. The problem of
constructing an adjunction from a monad has two well known
solutions --- the Kleisli construction
$\adjunction{\mathcal{A}}{\mathcal{A}_{\T}}{}{}$~\cite{Kleisli},
and the Eilenberg-Moore construction
$\adjunction{\mathcal{A}}{\mathcal{A}^{\T}}{}{}$~\cite{EM}. These
two solutions are the initial and terminal solution to the problem
of constructing such an adjunction, in the general sense explained
in~\cite{LS1}. Similarly, a comonoid in
$[\mathcal{A},\mathcal{A}]$ is known as a comonad, and these
constructions work equally well to create a pair of adjoint
functors where the functor $\mathcal{A}^{\T} \to \mathcal{A}$ is
now the left adjoint.

An interesting situation arises when the endofunctor $T \maps
\mathcal{A} \to \mathcal{A}$ defining the monad has a specified
right adjoint $G$.  In this case, Eilenberg and Moore showed that
$G$ can be equipped with the structure of a comonad $\G$, and that
the Eilenberg-Moore category of coalgebras $\mathcal{A}^{\G}$ for
the comonad $\G$ is isomorphic to the Eilenberg-Moore category of
algebras $\mathcal{A}^{\T}$ for the monad $\T$~\cite{EM}. The
isomorphism $\mathcal{A}^{\T} \cong \mathcal{A}^{\G}$ also has the
property that it commutes with the forgetful functors into
$\mathcal{A}$.

If the functor $T$ is equipped with a natural transformation from
$T$ to the identity of $\mathcal{A}$ such that precomposition with
the monad multiplication is the counit for a specified self
adjunction, then we call the resulting structure a Frobenius
monad. This turns out to be the same as a Frobenius object in
$[\mathcal{A},\mathcal{A}]$. Street uses these results to show
that a Frobenius monad always arises from a pair of adjoint
functors that are both left and right adjoints --- an ambijunction
in $\cat{Cat}$.

The approach outlined above is essentially the one taken in this
paper.  Monads and adjunctions can be defined in any 2-category,
and many of the properties of monads and adjunctions in
$\cat{Cat}$ carry over to this abstract context. For instance,
every adjunction in a 2-category $\mathcal{K}$ gives rise to a
monad on an object of $\mathcal{K}$. However, it is not always
true that one can find an adjunction generating a given monad.
This can be attributed to the lack of an object in $\mathcal{K}$
to play the role of the Eilenberg-Moore category of algebras (or
the lack of a Kleisli object, but we will focus on Eilenberg-Moore
objects in this paper). When such an object does exist it is
called an Eilenberg-Moore object for the monad $\T$. The existence
of Eilenberg-Moore objects in a 2-category $\mathcal{K}$ is a
completeness property of the 2-category in question. In
particular, $\mathcal{K}$ has Eilenberg-Moore objects if it is
finitely complete as a 2-category~\cite{LS1,Street3}.

Recall that every bicategory is biequivalent to a strict
2-category, and hence every monoidal category is biequivalent to a
strict monoidal category. Let $\mathcal{M}$ be a monoidal category
and denote as $\Sigma(M)$ the suspension of a strictification of
$M$. Then since the 2-category $\Sigma(M)$ has only one object,
say $\bullet$, a Frobenius object in $M$ is just a Frobenius monad
on the object $\bullet$ in the 2-category $\Sigma(M)$.  It is
tempting to use Eilenberg and Moore's theorem on adjoint monads to
conclude that this Frobenius monad arose from an ambijunction, but
their construction used the fact the 2-category $\mathcal{K}$ was
$\cat{Cat}$. Since $\cat{Cat}$ is finitely complete as a
2-category, this allows the construction of Eilenberg-Moore
objects which are a crucial ingredient in Eilenberg and Moore's
result. Considering Frobenius monads in $\Sigma(M)$, the
strictification of the suspension of the monoidal category $M$, it
is apparent that the required Eilenberg-Moore object is unlikely
to exist: the 2-category $\Sigma(M)$ has only one object!
Fortunately, there is a categorical construction that enlarges a
2-category into one that has Eilenberg-Moore objects. This is
known as the Eilenberg-Moore completion and it will be discussed
in greater detail later.  The important aspect to bear in mind is
that this construction produces a 2-category together with an
ambijunction generating our Frobenius object.

Frobenius objects have found tremendous use in topology,
particularly in the area of topological quantum field theory.  A
well known result going back to Dijkgraaf~\cite{Dij1} states that
2-dimensional topological quantum field theories are equivalent to
commutative Frobenius algebras, see also~\cite{Abrams1,Kock}. Our
result then indicates that:
\begin{quote}
{\it Every 2D topological quantum field theory arises from an
ambijunction in some 2-category.}
\end{quote}
More recently, higher-dimensional analogs of Frobenius algebras
have begun to appear in higher-dimensional topology.  For example,
instances of categorified Frobenius structure have appeared in 3D
topological quantum field theory~\cite{Tillmann}, Khovanov
homology --- the homology theory for tangle cobordisms
generalizing the Jones polynomial~\cite{khovanov}, and the theory
of thick tangles~\cite{Lau2}.

In all of the cases mentioned above, the higher-dimensional
Frobenius structures can be understood as instances of a single
unifying notion --- a `Frobenius pseudomonoid'. A Frobenius
pseudomonoid is a categorified Frobenius algebra --- a monoidal
category satisfying the axioms of a Frobenius algebra up to
coherent isomorphism.  Being inherently categorical, our approach
to solving the problem of constructing adjunctions from Frobenius
objects suggests a quite natural procedure for not only defining a
Frobenius pseudomonoid\footnote{Note that our definition of
Frobenius pseudomonoid nearly coincides with the notion given by
Street~\cite{StreetFrob}; the slight difference is that certain
isomorphisms in our definition are made explicit.}, but more
importantly, for showing that:
\begin{quote}
{\it Every Frobenius pseudomonoid in a semistrict monoidal
2-category arises from a pseudo ambijunction in a semistrict
3-category}.
\end{quote}

The categorified theorem as stated above takes place in the
context of a semistrict 3-category, also referred to as a
$\cat{Gray}$-category.  We take this as a sufficient context for
the generalization since every tricategory or weak 3-category is
triequivalent to a $\cat{Gray}$-category~\cite{GPS}. A
$\cat{Gray}$-category can be defined quite simply using enriched
category theory~\cite{kel}. Specifically, a $\cat{Gray}$-category
is a category enriched in $\cat{Gray}$.  Although a more explicit
definition of a $\cat{Gray}$-category can be given, see for
instance Marmolejo~\cite{mar}, we will not be needing it for this
paper. Adjunctions as well as monads generalize to this context
and are called pseudoadjunctions and pseudomonads, respectively.
They consist of the usual data, where the axioms now hold up to
coherent isomorphism. In the context of an arbitrary
$\cat{Gray}$-category we extend the notion of mateship under
adjunction to the notion of mateship under pseudoadjunction.
Eilenberg-Moore objects and the Eilenberg-Moore completion also
make sense in this context, so we are able to categorify Eilenberg
and Moore's theorem on adjoint monads, as well as our theorem
relating Frobenius objects to ambijunctions, to the context of an
arbitrary $\cat{Gray}$-category.

We remark that Street has demonstrated that the condition for a
monoidal category to be a Frobenius pseudomonoid is identical to
the condition of $*$-autonomy~\cite{StreetFrob}. These
$*$-autonomous monoidal categories are known to have an
interesting relationship with quantum groups and quantum groupoids
\cite{DS2}. Combined with our result relating Frobenius
pseudomonoids to pseudo ambijunctions, the relationship with
$*$-autonomous categories may have implications to quantum groups,
as well as the field of linear logic where $*$-autonomous
categories are used extensively.

\section{Adjoint monads and Frobenius objects}

\subsection{Preliminaries}

In this section we review the concepts of adjunctions and monads
in an arbitrary 2-category along with some of the general theory
needed later on. A good reference for much of the material
presented in this section is~\cite{ks1}.

\begin{defn}
An {\em adjunction} $i,e\maps F \dashv U \maps A \to B$ in a
2-category $\mathcal{K}$ consists of
\begin{itemize}
   \item morphisms $U \maps A \to B$ and $F \maps B \to A$, and
   \item 2-morphisms $i \maps 1 \To UF$ and $e \maps FU \To 1$,
\end{itemize}
such that
\[
 \vcenter{
 \xymatrix{
 & UFU \ar@{=>}[dr]^{Ue} \\
 U \ar@{=>}[ur]^{iU} \ar@{=>}[rr] && U
 }}
\qquad {\rm and} \qquad
 \vcenter{\xymatrix{
 & FUF \ar@{=>}[dr]^{eF} \\
 F \ar@{=>}[ur]^{Fi} \ar@{=>}[rr] && F
 }}
\]
commute.
\end{defn}

\begin{prop} \label{defcompadj}
If $i,e\maps F \dashv U\maps A \to B$ and $i',e' \maps F' \dashv
U'\maps B \to C$ are adjunctions in the 2-category $\mathcal{K}$,
then $FF' \dashv U'U$ with unit and counit:
 \ban
 \bar{i} &:=& \xymatrix@C=2.2pc{
 1 \ar@{=>}[r]^-{i'} & U'F' \ar@{=>}[r]^-{U'iF'} & U'UFF'
 } \\
 \bar{e} &:=& \xymatrix@C=2.2pc{
 FF'U'U \ar@{=>}[r]^-{Fe'U} & FU \ar@{=>}[r]^-{e} & 1
 }
 \ean
\end{prop}

\Proof  We must verify that the triangle identities are satisfied.
The proof is as follows:

 \[
 \xy
   (-40,-15)*+{U'U}="bl";
   (0,-15)*+{U'U}="bm";
   (40,-15)*+{U'U}="br";
   (-20,0)*+{U'F'U'U}="l";
   (20,0)*+{U'UFU}="r";
   (0,15)*+{U'UFF'U'U}="t";
        {\ar@{=>}_1 "bl";"bm"};
        {\ar@{=>}_1 "bm";"br"};
        {\ar@{=>}^{i'U'U} "bl";"l"};
        {\ar@{=>}^<<<<<<<{U'e'U} "l";"bm"};
        {\ar@{=>}^>>>>>>>{U'Ui} "bm";"r"};
        {\ar@{=>}^{U'eU} "r";"br"};
        {\ar@{=>}^{U'iF'U'U} "l";"t"};
        {\ar@{=>}^{U'UFe'U} "t";"r"};
 \endxy
\]

\[
 \xy
   (-40,-15)*+{FF'}="bl";
   (0,-15)*+{FF'}="bm";
   (40,-15)*+{FF'.}="br";
   (-20,0)*+{FF'U'F'}="l";
   (20,0)*+{FUFF'}="r";
   (0,15)*+{FF'U'UFF'}="t";
        {\ar@{=>}_1 "bl";"bm"};
        {\ar@{=>}_1 "bm";"br"};
        {\ar@{=>}^{FF'i'} "bl";"l"};
        {\ar@{=>}^<<<<<<<{Fe'F'} "l";"bm"};
        {\ar@{=>}^>>>>>>>{iFF'} "bm";"r"};
        {\ar@{=>}^{FeF'} "r";"br"};
        {\ar@{=>}^{FF'U'iF'} "l";"t"};
        {\ar@{=>}^{Fe'UFF'} "t";"r"};
 \endxy
\]
where the inner squares commute by the interchange law of the
2-category $\mathcal{K}$. \qed

We recall the notion of mateship under adjunction.

\begin{defn} \label{mates}
Let $i,e \maps F \dashv U\maps A \to B$ and $i',e' \maps F' \dashv
U' \maps A' \to B'$ in the 2-category $\mathcal{K}$.  It was shown
by Kelly and Street~\cite{ks1} that if $a \maps A \to A'$ and $b
\maps B \to B'$, then there is a bijection between 2-morphisms
$\xi\maps bU \To U'a$ and 2-morphisms $\zeta \maps F'b \To a F$,
where $\zeta$ is given in terms of $\xi$ by the composite:
\[
\zeta = \xymatrix@C=2.2pc{F'b \ar@{=>}[r]^-{F'bi} & F'bUF
\ar@{=>}[r]^-{F'\xi F} & F'U'aF \ar@{=>}[r]^-{e'aF} & aF }
\]
and $\xi$ is given in terms of $\zeta$ by the composite:
\[
\xi = \xymatrix@C=2.2pc{bU \ar@{=>}[r]^-{i'bU} & U'F'bU
\ar@{=>}[r]^-{U'\zeta U} & U'aFU \ar@{=>}[r]^-{U'ae} & U'a }.
\]
Under these circumstances we say that $\xi$ and $\zeta$ are {\em
mates under adjunction} and we sometimes write $\xi \dashv \zeta$.
\end{defn}

The naturality of this bijection can be expressed as an
isomorphism of certain double categories, see Proposition 2.2
\cite{ks1}. In both cases, the objects of the double categories
are those of $\mathcal{K}$. The horizontal arrows are the
morphisms of $\mathcal{K}$ with the usual composition and the
vertical arrows are the adjunctions in $\mathcal{K}$ with the
composition given in Proposition~\ref{defcompadj}. In the first
double category, a square with sides $a \maps A \to A'$, $b \maps
B \to B'$,$i,e \maps F \dashv U\maps A \to B$, and $i',e' \maps F'
\dashv U'\maps A' \to B'$ is a 2-cell $\xi\maps bU \To U'a$.  In
the second double category a square with the same sides is a
2-cell $\zeta \maps F'b \To a F$. The isomorphism between these
two double categories makes precise the idea that the association
of mateship under adjunction respects composites and identities
both of adjunctions and of morphisms in $\mathcal{K}$.

\begin{defn}
A {\em monad} $\T =(T,\mu, \eta)$ in a 2-category $\mathcal{K}$ on
the object $B$ of $\mathcal{K}$ consists of an endomorphism $T
\maps B \to B$ together with 2-morphisms:
\begin{itemize}
    \item {\em multiplication} for the monad: $\mu \maps T^2 \To T$, and
    \item {\em unit} for the monad: $\eta \maps 1 \To T$,
\end{itemize}
such that
\[
\vcenter{ \xymatrix{ T \ar@{=>}[r]^-{T\eta} \ar@{=}[dr] & TT
  \ar@{=>}[d]^-{\mu} & T \ar@{=>}[l]_-{\eta T} \ar@{=}[dl] \\ & T}}
\qquad {\rm and} \qquad
 \vcenter{\xymatrix{
 T^3 \ar@{=>}[r]^{\mu T} \ar@{=>}[d]_{T\mu} & T^2 \ar@{=>}[d]^{\mu}
 \\ T^2 \ar@{=>}[r]_{\mu} & T
 }}
\]
commute.
\end{defn}

A comonad is defined by reversing the directions of the 2-cells:

\begin{defn}
A {\em comonad} $\G =(G,\delta, \varepsilon)$ in a 2-category
$\mathcal{K}$ on the object $B$ of $\mathcal{K}$ consists of an
endomorphism $G \maps B \to B$ together with 2-morphisms:
\begin{itemize}
    \item {\em comultiplication} for the comonad: $\delta \maps G \To G^2$, and
    \item {\em counit} for the comonad: $\varepsilon \maps G \To 1$,
\end{itemize}
such that
\[
\vcenter{ \xymatrix{ G \ar@{=}[dr] &  \ar@{=>}[l]_-{G\varepsilon}
GG \ar@{=>}[r]^-{\varepsilon G}
   & G  \ar@{=}[dl] \\ & G \ar@{=>}[u]_-{\delta}}}
\qquad {\rm and} \qquad
 \vcenter{\xymatrix{
 G^3   & G^2 \ar@{=>}[l]_{\delta G}
 \\ G^2 \ar@{=>}[u]^{G\delta}  & G\ar@{=>}[l]^{\delta}\ar@{=>}[u]_{\delta}
 }}
\]
commute.
\end{defn}

A complete treatment of monads in this generality is presented in
\cite{LS2,LS1}. It is clear that if $i,e\maps F \dashv U \maps A
\to B$ is an adjunction in $\mathcal{K}$, then $(UF,UeF,i)$ is a
monad on $B$. We now recall a result due to Eilenberg-Moore
\cite{EM}, proven in the context $\mathcal{K}=\cat{Cat}$, that
easily generalizes to arbitrary $\mathcal{K}$.

\begin{prop} \label{defadjmonad}
Let $\T =(T,\mu, \eta)$ be a monad on an object $B$ in a
2-category $\mathcal{K}$ such that the endomorphism $T \maps B \to
B$ has a specified right adjoint $G$ with counit $\sigma \maps TG
\to 1$ and unit $\iota \maps 1 \to GT$. Then $\G =
(G,\varepsilon,\delta)$ is a comonad where $\varepsilon$ and
$\delta$ are the mates under adjunction of $\eta$ and $\mu$ with
the explicit formulas being:
 \ban
 \varepsilon &=& \sigma.\eta G \quad \\
  \delta &=&
 G^2\sigma.G^2\mu G. G\iota TG. \iota G
 \ean
and $\G$ is said to be a comonad {\em right adjoint} to the monad
$\T$, denoted $\T \dashv \G$.
\end{prop}

\Proof  This statement immediately follows from the composition
preserving property of the bijection of mates under adjunction.
Since $T \dashv G$, $\mu \dashv \delta$, $\eta \dashv
\varepsilon$, and mateship under adjunction preserves the various
composites, then because $(T,\mu,\eta)$ satisfies the monad
axioms, $(G,\delta,\varepsilon)$ will satisfy the comonad axioms.
\qed

\begin{defn}
A monad $\T$ in the 2-category $\mathcal{K}$ is called a {\em
Frobenius monad} if it is equipped with a morphism
$\varepsilon\maps T \to 1$ such that $\varepsilon.\mu$ is the
counit for an adjunction $T \dashv T$.
\end{defn}

The notion of a Frobenius monad (or Frobenius standard
construction as it was originally called) was first defined by
Lawvere~\cite{law}. In Street~\cite{StreetFrob} several
definitions of Frobenius monad are given and proven equivalent. If
one regards the monoidal category $\cat{Vect}$ as a one object
2-category $\Sigma(\cat{Vect})$, then a Frobenius monad in
$\Sigma(\cat{Vect})$ is just the usual notion of a Frobenius
algebra.

\begin{defn}
An {\em action} of the monad $\T$ on a morphism $s \maps A \to B$
in the 2-category $\mathcal{K}$ is a 2-morphism $\nu \maps T s \To
s$ such that
\[
 \vcenter{\xymatrix{
 s \ar@{=>}[r]^{\eta s} \ar@{=}[dr] & T s \ar@{=>}[d]^{\nu} \\
 & s
 }}
\qquad {\rm and} \qquad
 \vcenter{\xymatrix{
 T^2 s \ar@{=>}[r]^{\mu s} \ar@{=>}[d]_{T \nu} & T s
 \ar@{=>}[d]^{\nu} \\ T s \ar@{=>}[r]_{\nu} & s
 }}
\]
commute. A morphism $s$ together with an action is called a {\em
$\T$-algebra} (with domain $A$).
\end{defn}

For any morphism $s \maps A \to B$ in $\mathcal{K}$, $Ts$ with
action $\mu s \maps T^2 s \To T s$ is a $\T$-algebra. For reasons
that will soon become apparent we call the $\T$-algebra $(Ts,\mu
s)$ a {\em free $\T$-algebra}. The traditional notion of
$\T$-algebra corresponds to the notion presented above when
$\mathcal{K} = \cat{Cat}$ and $A$ is the one object category.  In
this case we identify the map $s \maps 1 \to B$ with its image.

\begin{defn}
Let $\T$ be a monad in $\mathcal{K}$. For each $A$ in
$\mathcal{K}$ define the category \cat{$\T$-Alg}$_{\bf A}$ whose
objects are $\T$-algebras, and whose morphisms between
$\T$-algebras $(s,\nu)$ and $(s',\nu')$ are those 2-morphisms $h
\maps s \To s'$ of $\mathcal{K}$ making
\[
 \xymatrix{
 Ts \ar@{=>}[r]^{T h} \ar@{=>}[d]_{\nu} &
 Ts' \ar@{=>}[d]^{\nu'} \\ s \ar@{=>}[r]_{h} & s'
 }
\]
commute.  We call the morphisms in \cat{$\T$-Alg}$_{\bf A}$ {\em
morphisms of $\T$-algebras}.
\end{defn}

Given a morphism $K\maps A' \to A$ in $\mathcal{K}$, one can
define a change of base functor $\hat{K}\maps$
\cat{$\T$-Alg}$_{\bf A}$ $\to$ \cat{$\T$-Alg}$_{\bf A'}$.  If $ h
\maps(s,\nu)\to(s',\nu')$ is in \cat{$\T$-Alg}$_{\bf A}$, then its
image under $\hat{K}$ is $ hK \maps(sK,\nu K)\to(s'K,\nu'K)$. If
$k \maps K \To K'$ in $\mathcal{K}$ then we get a natural
transformation $\hat{k}\maps \hat{K} \To \hat{K'}$ such that
$\hat{k}_{(s,\nu)}=sk$.   In fact, this shows that the
construction of $\T$-algebras defines a 2-functor
$\cat{$\T$-Alg}\maps \mathcal{K}^{\op} \to \cat{Cat}$.

As with the case when $\mathcal{K} = \cat{Cat}$ we have a
forgetful functor:
 \ban
  U_A^{\T} \maps \cat{$\T$-Alg}_{\bf A} &\to& \mathcal{K}(A,B) \\
   (s,\nu) &\mapsto& s \\
   h \maps s \To s' &\mapsto& h \maps s \To s'
 \ean
and this functor has a left adjoint:
 \ban
F_A^{\T} \maps \mathcal{K}(A,B) &\to& \cat{$\T$-Alg}_{\bf A} \\
   s &\mapsto& (Ts,\mu s) \\
   h \maps s \To s' &\mapsto& Th\maps (Ts,\mu s ) \To (Ts',\mu s').
 \ean
The unit of the adjunction is the natural transformation $i^{\T}_A
\maps 1 \To U^{\T}_AF^{\T}_A$ given by sending the morphism
$s\maps A \to B$ to the 2-morphism $\eta s$.  The counit $e^{\T}_A
\maps F^{\T}_AU^{\T}_A \To 1$ is the natural transformation that
assigns to each $\T$-algebra $(s,\nu)$ the morphism of
$\T$-algebras given by $\nu \maps Ts \To s$.  This adjunction
exists for every $A$ in $\mathcal{K}$.  In fact, we have the
following:

\begin{prop} \label{Tadj}
 The collection of adjunctions
\[
 i^{\T}_A,e^{\T}_A\maps F^{\T}_A \To U^{\T}_A  \maps \cat{$\T$-Alg}_{\bf
A}\to \mathcal{K}(A,B)
\]
defined for each $A$ in $\mathcal{K}$ defines an adjunction
\[
 i^{\T},e^{\T} \maps F^{\T} \To U^{\T}  \maps \cat{$\T$-Alg}\to \mathcal{K}(-,B)
\]
in the 2-category $[\mathcal{K}^{\op},\cat{Cat}]$ consisting of
2-functors $\mathcal{K}^{\op}\to\cat{Cat}$, 2-natural
transformations between them, and modifications.
\end{prop}

\proof We have already shown above that $\cat{$\T$-Alg}$ is a
2-functor.  It is also clear that the collection of natural
transformations $U^{\T}_A$ define a 2-natural transformation $
\cat{$\T$-Alg}\To\mathcal{K}(-,B)$.  To see that the collection of
$F^{\T}_A$ define a 2-natural transformation $F^{\T} \maps
\mathcal{K}(-,B) \to \cat{$\T$-Alg}$ we must verify the
1-naturality:
\[
 \xymatrix{
\mathcal{K}(A,B) \ar[r]^{F^{\T}_A} \ar[d]_{\mathcal{K}(K,B)} &
\cat{$\T$-Alg} \ar[d]^{\hat{K}} \\ \mathcal{K}(A',B)
\ar[r]_{F^{\T}_{A'}} & \cat{$\T$-Alg} }
\]
which commutes since both functors map $h \maps s \To s'$ to
$ThK\maps TsK \To Ts'K$, and we must verify the 2-naturality:
\[
 \xy
   (-32,0)*+{\mathcal{K}(A,B)}="l";
   (0,0)*+{\mathcal{K}(A',B)}="m";
   (25,0)*+{\cat{$\T$-Alg}_{\bf A'}}="r";
        {\ar@/^1.3pc/^{\mathcal{K}(K,B)} "l";"m"};
        {\ar@/_1.3pc/_{\mathcal{K}(K',B)}  "l";"m"};
        {\ar^{F^{\T}_{A'}}  "m";"r"};
        {\ar@{=>}^{\mathcal{K}(k,B)} (-20,3);(-20,-3) };
 \endxy
\quad = \quad
 \xy
   (-24,0)*+{\mathcal{K}(A,B)}="l";
   (0,0)*+{\cat{$\T$-Alg}_{\bf A}}="m";
   (28,0)*+{\cat{$\T$-Alg}_{\bf A'} .}="r";
        {\ar@/^1.3pc/^{\hat{K}} "m";"r"};
        {\ar@/_1.3pc/_{\hat{K}'}  "m";"r"};
        {\ar^{F^{\T}_{A}}  "l";"m"};
        {\ar@{=>}^{\hat{k}} (15,3);(15,-3) };
 \endxy
\]
This equality holds since both natural transformations assign to
the morphism $s \maps A \to B$ the morphism of $\T$-algebras $Tsk
\maps TsK \to TsK'$.

Next we verify that the collection of $i^{\T}_A$ define a
modification $i^{\T} \maps 1 \To U^{\T}F^{\T}$.  Consider the
diagrams below:
\[
 \xy
   (-32,0)*+{\mathcal{K}(A,B)}="l";
   (-4,0)*+{\mathcal{K}(A,B)}="m";
   (26,0)*+{\mathcal{K}(A',B)}="r";
        {\ar@/^1.3pc/^{1_\mathcal{K}(A,B)} "l";"m"};
        {\ar@/_1.3pc/_{U^{\T}_{A}F^{\T}_{A}}  "l";"m"};
        {\ar^-{\mathcal{K}(K,B)}  "m";"r"};
        {\ar@{=>}^{i^{\T}_{A}} (-20,3);(-20,-3) };
 \endxy
\quad  \quad
 \xy
   (-26,0)*+{\mathcal{K}(A,B)}="l";
   (4,0)*+{\mathcal{K}(A',B)}="m";
   (32,0)*+{\mathcal{K}(A',B)}="r";
        {\ar@/^1.3pc/^{1_{\mathcal{K}(A',B)}} "m";"r"};
        {\ar@/_1.3pc/_{U^{\T}_{A'}F^{\T}_{A'}}  "m";"r"};
        {\ar^-{\mathcal{K}(K,B)}  "l";"m"};
        {\ar@{=>}^{i^{\T}_{A'}} (16,3);(16,-3) };
 \endxy
\]
Both of these natural transformations assign to the morphism $s
\maps A \to B$ the morphism $\eta sK$ of $\T$-algebras.  Hence,
$i^{\T}$ is a modification.  Similarly one can check that the
collection of $e^{\T}_A$ define a modification $e^{\T} \maps
F^{\T}U^{\T}\To 1$.  Since the coherence axioms for an adjunction
are verified pointwise we have shown that  $i^{\T},e^{\T} \maps
F^{\T} \To U^{\T}  \maps \cat{$\T$-Alg}\to \mathcal{K}(-,B)$ is an
adjunction in $[\mathcal{K}^{\op},\cat{Cat}]$. \qed

\begin{defn} \label{defEM}
We say that an {\em Eilenberg-Moore object exists for a monad
$\T$} if the 2-functor $\cat{$\T$-Alg}\maps \mathcal{K}^{\op} \to
\cat{Cat}$ is representable.  An Eilenberg-Moore object for the
monad $T$ is then just a choice of representation for the
2-functor $\cat{$\T$-Alg}$, that is an object $B^{\T}$ of
$\mathcal{K}$ together with a specified 2-natural isomorphism from
\cat{$\T$-Alg} to the 2-functor $\mathcal{K}(-,B^{\T})$.
\end{defn}

If an Eilenberg-Moore object exists for a monad $\T$ then by the
enriched Yoneda lemma, or 2-categorical Yoneda lemma as it is
sometimes referred to in this context, the adjunction of
Proposition~\ref{Tadj} arises from an adjunction
$i^{\T},e^{\T}\maps F^{\T} \dashv U^{\T} \maps B \to B^{\T}$ in
$\mathcal{K}$.

Given a comonad $\G$ in $\mathcal{K}$ we can define the category
$\cat{$\G$-CoAlg}_A$ of $\G$-coalgebras $(s,\bar{\nu})$ and maps
of $\G$-coalgebras by reversing the directions of the 2-cells in
the definition of $\cat{$\T$-Alg}_A$ and substituting the
appropriate data for $\G$. If $K \maps A' \to A$ then we can
define a change of base functor $\hat{K}\maps \cat{$\G$-CoAlg}_A
\to \cat{$\G$-CoAlg}_{A'}$ sending $h\maps (s,\bar{\nu}) \to
(s,\bar{\nu})$ to $hK \maps (sK,\bar{\nu}K) \to (sK,\bar{\nu}K)$.
Further, if $k\maps K \to K'$ in $\mathcal{K}$ we can define the
natural transformation $\hat{k}$ which sends $(s,\bar{\nu})$ to
$sk \maps (sK,\bar{\nu}K) \to (sK',\bar{\nu}K')$. Hence we have a
2-functor $\cat{$\G$-CoAlg}\maps \mathcal{K}^{\op} \to \cat{Cat}$.
As before, there is a forgetful 2-natural transformation $U^{\G}
\maps \cat{$\G$-CoAlg} \to \mathcal{K}(-B)$.  However, in this
case, $U^{\G}$ has a right adjoint $F^{\G}$.  An Eilenberg-Moore
object for a comonad is just a choice of representation for the
2-functor \cat{$\G$-CoAlg}.  If an Eilenberg-Moore object for $\G$
does exist then, again by the 2-categorical Yoneda lemma, the
adjunction $i^{\G},e^{\G}\maps F^{\G} \vdash U^{\G}\maps
\cat{$\G$-CoAlg} \to \mathcal{K}(-,B)$ arises from an adjunction
 $i^{\G},e^{\G}\maps F^{\G} \vdash U^{\G}\maps
B \to B^{\G}$ in $\mathcal{K}$.

\subsection{Adjoint monads}

Given $\T \dashv \G$ in $\cat{Cat}$, Eilenberg and Moore~\cite{EM}
showed that mateship under adjunction of action and coaction
defines an isomorphism $B^{\T} \cong B^{\G}$ of categories between
the Eilenberg-Moore category of $\T$-algebras for the monad $\T$
and the Eilenberg-Moore category of $\G$-coalgebras for the
comonad $\G$.  In this section we continue the program for the
formal theory of monads begun by Street~\cite{LS1}.  In
particular, we extend the classical theory of adjoint monads
developed by Eilenberg and Moore to the context of an arbitrary
2-category.

\begin{lem} \label{lemXYZ}
 Let $\T$ be a monad on $B \in \mathcal{K}$.  If $\iota,\sigma \maps T \dashv
G$ and $G$ is equipped with the comonad structure $\G$ from
Proposition~\ref{defadjmonad}, then the category
$\cat{$\T$-Alg}_A$ is isomorphic to the category
$\cat{$\G$-CoAlg}_A$ and this isomorphism commutes with the
forgetful functors to $\mathcal{K}(A,B)$.
\end{lem}

\proof Define a functor $\mathcal{M}_A \maps \cat{$\T$-Alg}_A \to
\cat{$\G$-CoAlg}_A$ by sending $h\maps (s, \nu) \to (s',\nu')$ to
$h \maps (s,\bar{\nu}) \to (s',\bar{\nu}')$ where $\nu \dashv
\bar{\nu}$ and $\nu' \dashv \bar{\nu}'$.  The explicit formulas
for $\bar{\nu}$ and $\bar{\nu}'$ are $G\nu.\iota s$ and
$G\nu'.\iota s'$.  Note that $h\maps s \to s'$ corresponds via
mateship under adjunction to itself, and $Th$ corresponds via
mateship to $Gh$. Hence $(s,\bar{\nu})$ is a $\G$-coalgebra since
$(s, \nu)$ is a $\T$-algebra and the association of mates
preserves composites by the remarks following
Definition~\ref{mates}. Since $\mathcal{M}_A$ is the identity on
morphisms it is clear that composites and identities are
preserved.

We define the inverse functor $\overline{\mathcal{M}}_A \maps
\cat{$\G$-CoAlg}_A \to \cat{$\T$-Alg}_A$ again using mateship
under adjunction. The $\G$-coalgebra $(s,\bar{\nu})$ is sent to
the $\T$-algebra $(s,\tilde{\nu})$ where $\tilde{\nu} \dashv
\bar{\nu}$.  The explicit formula for $\tilde{\nu}$ is $\sigma
s.T\bar{\nu}$.  On morphisms $\overline{\mathcal{M}}_A$ is the
identity.  The fact that $\mathcal{M}_A$ and
$\overline{\mathcal{M}}_A$ are inverses follows from the triangle
identities for the adjunction $\iota,\sigma \maps T \dashv G$.
Clearly, $U^{\G}_A\mathcal{M}_A=U^{\T}_A$  since $\mathcal{M}_A$
maps $s\maps A \to B$ to itself and is the identity on morphisms.
\qed

\begin{thm} [The adjoint monad theorem] \label{EMthm} \label{preEMthm}
Let $\T$ be a monad in $\mathcal{K}$ with $T \dashv G$ and denote
the induced comonad of Proposition~\ref{defadjmonad} as $\G$. Then
there is a 2-natural isomorphism
$\mathcal{M}\maps\cat{$\T$-Alg}\to\cat{$\G$-CoAlg}$ 
making the following diagram
\[
 \xy
 (-20,20)*+{ \cat{$\T$-Alg}}="1";
 (20,20)*+{\cat{$\G$-CoAlg}}="2";
 (0,5)*+{\mathcal{K}(-,B)}="3";
    {\ar^-{\mathcal{M}} "1";"2" };
    {\ar_{U^{\T}} "1";"3" };
    {\ar^{U^{\G}} "2";"3" };
 \endxy
\]
commute. Furthermore, if one exists, an Eilenberg-Moore object
$B^{\T}$ for the monad $\T$ serves as an Eilenberg-Moore object
$B^{\G}$ for the comonad $\G$. So that the above diagram arises
via the 2-categorical Yoneda lemma from the commutative diagram:
\[
 \xy
 (-20,20)*+{B^{\T}}="1";
 (20,20)*+{B^{\G}}="2";
 (0,5)*+{B}="3";
    {\ar^-{\mathcal{M}} "1";"2" };
    {\ar_{U^{\T}} "1";"3" };
    {\ar^{U^{\G}} "2";"3" };
 \endxy
\]
in $\mathcal{K}$.
\end{thm}

\proof   We show that the collection of natural isomorphisms
$\mathcal{M}_A$ defined in Lemma~\ref{lemXYZ} define a 2-natural
isomorphism $\mathcal{M} \maps \cat{$\T$-Alg} \to
\cat{$\G$-CoAlg}$.  The 1-naturality of $\mathcal{M}$ follows from
the fact that if $K \maps A' \to A$, then
\[
\hat{K}\mathcal{M}_A(s,\nu)  =  \big( sK, (g\nu.\iota s)K\big)
  =  \big( sK, g\nu K.\iota sK\big)
  =  \mathcal{M}_{A'}\hat{K}(s,\nu) .
\]
The 2-naturality of $\mathcal{M}$ follows from the fact that
$\mathcal{M}_{A}$ is the identity on morphisms. Hence,
$\mathcal{M}$ is a 2-natural transformation.  From
Lemma~\ref{lemXYZ} it is clear that $\mathcal{M}$ commutes with
the forgetful functors since this is verified pointwise.

If $B^{\T}$ is an Eilenberg-Moore object for the monad $T$, then
we have a choice of 2-natural isomorphism $\mathcal{K}(-,B^{\T})
\cong \cat{$\T$-Alg}$.  Composing this 2-natural isomorphism with
the 2-natural isomorphism $\mathcal{M}$ equips $B^{\T}$ with the
structure of an Eilenberg-Moore object for the comonad $\G$.
Since the 2-natural isomorphism $\mathcal{M}$ commutes with the
forgetful 2-natural isomorphisms $U^{\T}$ and $U^{\G}$, it is
clear that their images under the 2-categorical Yoneda lemma will
make the required diagram commute. \qed

This theorem shows that if the monad $\T$ has an adjoint comonad
$\G$, and if the Eilenberg-Moore objects exists, then the
`forgetful morphism' $U^{\T} \maps B^{\T} \to B$ has not only a
left adjoint $F^{\T}$, but also a right adjoint
$\overline{\mathcal{M}} F^{\G}$. We can also extend the classical
converse of this theorem to show that if a morphism has both a
left and right adjoint, then the induced monad and comonad are
adjoint.

\begin{thm} \label{ambInduceadjmon}
Let
$
 \adjunction{B}{C}{L_1}{R}
\quad {\rm and} \quad
 \radjunction{B}{C}{L_2}{R}
$
(or $L_1 \dashv R \dashv L_2$)  be specified adjunctions in the
2-category $\mathcal{K}$. Also, let $\T_1$ be the monad on $B$
induced by the composite $R L_1$, and $\T_2$ be the comonad on $B$
induced by the composite $R L_2$. Then $\T_1 \dashv \T_2$ via a
specified adjunction determined from the data defining the
adjunctions $L_1 \dashv R \dashv L_2$.
\end{thm}

\Proof  Let $i_1,e_1 \maps L_1 \dashv R\maps C \to B$ and
$i_2,e_2\maps R \dashv L_2 \maps B \to C$, then it follows from
the composition of adjunctions that $RL_1 \dashv RL_2$ with $\iota
= Ri_2L_1.i_1\maps 1 \To RL_2RL_1$ and $\sigma = e_2.Re_1L_2 \maps
RL_1RL_2 \To 1$.  The triangle identities follow from the triangle
identities for the pairs $(i_1,e_1)$ and $(i_2,e_2)$.  It remains
to be shown that $\mu=L_1e_1R$ is mates under adjunction with
$\delta = Ri_2L_2$ and $\eta=i_1$ is mates with $\varepsilon=e_2$.

The mate to $e_2$ is given by the composite:
\[
\xymatrix@C=2.6pc{
   1        \ar@{=>}[r]^-{i_1}
 & RL_1     \ar@{=>}[r]^-{Ri_2L_1}
 & RL_2RL_1 \ar@{=>}[r]^-{e_2RL_1}
 & RL_1 }
\]
but, by one of the triangle identities, this is just the map
$i_1$. The mate to $i_1$ is given by the composite:
\[
 \xymatrix@C=2.6pc{
   RL_2        \ar@{=>}[r]^-{i_1RL_2}
 & RL_1RL_2     \ar@{=>}[r]^-{Re_1L_2}
 & RL_2 \ar@{=>}[r]^-{e_2}
 & 1 }
\]
and by the other triangle identity is equal to $e_2$. Hence $\eta
\dashv \varepsilon$. In a similar manner it can be shown that $\mu
\dashv \delta$ using multiple applications of the triangle
identities. \qed

As previously discussed, Street uses the classical version of this
theorem to show that a Frobenius monad in $\cat{Cat}$ is always
induced by ambijunction in $\cat{Cat}$. Furthermore, he shows that
the converse of this theorem is also true
--- corresponding to a Frobenius monad in $\cat{Cat}$ there always
exists an ambijunction generating it. We will now extend this
result to include the case of Frobenius monads in an arbitrary
2-category.  The most significant impedance to this extension is
the lack of Eilenberg-Moore objects in 2-categories that are not
finitely complete.  This leads us to the free completion of a
2-category under Eilenberg-Moore objects to be discussed in the
next section. Note that this free completion was not needed in
Street's construction since $\cat{Cat}$ already has
Eilenberg-Moore objects.

\subsection{Eilenberg-Moore completions} \label{secEMcompletion}

As we saw in the introduction, one of the aims of this paper is
show that every Frobenius object in \textit{any} monoidal category
arises from an ambijunction in some 2-category.  To prove this,
one is tempted to apply Theorem~\ref{EMthm}. However, when
regarding a Frobenius object in a monoidal category as a Frobenius
monad on the suspension of the monoidal category caution must be
exercised. The 2-category $\Sigma(M)$ has only one object.  Thus,
it is unlikely that the Eilenberg-Moore objects, supposed to exist
in Theorem~\ref{EMthm}, actually exists in $\Sigma(M)$.

Street~\cite{Street3} has shown that an Eilenberg-Moore object can
be considered as a certain kind of weighted limit. He has also
shown that the weight is finite in the sense of~\cite{kel2}. In
\textit{The Formal Theory of Monads.\ II.\ }\cite{LS2}, Lack and
Street use this result to show that one can define
\cat{EM}$(\mathcal{K})$; the free completion under Eilenberg-Moore
objects of the 2-category $\mathcal{K}$.  Since the free
completion under a class of colimits is more accessible than
completions under the corresponding limits, Lack and Street first
construct Kl$(\mathcal{K})$ --- the free completion under Kleisli
objects. They then take \cat{EM}$(\mathcal{K})$ to be
Kl$(\mathcal{K}^{\op})^{\op}$. Since a Kleisli object is a
colimit, to construct Kl$(\mathcal{K})$ one must complete
$\mathcal{K}$ embedded in $[\mathcal{K}^{\op},\cat{Cat}]$, by
Yoneda, under the class of $\Phi$-colimits, where $\Phi$ consists
of the weights for Kleisli objects. This amounts to taking the
closure of the representables under $\Phi$-colimits ~\cite{kel}.
By the theory of such completions, we obtain a 2-functor $Z \maps
\mathcal{K} \to \cat{EM}(\mathcal{K})$ with the property that for
any 2-category $\mathcal{L}$ with Eilenberg-Moore objects,
composition with $Z$ induces an equivalence of categories between
the 2-functor category $[\mathcal{K},\mathcal{L} ]$ and the full
subcategory of the 2-functor category
$[\cat{EM}(\mathcal{K}),\mathcal{L}]$ consisting of those
2-functors which preserve Eilenberg-Moore objects~\cite{LS2}.
Furthermore, the theory of completions under a class of colimits
also tells us that $Z$ will be fully faithful.

The Eilenberg-Moore completion can also be given a concrete
description. The object of \cat{EM}$(\mathcal{K})$ are the monads
in $\mathcal{K}$ and the morphisms are the usual morphisms of
monads.  Hence, a morphism from $\T =(T\maps B \to B,\mu, \eta)$
to $\T' =(T'\maps C \to C,\mu', \eta')$ in \cat{EM}$(\mathcal{K})$
is a morphism $F\maps B \to C$ and a 2-morphisms $\phi\maps T'F
\To FT$ of $\mathcal{K}$ satisfying two equations:
\[
\def\objectstyle{\scriptstyle}
\def\labelstyle{\scriptstyle}
 \xymatrix{
 T'T'F
    \ar@{=>}[d]_{\mu'F}
    \ar@{=>}[r]^{T'\phi }
 & T'FT
    \ar@{=>}[r]^{\phi T}
 &  FTT
    \ar@{=>}[d]^{F\mu} \\
  TF
    \ar@{=>}[rr]_{\phi}
 &&  FT
 }
\quad \quad
 \xymatrix{
 & F
    \ar@{=>}[dr]^{F\eta}
    \ar@{=>}[dl]_{\eta'F} \\
 T'F
    \ar@{=>}[rr]_{\phi}
 && FT
 }
\]
A crucial observation made by Lack and Street is that the
2-morphisms in \cat{EM}$(\mathcal{K})$ are \textit{not} the
2-morphisms of the 2-category of monads. Rather, a 2-morphism from
$(F,\phi)$ to $(F',\psi)$ in \cat{EM}$(\mathcal{K})$ consists of a
2-morphism $f\maps F \to F'T$ satisfying
\[
\def\objectstyle{\scriptstyle}
\def\labelstyle{\scriptstyle}
 \xymatrix@C=3pc{
   T'F
    \ar@{=>}[r]^{\phi}
    \ar@{=>}[d]_{T'f}
 & FT
    \ar@{=>}[r]^{fT}
 &  F'TT
    \ar@{=>}[d]^{F'\mu} \\
 T'F'T
    \ar@{=>}[r]_{\psi T }
 & F'TT
    \ar@{=>}[r]_{ F'\mu}
 & F'T
 }
\]

\subsection{$\cat{EM}$-Completions in \cat{Vect}} \label{EMVECT}

The Eilenberg-Moore completion may seem rather substantial, so in
order to gain some insight into this procedure we briefly discuss
the implications of this completion for $\cat{Vect}$.  We will
then construct the adjunction that generates a given monad in
$\Sigma(\cat{Vect})$.   The objects of
\cat{EM}$(\Sigma(\cat{Vect}))$ will be the monads in
$\Sigma(\cat{Vect})$. In this case, a monad in
$\Sigma(\cat{Vect})$ is an algebra in the traditional sense of
linear algebra --- a vector space equipped with an associative,
unital multiplication. For the duration of this example `algebra'
is to be interpreted in this sense; not in the sense of an algebra
for a monad.  A morphism in \cat{EM}$(\Sigma(\cat{Vect}))$ from an
algebra $A_1$ to an algebra $A_2$ amounts to a vector space $V$
together with a linear map $
 \xymatrix@C=1.4pc{
 V \ten A_2 \ar[r]^{\phi} &A_1\ten V
 }
$ such that
\[
\def\objectstyle{\scriptstyle}
\def\labelstyle{\scriptstyle}
 \xymatrix{
 V \ten A_2 \ten A_2
    \ar[d]_{V \ten m_2}
    \ar[r]^{\phi \ten A_2 }
 & A_1 \ten V \ten A_2
    \ar[r]^{A_1 \ten \phi}
 &  A_1 \ten A_1 \ten V
    \ar[d]^{ m_1\ten V} \\
  V \ten A_2
    \ar[rr]_{\phi}
 &&  A_1 \ten V
 }
\quad
 \xymatrix{
 & V
    \ar[dr]^{\iota_1 \ten V}
    \ar[dl]_{ V\ten\iota_2 } \\
 V \ten A_2
    \ar[rr]_{\phi}
 && A_1\ten V
 }
\]
commute, where $(m_1,\iota_1)$ and $(m_2,\iota_2)$ are the
multiplication and unit for the algebras $A_1$ and $A_2$
respectively.  This might be described as a left-free bimodule: a
vector space V with a right $A_2$ action on $A_1 \ten V$ given by
$\xymatrix{A_1 \ten V \ten A_2 \ar[r]^-{A_1 \ten \phi} & A_1 \ten
A_1 \ten V \ar[r]^-{m_1} & A_1 \ten V}$.  This action makes $A_1
\ten V$ into a $(A_1,A_2)$-bimodule.  The composite of morphisms
$(V,\phi)\maps A_1 \to A_2$ and $(V',\phi') \maps A_2 \to A_3$ is
given by $(V \ten V',\phi \ten V' \circ V \ten \phi') \maps A_1
\to A_3$ -- the left-free bimodule $A_1 \ten V \ten V'$.

A 2-morphism in \cat{EM}$(\Sigma(\cat{Vect}))$ from $(V,\phi) \To
(V',\psi)$ is a linear map $\rho \maps V \to A_1 \ten V'$ making
\[
\def\objectstyle{\scriptstyle}
\def\labelstyle{\scriptstyle}
 \xymatrix@C=3pc{
   V\ten A_2
    \ar[r]^{\phi}
    \ar[d]_{\rho \ten A_2}
 &A_1\ten  V
    \ar[r]^{A_1 \ten \rho}
 &  A_1 \ten A_1 \ten V'
    \ar[d]^{m_1 \ten V'} \\
 A_1 \ten V' \ten A_2
    \ar[r]_{A_1 \ten\psi }
 & A_1 \ten A_1 \ten V'
    \ar[r]_{ m_1\ten V'}
 & A_1 \ten V'
 }
\]
commute.  This amounts to saying that a 2-morphism is just a
bimodule homomorphism of left-free bimodules.  To summarize:
\begin{quote}
{\it Every Frobenius algebra in \cat{Vect} will be shown to arise
from an ambijunction in the 2-category
\cat{EM}$(\Sigma(\cat{Vect}))$ consisting of: algebras, left-free
bimodules, and bimodule homomorphisms}.
\end{quote}
Recall that the Eilenberg-Moore completion was obtained from the
Kleisli completion as Kl$(\mathcal{K}^{\op})^{\op}$.  Hence, a
similar description of the Kleisli completion of
$\Sigma(\cat{Vect})$ can be given in terms of right-free
bimodules\footnote{This description of the Eilenberg-Moore
completion and Kleisli completion was explained to the author by
Steve Lack.}.

Ambijunctions in the Eilenberg-Moore completion of
$\Sigma(\cat{Vect})$ correspond to the notion of a {\it Frobenius
extension} familiar to algebraists, see for
example~\cite{Kadison}. For an algebra $A$ over the field $k$ we
have the inclusion map $\iota \maps k \to A$.  The category of
$A$-modules corresponds to the category of algebras for the monad
$A$ in $\Sigma(\cat{Vect})$. The restriction functor ${\rm
Res}\maps A{\rm -mod} \to k{\rm -mod}$ has left and right adjoint
functors: the induction functor ${\rm Ind}(M)=A \otimes_k M$ and
coinduction ${\rm CoInd}(M)=\Hom_k(A,M)$.  When $A$ is a Frobenius
algebra in $\cat{Vect}$ these functors are isomorphic defining an
ambijunction generating $A$.

\subsection{Frobenius monads and ambijunctions}

\begin{lem} \label{lemFM}
Let $\T=(T,\mu,\eta,\varepsilon)$ be a Frobenius monad on
$\mathcal{K}$ with $\iota,\varepsilon.\mu\maps T \dashv T$. For
notational convenience, denote the induced comonad of
Proposition~\ref{defadjmonad} on $T$ as $\G$. Then the 2-natural
isomorphism $\mathcal{M}$ of Theorem~\ref{preEMthm} satisfies the
commuting diagram
\[
 \xy
 (-20,20)*++{\cat{$\T$-Alg}}="1";
 (20,20)*++{\cat{$\G$-CoAlg}}="2";
 (0,5)*++{\mathcal{K}(-,B)}="3";
    {\ar^-{\mathcal{M}} "1";"2" };
    {\ar@<-1ex>_{U^{\T}} "1";"3" };
    {\ar@<-1ex>_{F^{\T}} "3";"1" };
    {\ar@<1ex>^{U^{\G}} "2";"3" };
    {\ar@<1ex>^{F^{\G}} "3";"2" };
 \endxy
\]
\end{lem}

\proof By Theorem~\ref{preEMthm} all we must show is that
$\mathcal{M}F^{\T}=F^{\G}$.  This equality can be verified
pointwise. Let $s \maps A \to B$ so that
 \ban
 \mathcal{M}_AF_A^{\T}(s) & =& \left(Ts,T\mu s.\iota s\right) \\
 F^{\G}_A(s) & = & \left(Ts,\delta s\right) =
    \left(Ts,T^2(\varepsilon.\mu)s.T^2\mu Ts.T\iota T^2s.\iota Ts
    \right).
 \ean
The required equality follows by the commutativity of the
following diagram:
\[
 \xy
   (-60,0)*+{Ts}="1";
   (-40,10)*+{T^3s}="2t";
   (-40,-10)*+{T^3 s }="2b";
   (-20,10)*+{T^2s}="3t";
   (-20,-10)*+{T^5 s }="3b";
    (0,0)*+{T^4s}="4t";
   (0,-20)*+{T^4 s }="4b";
   (40,10)*+{T^2s}="5t";
   (20,-10)*+{T^3 s }="5b";
        {\ar@{=>}^{\iota Ts  } "1";"2t"};
        {\ar@{=>}_{\iota Ts  } "1";"2b"};
        {\ar@{=}_{} "2b";"2t"};
        {\ar@{=>}^{T \mu s  } "2t";"3t"};
        {\ar@{=>}_{T \iota T^2s } "2b";"3b"};
        {\ar@{=>}^{T \iota Ts } "3t";"4t"};
        {\ar@{=>}^{T^3 \mu Ts } "3b";"4t"};
        {\ar@{=>}_{T^2 \mu Ts } "3b";"4b"};
        {\ar@{=>}^{T^2s} "3t";"5t"};
        {\ar@{=>}^{T^2\mu s } "4t";"5b"};
        {\ar@{=>}_{T^2\mu s } "4b";"5b"};
        {\ar@{=>}_{T^2\varepsilon s } "5b";"5t"};
 \endxy
\]
together with the fact that if $h\maps s \To s'$, then $F^{\T}(h)=
Th = F^{\G}(h)$ and the fact that $\mathcal{M}$ is the identity on
morphisms. \qed

Compare the following two Theorems to Proposition 1.4 and
Proposition 1.5 of~\cite{StreetFrob}.

\begin{thm} \label{MainThm}
Given a Frobenius monad $(\T, \mu, \eta,\varepsilon)$ on an object
$B$ in $\mathcal{K}$, then in \cat{EM}$(\mathcal{K})$ the left
adjoint $F^{\T}\maps B \to B^{\T}$ to the forgetful functor
$U^{\T}\maps B^{\T} \to B$ is also  right adjoint to $U^{\T} $with
counit $\varepsilon$.  Hence, the Frobenius monad $\T$ is
generated by an ambidextrous adjunction in
\cat{EM}$(\mathcal{K})$.
\end{thm}

\Proof Identify $\T$ with its fully faithful image via the
2-functor $Z \maps \mathcal{K} \to \cat{EM}(\mathcal{K})$; then an
Eilenberg-Moore object $B^{\T}$ for the  monad $\T$ exists in
\cat{EM}$(\mathcal{K})$.  Let $\G$ denote the induced comonad
structure on $T$ given in Proposition~\ref{defadjmonad}.  Then by
the adjoint monad theorem, the object $B^{\T}$ serves as an
Eilenberg-Moore object $B^{\G}$ for the comonad $\G$ and we have
the commutative diagram:
\[
 \xy
 (-20,20)*+{B^{\T}}="1";
 (20,20)*+{B^{\G}}="2";
 (0,5)*+{B}="3";
    {\ar^-{\mathcal{M}} "1";"2" };
    {\ar_{U^{\T}} "1";"3" };
    {\ar^{U^{\G}} "2";"3" };
 \endxy
\]
By the remarks following Definition~\ref{defEM} we have
\[
i^{\T},e^{\T}\maps F^{\T} \dashv U^{\T} \maps B \to B^{\T}
\]
and
\[
 i^{\G},e^{\G}\maps F^{\G} \vdash U^{\G}\maps
B \to B^{\G}.
\]
Since $U^{\G}F^{\G}$ generates the comonad $\G$ and $\varepsilon$
is the counit for the comonad $\G$, it is clear that
$e^{\G}=\varepsilon$ above.  All that remains to be shown is  that
$\mathcal{M}F^{\T}=F^{\G}$. This follows by the 2-categorical
Yoneda lemma applied to Lemma~\ref{lemFM}.  Hence, the Frobenius
monad $\T = U^{\T}F^{\T}$ is generated by an ambijunction
$F^{\T}\dashv U^{\T} \dashv F^{\T}$ in \cat{EM}$(\mathcal{K})$.
\qed

\begin{thm} \label{Ambi_to_Frob}
Let $i,e,j,k\maps F \dashv U \dashv F\maps B \to C$ be an
ambidextrous adjunction in the 2-category $\mathcal{K}$. Then the
monad $(UF,UiF,e)$ generated by the adjunction is a Frobenius
monad with $\varepsilon = k$.
\end{thm}

\Proof All we must show is that $UF \dashv UF$ with counit
$k.UiF$.  Define the unit of the adjunction to be $UjF.i$. The
zig-zag identities follow from the zig-zag identities for $(i,e)$
and $(j,k)$. \qed

\begin{cor}
If $\ambijunction{B}{C}{F}{U}$ is a specified ambijunction in the
2-category $\mathcal{K}$, then $UF$ is a Frobenius object in the
strict monoidal category $\mathcal{K}(B,B)$.
\end{cor}

\Proof  By Theorem~\ref{Ambi_to_Frob}, $UF$ defines a Frobenius
monad on the object $B$ in $\mathcal{K}$.  As explained above,
this is simply a Frobenius object in the monoidal category
$\mathcal{K}(B,B)$. \qed

\begin{cor} \label{FrobjectMon}
A Frobenius object in a monoidal category $M$ yields an
ambijunction in \cat{EM}$(\Sigma(M))$, where $\Sigma(M)$ is the
2-category obtained by the strictification of the suspension of
$M$.
\end{cor}

\Proof  Recall that a monad on an object $B$ in a 2-category
$\mathcal{K}$ can be thought of as a monoid object in the monoidal
category $\mathcal{K}(B,B)$.  Similarly, a comonad on $B$ is just
a comonoid object in $\mathcal{K}(B,B)$.  Regarding $M$ as a one
object 2-category $\Sigma(M)$, a Frobenius object in $M$ is simply
a Frobenius monad in $\Sigma(M)$. Applying Theorem~\ref{MainThm}
completes the proof. \qed

\begin{cor} \label{xxcc}
Every Frobenius algebra in the category $\cat{Vect}$ arises from
an ambidextrous adjunction in the 2-category whose objects are
algebras, morphisms are bimodules of algebras, and whose
2-morphisms are bimodule homomorphisms.
\end{cor}

\Proof This follows immediately from Corollary \ref{FrobjectMon}
and the discussion in subsection \ref{EMVECT}. \qed

\begin{cor}
Every 2D topological quantum field theory, in the sense of
Atiyah~\cite{Atiyah}, arises from an ambidextrous adjunction in
the 2-category whose objects are algebras, morphisms are bimodules
of algebras, and whose 2-morphisms are bimodule homomorphisms.
\end{cor}

\Proof  Since a 2D topological quantum field theory is equivalent
to a commutative Frobenius algebra~\cite{Abrams1,Kock}, the proof
follows from Corollary \ref{xxcc}. \qed

\section{Categorification}

In this section we extend the theory of the previous section to
the context of $\cat{Gray}$-categories.  $\cat{Gray}$ is the
symmetric monoidal closed category whose underlying category is
\cat{2-Cat}; the category whose objects are 2-categories, and
whose morphisms are 2-functors. \cat{Gray} differs from
\cat{2-Cat} in that \cat{Gray} has a more interesting monoidal
structure than the usual cartesian monoidal structure on
\cat{2-Cat}.  A $\cat{Gray}$-category, also known as a semistrict
3-category, is defined using enriched category theory~\cite{kel}
as a category enriched in $\cat{Gray}$.  The unusual tensor
product in \cat{Gray}, or `Gray tensor product', has the effect of
equipping a \cat{Gray}-category $\mathcal{K}$ with a cubical
functor $M \maps \mathcal{K}(B,C) \times \mathcal{K}(A,B) \to
\mathcal{K}(A,C)$ for all objects $A$,$B$,$C$ in $\mathcal{K}$.
This means that if $f \maps  F \To F'$ in $\mathcal{K}(A,B)$, and
$g \maps G \To G'$ in $\mathcal{K}(B,C)$, then, rather than
commuting on the nose, we have an invertible 3-cell $M_{g,f}$ ---
denoted $g_f$ following Marmolejo~\cite{mar}
--- in the following square:
\[
 \xy
    (-10,8)*+{GF}="TL";
    (10,8)*+{G'F}="TR";
    (10,-8)*+{G'F'.}="BR";
    (-10,-8)*+{GF'}="BL";
    {\ar@{=>}^{gF} "TL";"TR"};
    {\ar@{=>}^{G'f} "TR";"BR"};
    {\ar@{=>}_{Gf} "TL";"BL"};
    {\ar@{=>}_{gF'} "BL";"BR"};
    (0,3)*{}="x";
    (0,-3)*{}="y";
    {\ar@3{->}^{g_f} "x";"y"};
 \endxy
\]
We take this notion to be a sufficiently general extension since
every tricategory or weak 3-category is triequivalent to a
$\cat{Gray}$-category~\cite{GPS}.

The proof of the adjoint monad theorem relied heavily on the
notion of mates under adjunction and the fact that this
relationship respected composites of morphisms and adjunctions. In
order to categorify this theorem we will first have to categorify
the notion of mates under adjunction to the notion of mates under
pseudoadjunction. In this case, rather than a bijection between
certain morphisms, we will have an equivalence of $\Hom$
categories.  The naturality of this equivalence will also be
discussed.

In Section~\ref{secpseudomonads} we define the notion of a
pseudomonad in a $\cat{Gray}$-category and review some of the
basic theory.  Using the notion of mateship under pseudoadjunction
it is shown that if a pseudomonad has a specified pseudoadjoint
$G$, then $G$ is a pseudocomonad.  All of the theorems from the
previous section are then extended into this context and the
notion of a Frobenius pseudomonad and Frobenius pseudomonoid are
given.  The main result that every Frobenius pseudomonoid arises
from a pseudo ambijunction is then proven as a corollary of the
categorified version of the Eilenberg-Moore adjoint monad theorem
in Section~\ref{pseudoadjsec}.

\subsection{Pseudoadjunctions}

We begin with the definition of a pseudoadjunction given by Verity
in~\cite{Verity} where they were called \textit{locally-adjoint
biadjoint pairs}.  For more details see also the discussion by
Lack where the `free living' or `walking' pseudoadjunction is
defined~\cite{Lack}.

\begin{defn}
A {\em pseudoadjunction} $I,E,i,e \maps F \dashv_p U\maps A \to B$
in a $\cat{Gray}$-category $\mathcal{K}$ consists of:
\begin{itemize}
   \item morphisms $U \maps A \to B$ and $F \maps B \to A$,
   \item 2-morphisms $i \maps 1 \To UF$ and $e \maps FU \To 1$,  and
   \item coherence 3-isomorphisms
   \[ \vcenter{
 \xy
   (-15,0)*+{U}="l";
   (0,15)*+{UFU}="t";
   (15,0)*+{U}="r";
    {\ar@{=>}^{iU} "l";"t"};
    {\ar@{=>}^{Ue} "t";"r"};
    {\ar@{=>}_1 "l";"r"};
        {\ar@3{->}^{I} (0,3);"t"+(0,-6)};
 \endxy}
\qquad {\rm and} \qquad
 \vcenter{\xy
   (-15,0)*+{F}="l";
   (0,15)*+{FUF}="t";
   (15,0)*+{F}="r";
    {\ar@{=>}^{Fi} "l";"t"};
    {\ar@{=>}^{eF} "t";"r"};
    {\ar@{=>}_1 "l";"r"};
    {\ar@3{->}_{E} "t"+(0,-6);(0,3)};
 \endxy}
\]
\end{itemize}
such that the following two diagrams are both identities:
\[
 \xy
   (-38,0)*+{FU}="ll";
   (0,14)*+{FU}="t";
   (-14,0)*+{FUFU}="l";
   (14,0)*+{1}="r";
   (0,-14)*+{FU}="b";
        {\ar@{=>}_{FUe} "l";"t"};
        {\ar@{=>}^{eFU} "l";"b"};
        {\ar@{=>}^{e} "t";"r"};
        {\ar@{=>}_{e} "b";"r"};
        {\ar@{=>}^{FiU} "ll";"l"};
        {\ar@/^1pc/@{=>}^1 "ll";"t"};
        {\ar@/_1pc/@{=>}_1 "ll";"b"};
        {\ar@3{->}^{e_e^{-1}} "t"+(2,-11);"b"+(2,11)};
        {\ar@3{->}_{FI} (-16,9);(-16,4)};
        {\ar@3{->}_{EU} (-16,-4);(-16,-9)};
 \endxy
\qquad \qquad
 \xy
   (38,0)*+{UF}="ll";
   (0,14)*+{UF}="t";
   (-14,0)*+{1}="l";
   (14,0)*+{UFUF}="r";
   (0,-14)*+{UF}="b";
        {\ar@{=>}^{i} "l";"t"};
        {\ar@{=>}_{i} "l";"b"};
        {\ar@{=>}_{iUF} "t";"r"};
        {\ar@{=>}^{UFi} "b";"r"};
        {\ar@{=>}^{UeF} "r";"ll"};
        {\ar@/^1pc/@{=>}^1 "t";"ll"};
        {\ar@/_1pc/@{=>}_1 "b";"ll"};
        {\ar@3{->}_{i_i^{-1}} "t"+(-2,-11);"b"+(-2,11)};
        {\ar@3{->}^{IF} (16,9);(16,4)};
        {\ar@3{->}^{UE} (16,-4);(16,-9)};
 \endxy
\]
\end{defn}

We will sometimes denote a pseudoadjunction as $F \dashv_p U$ and
say that the morphism $U$ is the {\em right pseudoadjoint} of $F$.
Likewise, $F$ is said to be the {\em left pseudoadjoint} of $U$.

\begin{prop} \label{defcomppseudoadj}
If $I,E,i,e\maps F \dashv_p U\maps A \to B$ and $I',E',i',e' \maps
F' \dashv_p U'\maps B \to C$, then $FF' \dashv_p U'U$ with
 \ban
 \bar{i} &:=& \xymatrix@C=2.2pc{
 1 \ar@{=>}[r]^-{i'} & U'F' \ar@{=>}[r]^-{U'iF'} & U'UFF'
 } \\
 \bar{e} &:=& \xymatrix@C=2.2pc{
 FF'U'U \ar@{=>}[r]^-{Fe'U} & FU \ar@{=>}[r]^-{e} & 1
 }
 \ean
 and
 \[
\bar{I} :=
 \xy
   (-40,-15)*+{U'U}="bl";
   (0,-15)*+{U'U}="bm";
   (40,-15)*+{U'U}="br";
   (-20,0)*+{U'F'U'U}="l";
   (20,0)*+{U'UFU}="r";
   (0,15)*+{U'UFF'U'U}="t";
        {\ar@{=>}_1 "bl";"bm"};
        {\ar@{=>}_1 "bm";"br"};
        {\ar@{=>}^{i'U'U} "bl";"l"};
        {\ar@{=>}^<<<<<<<{U'e'U} "l";"bm"};
        {\ar@{=>}^>>>>>>>{U'Ui} "bm";"r"};
        {\ar@{=>}^{U'eU} "r";"br"};
        {\ar@{=>}^{U'iF'U'U} "l";"t"};
        {\ar@{=>}^{U'UFe'U} "t";"r"};
   {\ar@3{->}^{I'U} (-20,-11);(-20,-5)};
   {\ar@3{->}^{U'I} (20,-11);(20,-5)};
   {\ar@3{->}^{U'i_{e'}U} (0,-2);(0,8)};
 \endxy
\]

\[
\bar{E} :=
 \xy
   (-40,-15)*+{FF'}="bl";
   (0,-15)*+{FF'}="bm";
   (40,-15)*+{FF'.}="br";
   (-20,0)*+{FF'U'F'}="l";
   (20,0)*+{FUFF'}="r";
   (0,15)*+{FF'U'UFF'}="t";
        {\ar@{=>}_1 "bl";"bm"};
        {\ar@{=>}_1 "bm";"br"};
        {\ar@{=>}^{FF'i'} "bl";"l"};
        {\ar@{=>}^<<<<<<<{Fe'F'} "l";"bm"};
        {\ar@{=>}^>>>>>>>{iFF'} "bm";"r"};
        {\ar@{=>}^{FeF'} "r";"br"};
        {\ar@{=>}^{FF'U'iF'} "l";"t"};
        {\ar@{=>}^{Fe'UFF'} "t";"r"};
   {\ar@3{->}^{FE'} (-20,-5);(-20,-11)};
   {\ar@3{->}^{EF'} (20,-5);(20,-11)};
   {\ar@3{->}^{Fi_{e'}^{-1}F'} (0,8);(0,-2)};
 \endxy
\]
\end{prop}

\Proof  The proof is given in~\cite{gray} although it is a routine
verification and can be checked directly. \qed

\begin{prop} \label{pseudomates}
Let
 \begin{itemize}
  \item $I,E,i,e \maps F \dashv_p U\maps A \to B$,  and
  \item $I',E',i',e' \maps
F' \dashv_p U' \maps A' \to B'$
 \end{itemize}
 in the $\cat{Gray}$-category
$\mathcal{K}$. If $a \maps A \to A'$ and $b \maps B \to B'$, then
there is an equivalence of categories $\mathcal{K}(bU,U'a) \simeq
\mathcal{K}(F'b,aF)$ given by:
 \ban
 \Theta\maps \mathcal{K}(bU,U'a) &\to& \mathcal{K}(F'b,aF) \\
 \xi &\mapsto& \zeta =
 \xymatrix@C=2.2pc{F'b \ar@{=>}[r]^-{F'bi} & F'bUF
    \ar@{=>}[r]^-{F'\xi F} & F'U'aF \ar@{=>}[r]^-{e'aF} & aF } \\
 \omega \maps \xi_1 \Rrightarrow \xi_2 &\mapsto&
 \xymatrix@C=2.8pc{F'b \ar@{=>}[r]^-{F'bi} & F'bUF
    \ar@/^1.2pc/@{=>}[r]^-{F'\xi_1 F}_>>>>>{}="1" \ar@/_1.2pc/@{=>}[r]_-{F'\xi_2 F}^>>>>>{}="2"
    & \ar@3{->}_{F'\omega F}"1";"2" F'U'aF \ar@{=>}[r]^-{e'aF} & aF }
 \ean
and
 \ban
 \Phi\maps  \mathcal{K}(F'b,aF)&\to&  \mathcal{K}(bU,U'a)\\
 \zeta &\mapsto& \xi =
\xymatrix@C=2.2pc{bU \ar@{=>}[r]^-{i'bU} & U'F'bU
\ar@{=>}[r]^-{U'\zeta U} & U'aFU \ar@{=>}[r]^-{U'ae} & U'a } \\
 \varrho\maps \zeta_1 \Rrightarrow \zeta_2 &\mapsto&
\xymatrix@C=2.8pc{bU \ar@{=>}[r]^-{i'bU} & U'F'bU
\ar@/^1.2pc/@{=>}[r]^-{U'\zeta_1 U}_>>>>>{}="1"
\ar@/_1.2pc/@{=>}[r]_-{U'\zeta_2 U}^>>>>>{}="2"
    & \ar@3{->}_{U'\varrho U}"1";"2" U'aFU \ar@{=>}[r]^-{U'ae} & U'a }.
 \ean
\end{prop}

\Proof It is clear that $\Theta$ is a functor from its definition
above. That is, $\Theta$ preserves composites of 3-morphisms along
2-morphisms in $\mathcal{K}$. Let $\xi$ be an object of
$\mathcal{K}(bU,U'a)$ so that $\Phi\Theta(\xi)$ is given by the
composite
\[
\xymatrix@C=2.3pc{bU \ar@{=>}[r]^-{i'bU} & U'F'bU
\ar@{=>}[r]^-{U'F'biU} & U'F'bUFU
    \ar@{=>}[r]^-{U'F'\xi FU} & U'F'U'aFU \ar@{=>}[r]^-{U'e'aFU} &  U'aFU \ar@{=>}[r]^-{U'ae} & U'a
    }.
\]
Define an isomorphism $\gamma_{\xi} \maps \xi
\Rrightarrow\Phi\Theta(\xi)$ by the diagram

\[
 \xy
   (-60,0)*+{bU}="1";
   (-40,-0)*+{U'F'bU}="2";
   (-30,-20)*+{U'F'bUFU}="3";
   (5,-20)*+{U'F'U'aFU}="4";
   (40,-20)*+{U'aFU}="5";
   (60,0)*+{U'a}="6";
   (-10,-0)*+{U'F'bU}="a";
   (20,-0)*+{U'F'U'a}="b";
   (0,20)*+{U'a}="z";
        {\ar@{=>}_-{i'bU} "1";"2"};
        {\ar@{=>}_{U'F'biU} "2";"3"};
        {\ar@{=>}_-{U'F'\xi FU} "3";"4"};
        {\ar@{=>}_-{U'e'aFU} "4";"5"};
        {\ar@{=>}_-{U'ae} "5";"6"};
        {\ar@{=>}^-{\xi} "1";"z"};
        {\ar@{=>}^-{1} "z";"6"};
        {\ar@{=>}^-{U'F'\xi} "a";"b"};
        {\ar@{=>}^-{1} "2";"a"};
        {\ar@{=>}_-{U'F'bUe} "3";"a"};
        {\ar@{=>}_-{U'F'U'ae} "4";"b"};
        {\ar@{=>}_-{U'e'a} "b";"6"};
        {\ar@{=>}_-{i'U'a} "z";"b"};
            {\ar@3{->}^{U'e'^{-1}_e} (30,-7);(30,-13)};
            {\ar@3{->}^{U'F'bE^{-1}} (-31,-3);(-31,-9)};
            {\ar@3{->}^{U'F'\xi_e^{-1}} (-2,-5);(-2,-11)};
            {\ar@3{->}_{i'^{-1}_{\xi}} (-10,11);(-10,5)};
            {\ar@3{->}^<<{I'a} (23,9);(21,4)};
 \endxy
\]
which is invertible because $I'$, $E$ and the structural maps in
the $\cat{Gray}$-category are invertible.  It is straight forward
to check the naturality of this isomorphism.  Let $\omega \maps
\xi_1 \Rrightarrow\xi_2$; then $\gamma_{\xi_2}\circ\omega =
\Phi\Theta(\omega)\circ\gamma_{\xi_1}$ by the invertibility of
$I'$, $E$ and the axioms of a $\cat{Gray}$-category. The
isomorphism $\bar{\gamma}_{\zeta} \maps\zeta
\Rrightarrow\Theta\Phi(\zeta)$, for $\zeta$ in
$\mathcal{K}(F'b,aF)$, is given by:
\[
 \xy
   (-60,0)*+{F'b}="1";
   (-40,-0)*+{F'bUF}="2";
   (-30,-20)*+{F'U'F'bUF}="3";
   (5,-20)*+{F'U'aFUF}="4";
   (40,-20)*+{F'U'aF}="5";
   (60,0)*+{aF}="6";
   (-10,-0)*+{F'bUF}="a";
   (20,-0)*+{aFUF}="b";
   (0,20)*+{aF}="z";
        {\ar@{=>}_-{F'bi} "1";"2"};
        {\ar@{=>}_{F'i'bUF} "2";"3"};
        {\ar@{=>}_-{F'U'\zeta UF} "3";"4"};
        {\ar@{=>}_-{F'U'aeF} "4";"5"};
        {\ar@{=>}_-{e'aF} "5";"6"};
        {\ar@{=>}^-{\zeta} "1";"z"};
        {\ar@{=>}^-{1} "z";"6"};
        {\ar@{=>}^-{\zeta UF} "a";"b"};
        {\ar@{=>}^-{1} "2";"a"};
        {\ar@{=>}_-{e'F'bUF} "3";"a"};
        {\ar@{=>}_-{e'aFUF} "4";"b"};
        {\ar@{=>}_-{aeF} "b";"6"};
        {\ar@{=>}_-{aiF} "z";"b"};
            {\ar@3{->}^{e'_{ae}F} (30,-7);(30,-13)};
            {\ar@3{->}^{E'^{-1}bUF} (-31,-3);(-31,-9)};
            {\ar@3{->}^{e'_{\zeta}UF} (-2,-5);(-2,-11)};
            {\ar@3{->}_{\zeta_i} (-10,11);(-10,5)};
            {\ar@3{->}^<<{aI} (23,9);(21,4)};
 \endxy
\]
By similar arguments as above this isomorphism is natural. \qed

Using this equivalence of categories we extend the notion of
mateship under adjunction to the notion of {\em mateship under
pseudoadjunction}.  We now express the naturality conditions this
equivalence satisfies:

\begin{prop} \label{comppsuedomatesONEcell}
Consider the collection of pseudoadjunctions and morphisms:
 \begin{itemize}
 \item $I,E,i,e\maps F \dashv_p U\maps A \to B$,
 \item $I',E',i',e' \maps F' \dashv_p U' \maps A' \to B'$,
 \item $I'',E'',i'',e'' \maps F'' \dashv_p U''\maps A'' \to B''$,
 and
 \item $a\maps A \to A', a' \maps A' \to A'' ,b\maps
B \to B', b' \maps B' \to B''$,
 \end{itemize}
in the $\cat{Gray}$-category $\mathcal{K}$.  Let
 \ban
 \Theta \maps \mathcal{K}(bU,U'a) \to
\mathcal{K}(F'b,aF)
 & \quad &
\Phi \maps \mathcal{K}(F'b,aF) \to \mathcal{K}(bU,U'a)                      \\
 \Theta' \maps \mathcal{K}(b'U',U''a')
\to \mathcal{K}(F''b',a'F')
 & \quad &
 \Phi' \maps \mathcal{K}(F''b',a'F')
\to  \mathcal{K}(b'U',U''a')        \\
 \bar{\Theta}\maps\mathcal{K}(b'bU,U''a'a) \to
\mathcal{K}(F''b'b,a'aF)
 & \quad &
 \bar{\Phi}\maps \mathcal{K}(F''b'b,a'aF)\to
\mathcal{K}(b'bU,U''a'a)
 \ean
be the functors from Proposition~\ref{pseudomates} defining the
relevant equivalences of categories.  Then there exists a natural
isomorphism $W$ between the following pasting composites of
functors:
 \ban
 a'\Theta (-).\Theta'(-)b &\maps& \mathcal{K}(bU,U'a) \times
 \mathcal{K}(b'U',U''a') \to \mathcal{K}(F''b'b,a'aF)        \\
 \bar{\Theta}(-a.b'-) &:& \mathcal{K}(bU,U'a) \times
 \mathcal{K}(b'U',U''a') \to \mathcal{K}(F''b'b,a'aF),
 \ean
and a natural isomorphism $Y$ between the pasting composites:
 \ban
 \Phi'(-)a.b'\Phi(-) &\maps& \mathcal{K}(F'b,aF) \times
 \mathcal{K}(F''b',a'F') \to \mathcal{K}(b'bU,U''a'a) \\
 \bar{\Phi}(a'-.-b) &\maps& \mathcal{K}(F'b,aF) \times
 \mathcal{K}(F''b',a'F') \to \mathcal{K}(b'bU,U''a'a) .
 \ean
\end{prop}

\Proof  Let $\xi\in\mathcal{K}(bU,U'a)$ and
$\xi'\in\mathcal{K}(b'U',U''a')$, then $W(\xi \times \xi')$ is
given by the following pasting composite of invertible
3-morphisms:
\[
\def\objectstyle{\scriptstyle}
\def\labelstyle{\scriptstyle}
 \xy
   (-75,10)*+{F''b'b}="t1";
   (-60,20)*+{F''b'U'F'b}="t2";
   (-30,30)*+{F''U''a'F'b}="t3";
   (0,30)*+{a'F'b}="t4";
   (30,30)*+{a'F'bUF}="t5";
   (60,20)*+{a'F'U'aF}="t6";
   (75,10)*+{a'aF}="t7";
   (-50,10)*+{F''b'U'F'bUF}="m2";
   (0,20)*+{F''U''a'F'bUF}="m3";
   (50,10)*+{F''U''a'F'U'aF}="m4";
   (-60,-20)*+{F''b'bUF}="b2";
   (-30,-20)*+{F''b'U'aF}="b3";
   (60,-20)*+{F''U''a'aF}="b4";
   (0,5)*+{F''b'U'F'U'aF}="l2";
   (15,-10)*+{F''b'U'aF}="l3";
        {\ar@{=>}^{F''b'i'b} "t1";"t2"};
        {\ar@{=>}^{F''\xi'F'b} "t2";"t3"};
        {\ar@{=>}^{e''a'F'b} "t3";"t4"};
        {\ar@{=>}^{a'F'bi} "t4";"t5"};
        {\ar@{=>}^{a'F' \xi F} "t5";"t6"};
        {\ar@{=>}^{a'e'aF} "t6";"t7"};
        {\ar@{=>}_{F''b'bi} "t1";"b2"};
        {\ar@{=>}_{F''b' \xi F} "b2";"b3"};
        {\ar@{=>}_{F'' \xi'aF} "b3";"b4"};
        {\ar@{=>}_{e''a'aF} "b4";"t7"};
        {\ar@{=>}^{F''b'U'F'bi} "t2";"m2"};
        {\ar@{=>}_{F''b'i'bUF} "b2";"m2"};
        {\ar@{=>}_{F''U''a'F'bi\;\;\;\;} "t3";"m3"};
        {\ar@{=>}^{F''\xi'F'bUF} "m2";"m3"};
        {\ar@{=>}_{\;\;\;\;\;e''a'F'bUF} "m3";"t5"};
        {\ar@{=>}^{F''U''a'F'\xi F} "m3";"m4"};
        {\ar@{=>}_{F''U''a'e'aF} "m4";"b4"};
        {\ar@{=>}^{e''a'F'U'aF} "m4";"t6"};
        {\ar@{=>}^{F''b'U'F'\xi F} "m2";"l2"};
        {\ar@{=>}^{F''b'i'U'aF} "b3";"l2"};
        {\ar@{=>}^{F''\xi'F'U'aF} "l2";"m4"};
        {\ar@{=>}^{F''bU'e'aF} "l2";"l3"};
        {\ar@{=>}^{F''\xi'aF} "l3";"b4"};
        {\ar@{=>}^{1} "b3";"l3"};
            {\ar@3{->}_{F''b'i'_{bi}} (-60,7);(-60,2)};
            {\ar@3{->}^{F''b'i'_{\xi}F} (-35,3);(-35,-2)};
            {\ar@3{->}_{F''b'I'aF} (0,-5);(0,-10)};
            {\ar@{=} (15,-13);(15,-16)};
            {\ar@3{->}^{F''\xi'_{e'}aF} (30,3);(30,-2)};
            {\ar@3{->}_{F''\xi'_{F'\xi}F} (0,15);(0,10)};
            {\ar@3{->}_<<{e''_{a'F'bi}} (0,28);(0,23)};
            {\ar@3{->}_{F''\xi'_{F'bi}} (-32,25);(-32,20)};
            {\ar@3{->}^{e''_{a'F'\xi}F} (31,25);(31,20)};
            {\ar@3{->}^{e''_{ae'}aF} (60,7);(60,2)};
 \endxy
\]
If $\omega\maps \xi_1 \TO \xi_2$ and $\omega'\maps \xi'_1 \TO
\xi'_2$ then the naturality of $W$ follows from the axioms of the
cubical functor defining the Gray tensor product. Given $\zeta \in
\mathcal{K}(F'b,aF)$ and $\zeta' \in \mathcal{K}(F''b',a'F')$ the
natural isomorphism $Y$ can be defined similarly:
\[
\def\objectstyle{\scriptstyle}
\def\labelstyle{\scriptstyle}
 \xy
   (-75,8)*+{b'bU}="t1";
   (-60,23)*+{b'U'F'bU}="t2";
   (-30,30)*+{b'U'aFU}="t3";
   (0,30)*+{b'U'a}="t4";
   (30,30)*+{U''F''b'U'a}="t5";
   (60,23)*+{U''a'F'U'a}="t6";
   (75,8)*+{U''a'a}="t7";
   (-50,10)*+{U''F''b'U'F'bU}="m2";
   (0,20)*+{U''F''b'U'aFU}="m3";
   (50,10)*+{U''a'F'U'aFU}="m4";
   (-65,-20)*+{U''F''b'bU}="u2";
   (-30,-20)*+{U''a'F'bU}="u3";
   (65,-20)*+{U''a'aFU}="u4";
   (0,5)*+{U''a'F'U'F'bU}="l2";
   (15,-10)*+{U''a'F'bU}="l3";
        {\ar@{=>}^{b'i'bU} "t1";"t2"};
        {\ar@{=>}^{b'U'\zeta U} "t2";"t3"};
        {\ar@{=>}^{b'U'ae} "t3";"t4"};
        {\ar@{=>}^{i''b'U'a} "t4";"t5"};
        {\ar@{=>}^{U''\zeta'U'a} "t5";"t6"};
        {\ar@{=>}^{U''a'e'a} "t6";"t7"};
        {\ar@{=>}_{i''b'bU} "t1";"u2"};
        {\ar@{=>}_{U''\zeta'bU} "u2";"u3"};
        {\ar@{=>}_{U''a'\zeta U} "u3";"u4"};
        {\ar@{=>}_{U''a'ae} "u4";"t7"};
        {\ar@{=>}^{i''b'U'F'bU} "t2";"m2"};
        {\ar@{=>}_{U''F''i'bU} "u2";"m2"};
        {\ar@{=>}_{i''b'U'aFU\;\;\;\;} "t3";"m3"};
        {\ar@{=>}^{U''F''b'U'\zeta U} "m2";"m3"};
        {\ar@{=>}_{\;\;\;\;\;U''F''b'U'ae} "m3";"t5"};
        {\ar@{=>}^{U''\zeta' F'U'aFU} "m3";"m4"};
        {\ar@{=>}_{U''a'e'aFU} "m4";"u4"};
        {\ar@{=>}^{U'a'F'U'ae} "m4";"t6"};
        {\ar@{=>}^{U''\zeta'U'F'bU} "m2";"l2"};
        {\ar@{=>}^{U''a'F'i'bU} "u3";"l2"};
        {\ar@{=>}^{U''a'F'U'\zeta'U} "l2";"m4"};
        {\ar@{=>}^{U''a'e'F'bU} "l2";"l3"};
        {\ar@{=>}^{U''a'\zeta U} "l3";"u4"};
        {\ar@{=>}^{1} "u3";"l3"};
            {\ar@3{->}_{(i''b'_{i'})^{-1}bU} (-55,7);(-55,2)};
            {\ar@3{->}^{U''(\zeta'_{i'})^{-1}} (-35,3);(-35,-2)};
            {\ar@3{->}_{U''a'E'^{-1}bU} (5,-5);(5,-10)};
            {\ar@{=} (15,-13);(15,-16)};
            {\ar@3{->}^{U''a'(e'_{F'U'\zeta'})^{-1}U} (27,3);(27,-2)};
            {\ar@3{->}_{U''(\zeta'_{U'\zeta})^{-1}U} (5,15);(5,10)};
            {\ar@3{->}_<<{(i''_{b'U'aFUe})^{-1}} (5,28);(5,23)};
            {\ar@3{->}_{(i''_{b'U'\zeta})^{-1}U} (-32,25);(-32,20)};
            {\ar@3{->}^{U''(\zeta'_{U'ae})^{-1}} (31,25);(31,20)};
            {\ar@3{->}^{U''a'(e'_{ae})^{-1}} (55,7);(55,2)};
 \endxy
\]
and by similar arguments, $Y$ can be shown to be natural. \qed

We will find it necessary later in this paper to refer to the
natural isomorphisms $W$ and $Y$ within the context of some
specified choices of composable pseudoadjunctions and morphisms
$a,a',b, b'$.  We will refer to these isomorphisms generically as
$W$ and $Y$, even though we will consider many different choices
of pseudoadjunctions and morphisms.  This is possible since these
natural isomorphisms exist for every possible choice of this data.
Furthermore, the specific pseudoadjoints and morphisms should be
clear from the context.  When no confusion is likely to arise we
will also denote $\Theta',\bar{\Theta},\Phi',$ and $\bar{\Phi}$
simply as $\Theta$ and $\Phi$, respectively. Note in particular
that when $a=1_A, b=1_B, a'=1_A, b'=1_B$, then we have
$\Theta(\xi).\Theta(\xi') \cong \Theta(\xi'.\xi)$, and similarly
for $\Phi$.

\begin{prop} \label{comppseudomatesPseudoadj}
Consider the collection of pseudoadjunctions and morphisms:
\begin{itemize}
  \item $I,E,i,e \maps F \dashv_p U\maps A \to B$,
  \item $I',E',i',e'\maps F' \dashv_p U' \maps A' \to B'$,
  \item $I_1,E_1,i_1,e_1\maps F_1 \dashv_p U_1 \maps B\to C$,
  \item  $I'_1,E'_1,i'_1,e'_1 \maps F'_1 \dashv_p U'_1 \maps B' \to C'$,
  and
  \item $a\maps A \to A',b\maps B \to
B',c\maps C \to C'$
\end{itemize}
in the $\cat{Gray}$-category $\mathcal{K}$. Let
 \ban
 \Theta \maps \mathcal{K}(bU,U'a)\to
\mathcal{K}(F'b,aF)
 & \quad &
 \Phi \maps \mathcal{K}(F'b,aF)\to \mathcal{K}(bU,U'a)                 \\
 \Theta_1 \maps \mathcal{K}(cU_1,U'_1b)\to
\mathcal{K}(F'_1c,bF_1) & \quad &
 \Phi_1 \maps \mathcal{K}(F'_1c,bF_1)\to\mathcal{K}(cU_1,U'_1b)
            \\
 \bar{\Theta}\maps\mathcal{K}(cU_1U,U'_1U'a)\to
\mathcal{K}(F'F'_1c,aFF_1)
 & \quad &
 \bar{\Phi}\maps\mathcal{K}(F'F'_1c,aFF_1)\to
\mathcal{K}(cU_1U,U'_1U'a)
 \ean
be the functors from Proposition~\ref{pseudomates} defining the
relevant equivalences of categories.  Here $\bar{\Theta}$ and
$\bar{\Phi}$ are the equivalence corresponding to the composite
pseudoadjunction defined in Proposition~\ref{defcomppseudoadj}.
Then there exists a natural isomorphism $V$ between the following
pasting composites of functors:
 \ban
\Theta(-)F_1.F'\Theta_1(-) &\maps& \mathcal{K}(bU,U'a) \times
\mathcal{K}(cU_1,U_1'b) \to \mathcal{K}(F'F_1'c,aFF_1)         \\
\bar{\Theta}(U_1'-.-U) &\maps& \mathcal{K}(bU,U'a) \times
\mathcal{K}(cU_1,U_1'b) \to \mathcal{K}(F'F_1'c,aFF_1) ,
 \ean
and a natural isomorphism $X$ between the pasting composites:
 \ban
 U_1'\Phi(-).\Phi_1(-)U &\maps& \mathcal{K}(F'b,aF) \times
 \mathcal{K}(F_1'c,bF_1) \to \mathcal{K}(cU_1U,U_1'Ua) \\
\bar{\Phi}(-F_1.F'-) &\maps& \mathcal{K}(F'b,aF) \times
 \mathcal{K}(F_1'c,bF_1) \to \mathcal{K}(cU_1U,U_1'Ua) .
 \ean
\end{prop}

\Proof Let $\xi\in\mathcal{K}(bU,U'a)$ and
$\xi_1\in\mathcal{K}(cU_1,U'_1b)$, then $V(\xi \times \xi_1)$ is
given by the following pasting composite of invertible
3-morphisms:
\[
\def\objectstyle{\scriptstyle}
\def\labelstyle{\scriptstyle}
  \xy
    (-75,0)*+{F'F_1'c}="t1";
    (-60,10)*+{F'F_1'cU_1F_1}="t2";
    (-45,20)*+{F'F_1'U_1'bF_1}="t3";
    (0,20)*+{F'bF_1}="t4";
    (45,20)*+{F'bUFF_1}="t5";
    (60,10)*+{F'U'aFF_1}="t6";
    (75,0)*+{aFF_1}="t7";
    (-60,-10)*+{F'F_1'cU_1F_1}="b2";
    (-45,-20)*+{F'F_1'cU_1UFF_1}="b3";
    (0,-20)*+{F'F_1'U_1'bUFF_1}="b4";
    (45,-20)*+{F'F_1U_1UaFF_1}="b5";
    (60,-10)*+{F'U'aFF_1}="b6";
        {\ar@{=>}^<<<<{FF_1'ci_1} "t1";"t2"};
        {\ar@{=>}^<<<{F'F_1'\xi_1F_1} "t2";"t3"};
        {\ar@{=>}^{F'e_1'bF_1} "t3";"t4"};
        {\ar@{=>}^{F'biF_1} "t4";"t5"};
        {\ar@{=>}^>>>{F'\xi FF_1} "t5";"t6"};
        {\ar@{=>}^>>>>{e'aFF_1} "t6";"t7"};
        {\ar@{=>}_<<<<{F'F_1'ci_1} "t1";"b2"};
        {\ar@{=>}_<<<{F'F_1'cU_1iF_1} "b2";"b3"};
        {\ar@{=>}_{F'F_1'\xi_1 UFF_1} "b3";"b4"};
        {\ar@{=>}_{F'F_1'U_1\xi FF_1} "b4";"b5"};
        {\ar@{=>}_>>>{\;\;F'e_1'U'aFF_1} "b5";"b6"};
        {\ar@{=>}_>>>>{e'aFF_1} "b6";"t7"};
        {\ar@{=} "t2";"b2"};
        {\ar@{=} "t6";"b6"};
        {\ar@{=>}_{FF_1U_1biF_1} "t3";"b4"};
        {\ar@{=>}_{F'e_1'bUFF_1} "b4";"t5"};
            {\ar@3{->}_{F'F_1'\xi_{1_i}F_1} (-41,5);(-41,0)};
            {\ar@3{->}^{F'e'_{1_{\xi}}FF_1} (41,5);(41,0)};
            {\ar@3{->}_{F'e'_{1_{bi}}F_1} (0,5);(0,0)};
  \endxy
\]
Since this 3-isomorphism is composed entirely of
$\cat{Gray}$-naturality isomorphisms, it is clear that $V$ is a
natural isomorphisms. Given $\zeta \in \mathcal{K}(F'b,aF)$ and
$\zeta_1 \in \mathcal{K}(F_1'c,bF_1)$ the natural isomorphism $X$
can be defined similarly:
\[
\def\objectstyle{\scriptstyle}
\def\labelstyle{\scriptstyle}
  \xy
    (-75,0)*+{cU_1U}="t1";
    (-60,10)*+{U_1'F_1'cU_1U}="t2";
    (-45,20)*+{U_1'bF_1U_1U}="t3";
    (0,20)*+{U_1'bU}="t4";
    (45,20)*+{U_1'U'F'bU}="t5";
    (60,10)*+{U_1'U'aFU}="t6";
    (75,0)*+{U_1'U'a}="t7";
    (-60,-10)*+{U_1'F_1'cU_1U}="b2";
    (-45,-20)*+{U_1'U'F'F_1'cU_1U}="b3";
    (0,-20)*+{U_1'U'F'bF_1U_1U}="b4";
    (45,-20)*+{U_1'U'aFF_1U_1U}="b5";
    (60,-10)*+{U_1'U'aFU}="b6";
        {\ar@{=>}^<<<<{i_1'cU_1U} "t1";"t2"};
        {\ar@{=>}^<<<{U_1'\zeta_1U_1U} "t2";"t3"};
        {\ar@{=>}^{U_1'be_1U} "t3";"t4"};
        {\ar@{=>}^{U_1'i'bU} "t4";"t5"};
        {\ar@{=>}^>>>{U_1'U'\zeta'U} "t5";"t6"};
        {\ar@{=>}^>>>>{U_1'U'ae} "t6";"t7"};
        {\ar@{=>}_<<<<{i_1'cU_1U} "t1";"b2"};
        {\ar@{=>}_<<<{U_1'i'F_1'cU_1U} "b2";"b3"};
        {\ar@{=>}_{U_1'UF'\zeta_1U_1U} "b3";"b4"};
        {\ar@{=>}_{U_1'U'\zeta F_1U_1U} "b4";"b5"};
        {\ar@{=>}_>>>{\;\;U_1'U'aFe_1U} "b5";"b6"};
        {\ar@{=>}_>>>>{U_1'U'ae} "b6";"t7"};
        {\ar@{=} "t2";"b2"};
        {\ar@{=} "t6";"b6"};
        {\ar@{=>}_>>>>>>>>>>>>>>{U_1'i'bF_1U_1U} "t3";"b4"};
        {\ar@{=>}_<<<<<<<<<<<<<<{U_1'U'F'be_1U} "b4";"t5"};
            {\ar@3{->}_{U_1'(i'_{\zeta_1})^{-1}U_1U} (-35,2);(-35,-3)};
            {\ar@3{->}^{U_1'U'(\zeta_{e_1})^{-1}U} (35,2);(35,-3)};
            {\ar@3{->}_{U_1'(i'b_{e_1U})^{-1}} (0,5);(0,0)};
  \endxy
\]
 \qed

As with the isomorphisms $W$ and $Y$ in
Proposition~\ref{comppsuedomatesONEcell}, we will find it
necessary to generically refer to the natural isomorphisms $V$ and
$X$ even though we may consider many different choices of
pseudoadjunctions and composable morphisms $a,b,c$.  Again, this
is allowed because these natural isomorphisms exist for every
choice of this data.

Before moving on to the theory of pseudomonads we first collect a
result about the functors $\Theta$ and $\Phi$.

\begin{prop}
Let $\Theta$ and $\Phi$ be as in Proposition~\ref{pseudomates}
with $I,E,i,e\maps F \dashv_p U = I',E',i',e'\maps F' \dashv_p U'$
and $a=1_A$ and $b=1_B$. Then in the category $\mathcal{K}(U,U)$
the object $\Phi(1_F)$ is isomorphic to the object $1_U$, and in
the category $\mathcal{K}(F,F)$ the object $\Theta(1_U)$ is
isomorphic to the object $1_F$.
\end{prop}

\Proof The isomorphisms are $I^{-1}$ and $E$ respectively. \qed

\subsection{Pseudomonads} \label{secpseudomonads}

Here we present the theory of pseudomonads.  For more details see
\cite{Lack,mar,mar2}.

\begin{defn}
A {\em pseudomonad} $\T =(T,\mu, \eta, \lambda,\rho,\alpha)$ on an
object $B$ of the $\cat{Gray}$-category $\mathcal{K}$ consists of
an endomorphism $T \maps B \to B$ together with:
\begin{itemize}
    \item {\em multiplication} for the pseudomonad: $\mu \maps T^2 \To T$,
    \item {\em unit} for the pseudomonad: $\eta \maps 1 \To T$, and
    \item coherence 3-isomorphisms
\[
 \xy
   (-20,10)*+{T}="L";
   (0,10)*+{TT}="M";
   (20,10)*+{T}="R";
   (0,-10)*+{T}="B";
    {\ar@{=>}^{\eta T} "L";"M"};
    {\ar@{=>}_{T\eta} "R";"M"};
    {\ar@{=>} "L";"B"};
    {\ar@{=>} "R";"B"};
    {\ar@{=>}^{\mu} "M";"B"};
    {\ar@3{->}_{\lambda} "M"+(-4,-4);(-8,2)  };
    {\ar@3{->}_{\rho} (8,2);"M"+(4,-4)  };
 \endxy
\qquad {\rm and} \qquad
 \xy
   (0,11)*+{T^3}="T";
   (-12,0)*+{T^2}="L";
   (12,0)*+{T^2}="R";
   (0,-11)*+{T}="B";
   {\ar@{=>}_{T\mu} "T";"L"};
   {\ar@{=>}^{\mu T} "T";"R"};
   {\ar@{=>}_{\mu} "L";"B"};
   {\ar@{=>}^{\mu} "R";"B"};
   {\ar@3{->}^{\alpha} "L"+(8,0); "R"+(-8,0) };
 \endxy
\]
\end{itemize}
such that the following two equations are satisfied:
\[
 \vcenter{\xy
   (-24,22)*+{T^4}="tl";
   (-4,22)*+{T^3}="tr";
   (-10,10)*+{T^3}="ml";
   (10,10)*+{T^2}="mr";
   (-24,2)*+{T^3}="o";
   (-10,-10)*+{T^2}="bl";
   (10,-10)*+{T}="br";
        {\ar@{=>}^{T^2\mu} "tl";"tr"};
        {\ar@{=>}_{\mu T^2} "tl";"o"};
        {\ar@{=>}^-{T\mu T} "tl";"ml"};
        {\ar@{=>}_{T\mu} "ml";"mr"};
        {\ar@{=>}_{\mu T} "ml";"bl"};
        {\ar@{=>}^{T\mu} "tr";"mr"};
        {\ar@{=>}_{\mu T} "o";"bl"};
        {\ar@{=>}^{\mu} "mr";"br"};
        {\ar@{=>}_{\mu} "bl";"br"};
          {\ar@3{->}^{T\alpha} "tr"+(-2,-4); "ml"+(2,4) };
          {\ar@3{->}_{\alpha T} "ml"+(-4.5,-2.5); "o"+(4.5,2.5) };
          {\ar@3{->}_{\alpha} "mr"+(-8,-8); "bl"+(8,8) };
 \endxy}
\qquad = \qquad \vcenter{ \xy
   (10,-10)*+{T}="tl";
   (-10,-10)*+{T^2}="tr";
   (-4,2)*+{T^2}="ml";
   (-24,2)*+{T^3}="mr";
   (10,10)*+{T^2}="o";
   (-4,22)*+{T^3}="bl";
   (-24,22)*+{T^4}="br";
        {\ar@{=>}_{\mu} "tl";"tr";};
        {\ar@{=>}^{\mu} "tl";"o";};
        {\ar@{=>}_-{\mu} "tl";"ml";};
        {\ar@{=>}^{T\mu} "ml";"mr";};
        {\ar@{=>}^{\mu T} "ml";"bl";};
        {\ar@{=>}_{\mu T} "tr";"mr";};
        {\ar@{=>}^{T\mu} "o";"bl";};
        {\ar@{=>}_{\mu T^2} "mr";"br";};
        {\ar@{=>}^{T^2\mu} "bl";"br";};
          {\ar@3{->}_{\alpha} "ml"+(-2,-4);  "tr"+(2,4)};
          {\ar@3{->}^{\alpha} "o"+(-4.5,-2.5);  "ml"+(4.5,2.5)};
          {\ar@3{->}_<<{\mu_{\mu}^{-1}} "bl"+(-8,-8);  "mr"+(8,8)};
 \endxy}
\]

\[
 \xy
   (-30,0)*+{T^2}="o";
   (0,12)*+{T^2}="t";
   (-12,0)*+{T^3}="l";
   (12,0)*+{T}="r";
   (0,-12)*+{T^2}="b";
        {\ar@{=>}^{T\eta T} "o";"l" };
        {\ar@{=>}^{T\mu} "l";"t" };
        {\ar@{=>}^{\mu} "t";"r" };
        {\ar@{=>}_{\mu} "b";"r" };
        {\ar@{=>}_{\mu T} "l";"b" };
        {\ar@3{->}_{\alpha} "t"+(0,-7);  "b"+(0,7)};
 \endxy
\qquad = \qquad
 \xy
   (30,0)*+{T}="o";
   (0,12)*+{T^3}="t";
   (-12,0)*+{T^2}="l";
   (12,0)*+{T^2}="r";
   (0,-12)*+{T^3}="b";
        {\ar@{=>}^{\mu} "r";"o" };
        {\ar@{=>}^{T\eta T } "l";"t" };
        {\ar@{=>} "l";"r" };
        {\ar@{=>}^{T\mu} "t";"r" };
        {\ar@{=>}_{\mu T} "b";"r" };
        {\ar@{=>}_{T\eta T} "l";"b" };
         {\ar@3{->}_{ T\lambda} "t"+(0,-4);  "t"+(0,-10)};
          {\ar@3{->}_{\rho T} "b"+(0,10);  "b"+(0,4)};
 \endxy
\]
\end{defn}

This definition was given by F. Marmolejo in~\cite{mar} and can be
understood as a pseudomonoid (in the sense of~\cite{DS}) in
$\mathcal{K}(B,B)$. An elegant treatment of pseudomonads is
presented in~\cite{Lack} where the `free living' or `walking'
pseudomonad is defined.  A {\em pseudocomonad} $\G
=(G,\delta,\varepsilon,\bar{\lambda},\bar{\rho},\bar{\alpha})$ is
defined by reversing the directions of the 2-cells in the
definition of a pseudomonad. A pseudocomonad can also be
understood as a pseudocomonoid in $\mathcal{K}(B,B)$.

\begin{prop}[Lack {\cite{Lack}}]
A pseudoadjunction $F \dashv_p U\maps B \to C$ in the
\cat{Gray}-category $\mathcal{K}$ induces a pseudomonad
$(UF,i,UeF,IF,UE,Ue_e^{-1})$ on the object $B$ in $\mathcal{K}$.
\end{prop}

\begin{prop} \label{PROPaa}
Let $\T =(T,\mu, \eta, \lambda,\rho,\alpha)$ be a pseudomonad on
an object $B$ in a $\cat{Gray}$-category $\mathcal{K}$ such that
the endomorphism $T \maps B \to B$ has a specified right
pseudoadjoint $G$ with counit $\sigma \maps TG \to 1$, unit $\iota
\maps 1 \to GT$, and coherences $\Upsilon \maps \sigma T.T\iota
\to 1$ and $\Sigma \maps 1 \to G\sigma.\iota G$.  Then mateship
under pseudoadjunction, together with the natural isomorphisms in
Propositions \ref{comppsuedomatesONEcell} and
\ref{comppseudomatesPseudoadj}, define a pseudocomonad
$\G=(G,\varepsilon,\delta,\bar{\lambda},\bar{\rho},\bar{\alpha})$
on $G$ with explicit formulas:
 \ban
  \varepsilon &:=& \Phi(\eta) = \sigma.\eta G
\\
  \delta &:=& \Phi(\mu) =
 G^2\sigma.G^2\mu G. G\iota TG. \iota G
\\
  \bar{\lambda} &:=&
  \def\objectstyle{\scriptstyle}
  \def\labelstyle{\scriptstyle}
  \xy
    (-65,0)*+{G\Phi(\eta).\Phi(\mu)}="1";
    (-20,0)*+{G\Phi(\eta).\Phi(1_T)G.\Phi(\mu)}="2";
    (12,0)*+{\Phi(\eta T).\Phi(\mu)}="3";
    (31,0)*+{\Phi(\mu.\eta T)}="4";
    (50,0)*+{\Phi(1_T)}="5";
    (61,0)*+{G}="6";
        {\ar@3{->}^-{\Sigma^{-1}} "5";"6"};
        {\ar@3{->}^-{\Phi(\lambda)} "4";"5"};
        {\ar@3{->}^-{Y} "3";"4"};
        {\ar@3{->}^-{X.\Phi(\mu)} "2";"3"};
        {\ar@3{->}^-{G\Phi(\eta).\Sigma G.\Phi(\mu)} "1";"2"};
  \endxy
\\
 \bar{\rho} &:=&
  \def\objectstyle{\scriptstyle}
  \def\labelstyle{\scriptstyle}
  \xy
    (68,0)*+{\Phi(\eta)G.\Phi(\mu)}="6";
    (22,0)*+{\Phi(1_T).\Phi(\eta)G.\Phi(\mu)}="5";
    (-6,0)*+{\Phi(T\eta).\Phi(\mu)}="4";
    (-28,0)*+{\Phi(\mu.T\eta)}="3";
    (-46,0)*+{\Phi(1_T)}="2";
    (-58,0)*+{G}="1";
        {\ar@3{->}^-{\Sigma} "1";"2"};
        {\ar@3{->}^-{\Phi(\rho)} "2";"3"};
        {\ar@3{->}^-{Y^{-1}} "3";"4"};
        {\ar@3{->}^-{X^{-1}} "4";"5"};
        {\ar@3{->}^-{\Sigma^{-1}.\Phi(\eta)G.\Phi(\mu)} "5";"6"};
  \endxy
    \\
 \bar{\alpha} &:=&
  \def\objectstyle{\scriptstyle}
  \def\labelstyle{\scriptstyle}
  \xy
    (70,-10)*+{G\Phi(\mu).\Phi(\mu)}="8";
    (21,-10)*+{G\Phi(\mu).\Phi(1_T)G.\Phi(\mu)}="7";
    (-18,-10)*+{\Phi(\mu T).\Phi(\mu)}="6";
    (-40,-10)*+{\Phi(\mu.\mu T)}="5";
    (55,0)*+{\Phi(\mu.T\mu)}="4";
    (35,0)*+{\Phi(T\mu).\Phi(\mu)}="3";
    (0,0)*+{GG\Phi(1_T).\Phi(\mu)G.\Phi(\mu)}="2";
    (-50,0)*+{\Phi(\mu)G.\Phi(\mu)}="1";
        {\ar@3{->}^-{GG\Sigma.\Phi(\mu)G.\Phi(\mu)} "1";"2"};
        {\ar@3{->}^-{X.\Phi(\mu)} "2";"3"};
        {\ar@3{->}^-{Y} "3";"4"};
        {\ar@3{->}^-{\Phi(\alpha)} (-55,-10);"5"};
        {\ar@3{->}^-{Y^{-1}} "5";"6"};
        {\ar@3{->}^-{X^{-1}.\Phi(\mu)} "6";"7"};
        {\ar@3{->}^-{G\Phi(\mu).\Sigma^{-1}G.\Phi(\mu)} "7";"8"};
  \endxy
 \ean
Under these circumstances $\G$ is said to be a pseudocomonad right
pseudoadjoint to the pseudomonad $\T$, denoted $\T \dashv_p \G$.
\end{prop}

\Proof Mateship under pseudoadjunction preserves composites along
morphisms and pseudoadjoints up to natural isomorphism by
Propositions \ref{comppsuedomatesONEcell} and
\ref{comppseudomatesPseudoadj}. Therefore because $\T =(T,\mu,
\eta, \lambda,\rho,\alpha)$ is a pseudomonad,
$\G=(G,\varepsilon,\delta,\bar{\lambda},\bar{\rho},\bar{\alpha})$
defines a pseudocomonad on $B$. \qed

\begin{defn}
A pseudomonad $\T$ in the $\cat{Gray}$-category $\mathcal{K}$ is
called a {\em Frobenius pseudomonad} if it is equipped with a map
$\varepsilon\maps T \to 1$ such that $\varepsilon.\mu$ is the
counit for a specified pseudoadjunction $T \dashv_p T$.
\end{defn}

We use this notion of Frobenius pseudomonad to define a Frobenius
pseudomonoid in a $\cat{Gray}$-monoid or semistrict monoidal
2-category.  A $\cat{Gray}$-monoid is just a one object
$\cat{Gray}$-category. In particular, if $\mathcal{K}$ is a
$\cat{Gray}$-category and $B$ is an object of $\mathcal{K}$, then
$\mathcal{K}(B,B)$ is a $\cat{Gray}$-monoid.  A Frobenius
pseudomonad on $B$ is then just a Frobenius pseudomonoid in the
$\cat{Gray}$-monoid $\mathcal{K}(B,B)$. This definition of
Frobenius pseudomonoid takes the minimalists approach, a
pseudomonoid equipped with the specified pseudoadjoint structure
that enables one to construct a pseudocomonoid structure.  For a
more explicit description of this definition see Street's work
\cite{StreetFrob}. One may prefer the definition of a Frobenius
pseudomonoid to be symmetrical: a pseudomonoid, and a
pseudocomonoid subject to compatibility conditions.  In the sequel
to this paper we explain the relationship between these two
perspectives which turn out to be equivalent in a precise sense
\cite{Lau2}.

We now describe the generalization of algebras for a monad and
construct the 2-category of pseudoalgebras based at $A$ for a
pseudomonad $\T$.  Pseudoalgebras for a 2-monad were first
explicitly defined by Street~\cite{Street4} and were well known to
the Australian category theory community at that time
\cite{CatSem}. For a treatment using the powerful machinery of
Blackwell-Kelly-Power~\cite{BKP}, see~\cite{BKPS}.  The treatment
we give here follows Marmolejo~\cite{mar}.

\begin{defn}
Let $\T$ be a pseudomonad in the $\cat{Gray}$-category
$\mathcal{K}$ and let $A$ be an object of $\mathcal{K}$. We define
a {\em pseudoalgebra based at $A$} for the pseudomonad $\T$ to
consist of:
\begin{itemize}
    \item a morphism $s \maps A \to B$,
    \item 2-morphisms $\nu \maps T s \To s$, and
    \item 3-isomorphisms
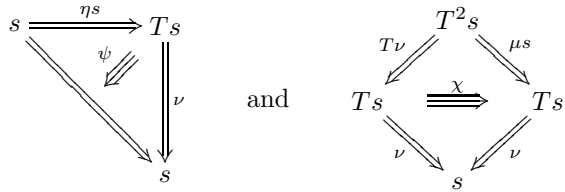
\[
 \xy
   (-20,10)*+{s}="L";
   (0,10)*+{T s}="M";
   (0,-10)*+{s}="B";
    {\ar@{=>}^{\eta s} "L";"M"};
    {\ar@{=>} "L";"B"};
    {\ar@{=>}^{\nu} "M";"B"};
    {\ar@3{->}_{\psi} "M"+(-4,-4);(-8,2)  };
 \endxy
\qquad {\rm and} \qquad
 \xy
   (0,11)*+{T^2s}="T";
   (-12,0)*+{T s}="L";
   (12,0)*+{T s}="R";
   (0,-11)*+{s}="B";
   {\ar@{=>}_{T \nu} "T";"L"};
   {\ar@{=>}^{\mu s} "T";"R"};
   {\ar@{=>}_{\nu} "L";"B"};
   {\ar@{=>}^{\nu} "R";"B"};
   {\ar@3{->}^{\chi} "L"+(8,0); "R"+(-8,0) };
 \endxy
\]
\end{itemize}
such that the following two equations are satisfied:
\[
 \vcenter{\xy
   (-24,22)*+{T^3s}="tl";
   (-4,22)*+{T^2s}="tr";
   (-10,10)*+{T^2s }="ml";
   (10,10)*+{T s}="mr";
   (-24,2)*+{T^2 s}="o";
   (-10,-10)*+{T s}="bl";
   (10,-10)*+{s}="br";
        {\ar@{=>}^{T^2 \nu} "tl";"tr"};
        {\ar@{=>}_{\mu T^2} "tl";"o"};
        {\ar@{=>}^-{T\mu T} "tl";"ml"};
        {\ar@{=>}_{T\mu} "ml";"mr"};
        {\ar@{=>}_{\mu \nu} "ml";"bl"};
        {\ar@{=>}^{T \nu} "tr";"mr"};
        {\ar@{=>}_{\mu \nu} "o";"bl"};
        {\ar@{=>}^{\nu} "mr";"br"};
        {\ar@{=>}_{\nu} "bl";"br"};
          {\ar@3{->}^{T \chi} "tr"+(-2,-4); "ml"+(2,4) };
          {\ar@3{->}_{\alpha s} "ml"+(-4.5,-2.5); "o"+(4.5,2.5) };
          {\ar@3{->}_{\chi} "mr"+(-8,-8); "bl"+(8,8) };
 \endxy}
\qquad = \qquad \vcenter{ \xy
   (10,-10)*+{s}="tl";
   (-10,-10)*+{T s}="tr";
   (-4,2)*+{T s}="ml";
   (-24,2)*+{T^2 s}="mr";
   (10,10)*+{T s}="o";
   (-4,22)*+{T^2 s}="bl";
   (-24,22)*+{T^3 s}="br";
        {\ar@{=>}_{\nu} "tl";"tr";};
        {\ar@{=>}^{\nu} "tl";"o";};
        {\ar@{=>}_-{\nu} "tl";"ml";};
        {\ar@{=>}^{T\nu} "ml";"mr";};
        {\ar@{=>}^{\mu s} "ml";"bl";};
        {\ar@{=>}_{\mu s} "tr";"mr";};
        {\ar@{=>}^{T\nu} "o";"bl";};
        {\ar@{=>}_{\mu T s} "mr";"br";};
        {\ar@{=>}^{T^2\nu} "bl";"br";};
          {\ar@3{->}_{\chi} "ml"+(-2,-4);  "tr"+(2,4)};
          {\ar@3{->}^{\chi} "o"+(-4.5,-2.5);  "ml"+(4.5,2.5)};
          {\ar@3{->}_<<{\mu_{\nu}^{-1}} "bl"+(-8,-8);  "mr"+(8,8)};
 \endxy}
\]

\[
 \xy
   (-30,0)*+{T s}="o";
   (0,12)*+{T s}="t";
   (-12,0)*+{T^2 s}="l";
   (12,0)*+{s}="r";
   (0,-12)*+{T s}="b";
        {\ar@{=>}^{T\eta s} "o";"l" };
        {\ar@{=>}^{T s} "l";"t" };
        {\ar@{=>}^{\nu} "t";"r" };
        {\ar@{=>}_{\nu} "b";"r" };
        {\ar@{=>}_{\mu s} "l";"b" };
        {\ar@3{->}_{\chi} "t"+(0,-7);  "b"+(0,7)};
 \endxy
\qquad = \qquad
 \xy
   (30,0)*+{s .}="o";
   (0,12)*+{T^2 s}="t";
   (-12,0)*+{T s}="l";
   (12,0)*+{T s}="r";
   (0,-12)*+{T^2 s}="b";
        {\ar@{=>}^{\nu} "r";"o" };
        {\ar@{=>}^{T\eta s} "l";"t" };
        {\ar@{=>} "l";"r" };
        {\ar@{=>}^{T\nu} "t";"r" };
        {\ar@{=>}_{\mu s} "b";"r" };
        {\ar@{=>}_{T\eta s} "l";"b" };
         {\ar@3{->}_{T\psi} "t"+(0,-4);  "t"+(0,-10)};
          {\ar@3{->}_{\rho s} "b"+(0,10);  "b"+(0,4)};
 \endxy
\]
\end{defn}

It is clear that for any morphism $r \maps A \to B$ in
$\mathcal{K}$, $T r$ with action $\mu r \maps T^2r \To Tr$ and
coherence $\lambda r \maps \mu r.\eta T r \Rrightarrow T r$ and
$\alpha r \maps \mu r . T\mu t \Rrightarrow \mu r. \mu Tr$ is a
pseudoalgebra based at $A$. We call the pseudoalgebra $Tr$ a {\em
free pseudoalgebra}.

\begin{defn} \label{def22}
Let \cat{$\T$-Alg}$_{\bf A}$ be the 2-category whose objects are
pseudoalgebras based at $A$ for the pseudomonad $\T$. A morphism
$(h,\varrho)\maps (s,\nu,\psi,\chi) \to (s',\nu',\psi',\chi')$ in
\cat{$\T$-Alg}$_{\bf A}$ consists of a 2-morphism $h \maps s \To
s'$ in $\mathcal{K}$ (a morphism in $\mathcal{K}(A,B)$), together
with an invertible 3-morphism
\[
 \xy
   (0,11)*+{T s}="T";
   (-12,0)*+{T s'}="L";
   (12,0)*+{ s}="R";
   (0,-11)*+{s'}="B";
   {\ar@{=>}_{T h} "T";"L"};
   {\ar@{=>}^{\nu} "T";"R"};
   {\ar@{=>}_{\nu'} "L";"B"};
   {\ar@{=>}^{h} "R";"B"};
   {\ar@3{->}^{\varrho} "L"+(8,0); "R"+(-8,0) };
 \endxy
\]
satisfying the following two equations:
 \ba
 \xy 0;/r.20pc/:
   (-20,10)*+{s}="t1";
   (0,10)*+{T s}="t2";
   (20,10)*+{T s'}="t3";
   (0,-5)*+{s}="b1";
   (20,-5)*+{s'}="b2";
        {\ar@{=>}^{\eta s} "t1";"t2"};
        {\ar@{=>}^{T h} "t2";"t3"};
        {\ar@{=>}_{1} "t1";"b1"};
        {\ar@{=>}^{\nu} "t2";"b1"};
        {\ar@{=>}^{\nu'} "t3";"b2"};
        {\ar@{=>}_{h} "b1";"b2"};
            {\ar@3{->}_{\varrho} "t3"+(-7,-5); "b1"+(8,5) };
            {\ar@3{->}^{\psi} "t2"+(-4,-3); (-8,3) };
 \endxy
\qquad &=& \qquad
 \xy 0;/r.20pc/:
   (0,10)*+{ s}="t1";
   (20,10)*+{T s}="t2";
   (0,-5)*+{s'}="b1";
   (20,-5)*+{T s'}="b2";
   (20,-20)*+{s'}="bb";
        {\ar@{=>}^{\eta s} "t1";"t2"};
        {\ar@{=>}_{h} "t1";"b1"};
        {\ar@{=>}^{T h} "t2";"b2"};
        {\ar@{=>}^{\eta s'} "b1";"b2"};
        {\ar@{=>}_{1} "b1";"bb"};
        {\ar@{=>}^-{\nu'} "b2";"bb"};
         {\ar@3{->}_{\eta_h} "t2"+(-7,-5); "b1"+(8,5) };
            {\ar@3{->}^{\psi'} "b2"+(-4,-3); (12,-13) };
 \endxy
  \ea
  \ba
\vcenter{\xy
   (-24,22)*+{T^2s}="tl";
   (-4,22)*+{T^2s'}="tr";
   (-10,10)*+{T s }="ml";
   (10,10)*+{T s'}="mr";
   (-24,2)*+{T s}="o";
   (-10,-10)*+{ s}="bl";
   (10,-10)*+{s'}="br";
        {\ar@{=>}^{T^2 h} "tl";"tr"};
        {\ar@{=>}_{\mu s} "tl";"o"};
        {\ar@{=>}^-{T\nu} "tl";"ml"};
        {\ar@{=>}_{T h} "ml";"mr"};
        {\ar@{=>}_{\nu} "ml";"bl"};
        {\ar@{=>}^{T \nu'} "tr";"mr"};
        {\ar@{=>}_{\nu} "o";"bl"};
        {\ar@{=>}^{\nu'} "mr";"br"};
        {\ar@{=>}_{h} "bl";"br"};
          {\ar@3{->}^{T \varrho} "tr"+(-2,-4); "ml"+(2,4) };
          {\ar@3{->}_{\chi} "ml"+(-4.5,-2.5); "o"+(4.5,2.5) };
          {\ar@3{->}_{\varrho} "mr"+(-8,-8); "bl"+(8,8) };
 \endxy}
\qquad &=& \qquad \vcenter{ \xy
   (10,-10)*+{s'}="tl";
   (-10,-10)*+{s}="tr";
   (-4,2)*+{T s'}="ml";
   (-24,2)*+{T s}="mr";
   (10,10)*+{T s'}="o";
   (-4,22)*+{T^2 s'}="bl";
   (-24,22)*+{T^2 s}="br";
        {\ar@{=>}_{h} "tl";"tr";};
        {\ar@{=>}^{\nu'} "tl";"o";};
        {\ar@{=>}_-{\nu'} "tl";"ml";};
        {\ar@{=>}^{T h} "ml";"mr";};
        {\ar@{=>}^{\mu s'} "ml";"bl";};
        {\ar@{=>}_{\nu} "tr";"mr";};
        {\ar@{=>}^{T\nu'} "o";"bl";};
        {\ar@{=>}_{\mu s} "mr";"br";};
        {\ar@{=>}^{T^2h} "bl";"br";};
          {\ar@3{->}_{\varrho} "ml"+(-2,-4);  "tr"+(2,4)};
          {\ar@3{->}^{\chi'} "o"+(-4.5,-2.5);  "ml"+(4.5,2.5)};
          {\ar@3{->}_<<{\mu_{h}^{-1}} "bl"+(-8,-8);  "mr"+(8,8)};
 \endxy}
 \ea

A 2-morphism $\xi \maps (h,\varrho) \To (h',\varrho')\maps
(\psi,\chi) \to (\psi',\chi')$ in \cat{$\T$-Alg}$_{\bf A}$ is a
3-morphisms $\xi \maps h \Rrightarrow h'$ such that the following
condition is satisfied:
 \ba
 \xy
   (-15,10)*+{T s}="tl";
   (15,10)*+{T s'}="tr";
   (-15,-10)*+{s}="bl";
   (15,-10)*+{s'}="br";
        {\ar@/^1pc/@{=>}^{T h} "tl";"tr"};
        {\ar@/_1pc/@{=>}_{T h'} "tl";"tr"};
        {\ar@{=>}_{\nu} "tl";"bl"};
        {\ar@{=>}^{\nu'} "tr";"br"};
        {\ar@{=>}_{h'} "bl";"br"};
            {\ar@3{->}_{T \xi} (0,12);(0,7)};
            {\ar@3{->}_{\varrho'} (0,-1);(0,-6)};
 \endxy
\qquad = \qquad
 \xy
   (-15,10)*+{T s}="tl";
   (15,10)*+{T s'}="tr";
   (-15,-10)*+{s}="bl";
   (15,-10)*+{s'}="br";
        {\ar@/^1pc/@{=>}^{h} "bl";"br"};
        {\ar@/_1pc/@{=>}_{ h'} "bl";"br"};
        {\ar@{=>}_{\nu} "tl";"bl"};
        {\ar@{=>}^{\nu'} "tr";"br"};
        {\ar@{=>}^{T h} "tl";"tr"};
            {\ar@3{->}_{ \xi} (0,-7);(0,-12)};
            {\ar@3{->}_{\varrho} (0,6);(0,1)};
 \endxy
 \ea
\end{defn}

Marmolejo has shown that given a morphism $K\maps A' \to A$ in
$\mathcal{K}$, one can define a change of base 2-functor
$\hat{K}\maps$ \cat{$\T$-Alg}$_{\bf A}$ $\to$ \cat{$\T$-Alg}$_{\bf
A'}$.  If $\xi \maps (h,\varrho)\To
(h',\varrho')\maps(s,\nu,\psi,\chi)\to(s',\nu',\psi',\chi')$ is in
\cat{$\T$-Alg}$_{\bf A}$, then its image under $\hat{K}$ is $\xi
K\maps (hK,\varrho K)\To (h'K,\varrho'K)\maps(sK,\nu K,\psi K,\chi
K)\to(s'K,\nu'K,\psi'K,\chi'K)$. If $k \maps K \To K'$ in
$\mathcal{K}$ then we get a pseudo natural transformation
$\hat{k}\maps \hat{K} \To \hat{K'}$ such that
$\hat{k}_{(s,\nu,\psi,\chi)}=(sk,\nu_{k}^{-1})$ and
$\hat{k}_{(h,\varrho)}=h_k^{-1}$.  If $\kappa \maps k \TO k'\maps
K \To K'$, then $\kappa_{(s,\nu,\psi,\chi)}=s\kappa$ defines a
modification $\hat{\kappa}\maps \hat{k} \TO \hat{k}$.  In fact,
this shows that the construction of $\T$-pseudoalgebras defines a
$\cat{Gray}$-functor $\cat{$\T$-Alg}\maps \mathcal{K}^{\op} \to
\cat{Gray}$.

For every object $A$ in $\mathcal{K}$ there is a forgetful
2-functor $U_A^{\T} \maps \cat{$\T$-Alg}_{\bf A} \to
\mathcal{K}(A,B)$
 \ban
 U_A^{\T} \maps \cat{$\T$-Alg}_{\bf A} &\to& \mathcal{K}(A,B) \\
 (s,\nu, \psi,\chi) &\mapsto& s \\
 (h,\varrho) &\mapsto& h \\
 \xi \maps h \Rrightarrow h' &\mapsto& \xi .
 \ean
This assignment extends to a $\cat{Gray}$-natural transformation
$U^{\T}\maps \cat{$\T$-Alg} \to \mathcal{K}(-,B)$. In Proposition
\ref{FUApseudoadjoints} we will define a left pseudoadjoint
$F_A^{\T}$ to the 2-functor $U_A^{\T}$, see also~\cite{mar}. In
Theorem~\ref{FextendsGray} we will show that this left
pseudoadjoint $F^{\T}_A$ extends to \cat{Gray}-natural
transformation $F^{\T}\maps \mathcal{K}(-,B)\to\cat{$\T$-Alg}$
left pseudoadjoint to $U^{\T}$ in the \cat{Gray}-category
$[\mathcal{K}^{\op},\cat{Gray}]$ described below.

Recall that $\cat{Gray}$ is the symmetric monoidal closed category
whose closed structure is given by the internal $\hom$ in
$\cat{Gray}$.  Hence, for $\cat{Gray}$-functors $F,G \maps
\mathcal{K} \to \mathcal{L}$ the internal $\hom$ $\cat{Gray}(F,G)$
in \cat{Gray} is the 2-category consisting of 2-functors, pseudo
natural transformations, and modifications. It is a standard
result from enriched category theory that $\cat{Gray}$-categories,
$\cat{Gray}$-functors, and $\cat{Gray}$-natural transformations
form a 2-category written $\cat{Gray-Cat}$~\cite{handbookII,kel}.
Furthermore, since $\cat{Gray}$ is a complete symmetric monoidal
closed category, if $\mathcal{K}$ is small, then the category of
$\cat{Gray}$-functors and $\cat{Gray}$-natural transformations can
be provided with the structure of a $\cat{Gray}$-category, written
$[\mathcal{K},\mathcal{L}]$.

The objects of $[\mathcal{K},\mathcal{L}]$ are the
\cat{Gray}-functors $F,G \maps \mathcal{K} \to \mathcal{L}$, and
the morphisms are the \cat{Gray}-natural transformations between
them.  The 2-category $\cat{Gray-Nat}(F,G)$ of \cat{Gray}-natural
transformations is given by the following equalizer:
\[
 \xy
    (-50,0)*+{\cat{Gray-Nat}(F,G)}="1";
    (-5,0)*+{\prod_{A \in \mathcal{K}} \mathcal{L}(FA,GA)}="2";
    (60,0)*+{\prod_{A',A'' \in \mathcal{K}} [ \mathcal{K}(A',A''), \mathcal{L}(FA',GA'')]}="3";
        {\ar@{>->} "1";"2"};
        {\ar^-{u} "2"+(17,2);"3"+(-32,2)};
        {\ar_-{v} "2"+(17,-2);"3"+(-32,-2)};
 \endxy
\]
where $u$ and $v$ are the morphisms corresponding via adjunction
and symmetry to the morphisms\footnote{Here we are using the
notation for $\mathcal{V}$-functors and $\mathcal{V}$-natural
transformations given in~\cite{handbookII}.}:
\[
 \xy
    (0,16)*+{\big(\prod_A \mathcal{L}(FA,GA)\big) \ten \mathcal{K}(A',A'')}="1";
    (0,0)*+{\mathcal{L}(FA',GA') \ten \mathcal{L}(GA',GA'')}="2";
    (0,-16)*+{\mathcal{L}(FA',GA'')}="3";
        {\ar@{->}^{p_{A'} \ten G_{A'A''}} "1";"2"};
        {\ar^-{c_{FA',GA',GA''}} "2";"3"};
 \endxy
    \qquad \qquad
  \xy
    (0,16)*+{\mathcal{K}(A',A'') \ten \big(\prod_A \mathcal{L}(FA,GA)\big)}="1";
    (0,0)*+{\mathcal{L}(FA',FA'') \ten \mathcal{L}(FA'',GA'')}="2";
    (0,-16)*+{\mathcal{L}(FA',GA'')}="3";
        {\ar@{->}^{F_{A'A'' \ten p_{A''}}} "1";"2"};
        {\ar^-{c_{FA',FA'',GA''}} "2";"3"};
 \endxy
\]
We will refer to the morphisms and 2-morphisms of the 2-category
$\cat{Gray-Nat}(F,G)$ as \cat{Gray}-modifications and
\cat{Gray}-perturbations respectively.  This terminology should
not be interpreted to mean some sort of `\cat{Gray} enriched
modification' or `\cat{Gray} enriched perturbation' since there is
no such notion as a $\mathcal{V}$-modification or
$\mathcal{V}$-perturbation for arbitrary enriching category
$\mathcal{V}$.

Let $\alpha, \beta \maps F \To G \maps \mathcal{K} \to \cat{Gray}$
be \cat{Gray}-natural transformations with $\mathcal{K}$ a small
\cat{Gray}-category.  A \cat{Gray}-modification $\theta \maps
\alpha \To \beta$ assigns to each object $A$ of $\mathcal{K}$ a
pseudo natural transformation $\theta_A \maps \alpha_A \to
\beta_A$ such that if $k \maps K \To K'\maps A' \to A''$ in
$\mathcal{K}$, then the following equality holds:
\[
 \xy
 (-20,10)*+{ \alpha_{A''}FK }="tl";
 (20,10)*+{  \alpha_{A''}FK'}="tr";
 (-20,-10)*+{ \beta_{A''}FK }="bl";
 (20,-10)*+{  \beta_{A''}FK'}="br";
        {\ar@{=>}^{\alpha_{A''}Fk} "tl";"tr"};
        {\ar@{=>}_{ \theta_{A''}FK} "tl";"bl"};
        {\ar@{=>}^{ \theta_{A''}FK'} "tr";"br"};
        {\ar@{=>}_{ \beta_{A''}Fk} "bl";"br"};
            {\ar@3{->}_-{\sim} (5,2.5);(-5,-2.5)};
 \endxy
\quad = \quad
 \xy
 (-20,10)*+{  GK\alpha_{A'}}="tl";
 (20,10)*+{GK'\alpha_{A'}}="tr";
 (-20,-10)*+{ GK\beta_{A'}}="bl";
 (20,-10)*+{  GK'\beta_{A'}}="br";
        {\ar@{=>}^{  Gk\alpha_{A'}} "tl";"tr"};
        {\ar@{=>}_{ GK \theta_{A'}} "tl";"bl"};
        {\ar@{=>}^{ GK'\theta_{A'}} "tr";"br"};
        {\ar@{=>}_{ Gk\beta_{A'} } "bl";"br"};
            {\ar@3{->}_-{\sim} (5,2.5);(-5,-2.5)};
 \endxy
\]
If $\Omega \maps \theta,\varphi \maps \alpha \To \beta \maps F \to
G \maps \mathcal{K} \to \cat{Gray}$ are \cat{Gray}-modifications,
a \cat{Gray}-perturbation assigns to each object $A \in
\mathcal{K}$ a modification $\Omega_A \maps \theta_A \to
\varphi_A$ such that if $\kappa \maps k \TO k' \maps K \To K'\maps
A' \to A''$ in $\mathcal{K}$, then the following equality holds:
\[
 \xy
 (-20,16)*+{ \alpha_{A''}FK}="tl";
 (20,16)*+{    \alpha_{A''}FK'}="tr";
 (-20,-14)*+{\beta_{A''}FK }="bl";
 (20,-14)*+{  \beta_{A''}FK'}="br";
        {\ar@/^1.5pc/@{=>}^{ \alpha_{A''}Fk} "tl";"tr"};
        {\ar@/_1.5pc/@{=>}_{\alpha_{A''}Fk'} "tl";"tr"};
        {\ar@/_1.8pc/@{=>}_{\varphi_{A''}FK } "tl";"bl"};
        {\ar@/^1.8pc/@{=>}^{\theta_{A''}FK } "tl";"bl"};
        {\ar@{=>}^{ \theta_{A''}FK'} "tr";"br"};
        {\ar@{=>}_{ \beta_{A''}Fk'} "bl";"br"};
            {\ar@3{->}_-{\sim} (10,2.5);(2,-4.5)};
            {\ar@3{->}_-{  \alpha_{A''}F\kappa} (5,18);(5,12)};
            {\ar@3{->}_-{   \Omega_{A''}FK} (-15,0);(-25,0)};
 \endxy
\quad = \quad
 \xy
 (-20,16)*+{GK\alpha_{A'}}="tl";
 (20,16)*+{GK'\alpha_{A'}}="tr";
 (-20,-14)*+{GK\beta_{A'}}="bl";
 (20,-14)*+{GK'\beta_{A'}}="br";
        {\ar@{=>}^{Gk\alpha_{A'}} "tl";"tr"};
        {\ar@{=>}_{GK\varphi_{A'} } "tl";"bl"};
        {\ar@/^1.8pc/@{=>}^{GK'\theta_{A'}} "tr";"br"};
        {\ar@/_1.85pc/@{=>}_{ GK'\varphi_{A'}} "tr";"br"};
        {\ar@/^1.5pc/@{=>}^{Gk\beta_{A'} } "bl";"br"};
        {\ar@/_1.5pc/@{=>}_{Gk'\beta_{A'} } "bl";"br"};
            {\ar@3{->}_-{\sim} (-2,8.5);(-10,1.5)};
            {\ar@3{->}_-{ G\kappa\beta{A'}} (5,-12);(5,-18)};
            {\ar@3{->}_<<<<<<<{GK'\Omega_{A'} } (26,0);(13,0)};
 \endxy
\]

\begin{prop} \label{FUApseudoadjoints}
Let $\T$ be a pseudomonad in the $\cat{Gray}$-category
$\mathcal{K}$. Then the forgetful 2-functor $U_A^{\T} \maps
\cat{$\T$-Alg}_{\bf A} \to \mathcal{K}(A,B)$ has a left
pseudoadjoint $F_A^{\T} \maps \mathcal{K}(A,B) \to
\cat{$\T$-Alg}_{\bf A}$ in the $\cat{Gray}$-category $\cat{Gray}$
given by sending each object $r$ of $\mathcal{K}(A,B)$ to the
corresponding free pseudoalgebra $(Tr,\mu r, \lambda r, \alpha
r)$, each morphism $h \maps r \to r'$ to $(Th, \mu_h^{-1})$, and
each 2-morphism $\xi \maps h \Rrightarrow h'$ to $T\xi \maps Th
\Rrightarrow Th'$.
\end{prop}

\Proof  First we show that $F_A^{\T}$ is a 2-functor. It is clear
that if $r \maps A \to B$ in $\mathcal{K}$, then $(Tr,\mu r,
\lambda r, \alpha r)$ is a pseudoalgebra based at $A$.  If $h
\maps r \To r'\maps A \to B$ then $Th$ is a morphism of
pseudoalgebras with $\mu_h^{-1}$ playing the role of the
invertible 2-morphism $\varrho$.  Indeed, the axioms:
 \ban
 \xy 0;/r.24pc/:
   (-20,10)*+{T r}="t1";
   (0,10)*+{T T r}="t2";
   (20,10)*+{T T r'}="t3";
   (0,-5)*+{T r}="b1";
   (20,-5)*+{T r'}="b2";
        {\ar@{=>}^{\eta T r} "t1";"t2"};
        {\ar@{=>}^{T T h} "t2";"t3"};
        {\ar@{=>}_{1} "t1";"b1"};
        {\ar@{=>}^{\mu r} "t2";"b1"};
        {\ar@{=>}^{\mu r'} "t3";"b2"};
        {\ar@{=>}_{T h} "b1";"b2"};
            {\ar@3{->}^{\mu_{h}^{-1}} "t3"+(-7,-5); "b1"+(8,5) };
            {\ar@3{->}^{\lambda r} "t2"+(-4,-3); (-8,3) };
 \endxy
\qquad &=& \qquad
 \xy 0;/r.24pc/:
   (0,10)*+{ T r}="t1";
   (20,10)*+{T T r}="t2";
   (0,-5)*+{T r'}="b1";
   (20,-5)*+{T T r'}="b2";
   (20,-20)*+{T r'}="bb";
        {\ar@{=>}^{\eta T r} "t1";"t2"};
        {\ar@{=>}_{T h} "t1";"b1"};
        {\ar@{=>}^{T T h} "t2";"b2"};
        {\ar@{=>}^{\eta T r'} "b1";"b2"};
        {\ar@{=>}_{1} "b1";"bb"};
        {\ar@{=>}^{\mu r'} "b2";"bb"};
         {\ar@3{->}_{\mu_{ h}^{-1}} "t2"+(-7,-5); "b1"+(8,5) };
            {\ar@3{->}^{\lambda r'} "b2"+(-4,-3); (12,-13) };
 \endxy \\ & & \\
\vcenter{\xy 0;/r.24pc/:
   (-24,22)*+{T^2T r}="tl";
   (-4,22)*+{T^2T r'}="tr";
   (-10,10)*+{T T r }="ml";
   (10,10)*+{T T r'}="mr";
   (-24,2)*+{T T r}="o";
   (-10,-10)*+{ T r}="bl";
   (10,-10)*+{T r'}="br";
        {\ar@{=>}^-{T^2 T h} "tl";"tr"};
        {\ar@{=>}_{\mu T r} "tl";"o"};
        {\ar@{=>}^-{T\mu r} "tl";"ml"};
        {\ar@{=>}_{T T h} "ml";"mr"};
        {\ar@{=>}_{\mu r} "ml";"bl"};
        {\ar@{=>}^{T \mu r'} "tr";"mr"};
        {\ar@{=>}_{\mu r} "o";"bl"};
        {\ar@{=>}^{\mu r'} "mr";"br"};
        {\ar@{=>}_{T h} "bl";"br"};
          {\ar@3{->}^{T \mu_{h}^{-1}} "tr"+(-2,-4); "ml"+(2,4) };
          {\ar@3{->}_{\alpha r} "ml"+(-4.5,-2.5); "o"+(4.5,2.5) };
          {\ar@3{->}_{\mu_{h}^{-1}} "mr"+(-8,-8); "bl"+(8,8) };
 \endxy}
\qquad &=& \qquad \vcenter{ \xy 0;/r.24pc/:
   (10,-10)*+{T r'}="tl";
   (-10,-10)*+{T r}="tr";
   (-4,2)*+{T T r'}="ml";
   (-24,2)*+{T T r}="mr";
   (10,10)*+{T T r'}="o";
   (-4,22)*+{T^2 T r'}="bl";
   (-24,22)*+{T^2 T r}="br";
        {\ar@{=>}_{T h} "tl";"tr";};
        {\ar@{=>}^{\mu r'} "tl";"o";};
        {\ar@{=>}_-{\mu r'} "tl";"ml";};
        {\ar@{=>}^{T T h} "ml";"mr";};
        {\ar@{=>}^{\mu T r'} "ml";"bl";};
        {\ar@{=>}_{\mu r} "tr";"mr";};
        {\ar@{=>}^{T\mu r'} "o";"bl";};
        {\ar@{=>}_{\mu T r} "mr";"br";};
        {\ar@{=>}^-{T^2T h} "bl";"br";};
          {\ar@3{->}_{\mu_{ h}^{-1}} "ml"+(-2,-4);  "tr"+(2,4)};
          {\ar@3{->}^{\alpha r'} "o"+(-4.5,-2.5);  "ml"+(4.5,2.5)};
          {\ar@3{->}_<<{\mu_{h}^{-1}} "bl"+(-8,-8);  "mr"+(8,8)};
 \endxy}
 \ean
 follow directly from the axioms of a $\cat{Gray}$-category. If $\xi \maps h \Rrightarrow
 h'$ is a 2-morphism in $\mathcal{K}(A,B)$ then $T\xi$ is a 2-morphism of
 pseudoalgebras since the equation
  \[
 \xy
   (-15,10)*+{T T r}="tl";
   (15,10)*+{T T r'}="tr";
   (-15,-10)*+{T r}="bl";
   (15,-10)*+{T r'}="br";
        {\ar@/^1pc/@{=>}^{T h} "tl";"tr"};
        {\ar@/_1pc/@{=>}_{T h'} "tl";"tr"};
        {\ar@{=>}_{\mu r} "tl";"bl"};
        {\ar@{=>}^{\mu r'} "tr";"br"};
        {\ar@{=>}_{h'} "bl";"br"};
            {\ar@3{->}_{T \xi} (0,12);(0,7)};
            {\ar@3{->}_{\mu^{-1}_{h'}} (0,-1);(0,-6)};
 \endxy
\qquad = \qquad
 \xy
   (-15,10)*+{T T r}="tl";
   (15,10)*+{T T r'}="tr";
   (-15,-10)*+{T r}="bl";
   (15,-10)*+{T r'}="br";
        {\ar@/^1pc/@{=>}^{h} "bl";"br"};
        {\ar@/_1pc/@{=>}_{ h'} "bl";"br"};
       {\ar@{=>}_{\mu r} "tl";"bl"};
        {\ar@{=>}^{\mu r'} "tr";"br"};
        {\ar@{=>}^{T h} "tl";"tr"};
            {\ar@3{->}_{ \xi} (0,-7);(0,-12)};
            {\ar@3{->}_{\mu_h^{-1}} (0,6);(0,1)};
 \endxy
\]
is satisfied, again, because of the axioms of a
$\cat{Gray}$-category. It is straightforward to verify that all
composites and identities are preserved.

Since $\mathcal{K}$ is a $\cat{Gray}$-category, $\mathcal{K}(A,B)$
is an object of $\cat{Gray}$, i.e., a 2-category.   It is clear
that if $\T$ is a pseudomonad on $B$ in $\mathcal{K}$, then
$\mathcal{K}(A,\T)$ is a pseudomonad on the 2-category
$\mathcal{K}(A,B)$ in the $\cat{Gray}$-category $\cat{Gray}$.  We
now show that $F_A^{\T}$ is the left pseudoadjoint of $U_A^{\T}$.
Define the pseudo natural transformation $i_A^{\T}\maps 1 \To
U_A^{\T}F_A^{\T}$ to be the unit $\mathcal{K}(A,\eta)$ for the
pseudomonad $\mathcal{K}(A,\T)$ in $\cat{Gray}$.  That is,
$i^{\T}_A (s) = \eta s$ and $i^{\T}_A(h)= h_{\eta}^{-1}$. We
define $e_A^{\T}\maps F_A^{\T}U_A^{\T} \To 1$ on any pseudoalgebra
$(s,\nu,\psi,\chi)$ to be the morphism of pseudoalgebras given by
setting $h=\nu$ and $\varrho = \chi$. The relevant axioms for a
morphism of pseudoalgebras follow using the axioms of the
pseudoalgebra $(s,\nu,\psi,\chi)$ and the pseudomonad $\T$. To
establish the pseudo naturality of $e_A^{\T}$ let
$(h,\varrho)\maps (s,\nu,\psi,\xi) \to (s',\nu',\psi',\xi')$ be a
map of pseudoalgebras.  Then the pseudo naturality 2-cell filling
the square:
\[
 \xymatrix{ (Ts,\mu s, \lambda s, \alpha s) \ar[d]_{(\nu,\chi)}\ar[rr]^{(Th, \mu_h^{-1})}
 && (Ts,\mu s', \lambda s', \alpha s') \ar[d]^{(\nu',\chi')}
 \\
 (s, \nu,\psi,\chi) \ar[rr]_{(h,\varrho)}&& (s', \nu',\psi',\chi')
 }
\]
is given by the 3-isomorphism $\varrho$.

We now describe the coherence modifications for this
pseudoadjunction in $\cat{Gray}$. Define $I_A^{\T} \maps
1_{U^{\T}_A} \Rrightarrow U_A^{\T}e_A^{\T}.i_A^{\T}U_A^{\T}$ on
the pseudoalgebra $(s,\nu,\psi,\chi)$ to be $\psi^{-1}$.  It is
easy to check that this map defines a modification between pseudo
natural transformations $1_{U^{\T}_A}$ and
$U_A^{\T}e_A^{\T}.i_A^{\T}U_A^{\T}$.  We define $E_A^{\T} \maps
e_A^{\T}F_A^{\T}.F_A^{\T}i_A^{\T} \Rrightarrow 1_{F^{\T}_A}$ on
the map $s$ to be the 2-homomorphism of pseudoalgebras
$\rho^{-1}s$. The fact that this map is a modification follows
from the pseudo naturality of $\mathcal{K}(A,\rho)$.

To establish the coherence of $I_A^{\T}$ and $E_A^{\T}$ consider
the diagram below:
\[
 \def\objectstyle{\scriptstyle}
\def\labelstyle{\scriptstyle}
 \xy 0;/r.30pc/:
   (-44,0)*+{\left(T s,\mu s, \lambda s, \alpha s\right)}="ll";
   (0,14)*+{\left(T s,\mu s, \lambda s, \alpha s\right)}="t";
   (-14,0)*+{\left(TT s,\mu T s, \lambda T s, \alpha T s\right)}="l";
   (14,0)*+{\left(s,\nu, \psi, \chi \right)}="r";
   (0,-14)*+{\left(T s,\mu s, \lambda s, \alpha s\right)}="b";
        {\ar@{=>}_>>>>>>>{(T \nu, \mu_{\nu}^{-1})} "l";"t"};
        {\ar@{=>}^{(\mu s, \alpha s)} "l";"b"};
        {\ar@{=>}^{(\nu,\chi)} "t";"r"};
        {\ar@{=>}_{(\nu,\chi)} "b";"r"};
        {\ar@{=>}^-{(T \eta s, \mu_{\eta s}^{-1})} "ll";"l"};
        {\ar@/^1pc/@{=>}^1 "ll";"t"};
        {\ar@/_1pc/@{=>}_1 "ll";"b"};
        {\ar@3{->}^{\chi} "t"+(2,-11);"b"+(2,11)};
        {\ar@3{->}_{T \psi^{-1}} (-16,9);(-16,4)};
        {\ar@3{->}_{\rho^{-1}T s} (-16,-4);(-16,-9)};
 \endxy
\]
The commutativity of this diagram can be deduced from the second
coherence condition in the definition of a pseudoalgebra.  The
other coherence law for a pseudoadjunction,
\[
 \xy
   (38,0)*+{s}="ll";
   (0,14)*+{T s}="t";
   (-14,0)*+{T^2 s}="l";
   (14,0)*+{T s}="r";
   (0,-14)*+{T s}="b";
        {\ar@{=>}^{\eta s} "l";"t"};
        {\ar@{=>}_{\eta s} "l";"b"};
        {\ar@{=>}_{\eta T s} "t";"r"};
        {\ar@{=>}^{T \eta s} "b";"r"};
        {\ar@{=>}^{\mu s} "r";"ll"};
        {\ar@/^1pc/@{=>}^1 "t";"ll"};
        {\ar@/_1pc/@{=>}_1 "b";"ll"};
        {\ar@3{->}_{\eta_{\eta}^{-1}} "t"+(-2,-11);"b"+(-2,11)};
        {\ar@3{->}^{\lambda^{-1}} (16,9);(16,4)};
        {\ar@3{->}^{\rho^{-1}} (16,-4);(16,-9)};
 \endxy
\]
can be deduced from the axioms of a pseudomonad, see Proposition
8.1 of~\cite{mar}. \qed

\begin{thm} \label{FextendsGray}
The collection of pseudoadjunctions:
$$I^{T}_A,E^{\T}_A,i^{T}_A,e^{T}_A\maps F^{T}_A \dashv_p U^{T}_A
\maps \cat{$\T$-Alg}_{A} \to \mathcal{K}(A,B)$$ defined for each
$A$ in Proposition~\ref{FUApseudoadjoints} extend to a
pseudoadjunction
$$I^{\T},E^{\T},i^{\T},e^{\T}\maps F^{\T} \dashv_p U^{\T}\maps
\cat{$\T$-Alg} \to \mathcal{K}(-,B)$$ in the \cat{Gray}-category
$[\mathcal{K}^{\op},\cat{Gray}]$. In particular, $F^{\T}$ is a
$\cat{Gray}$-natural transformation, $i^{\T}, e^{\T}$ are
\cat{Gray}-modifications, and $I^{\T}, E^{\T}$ are
\cat{Gray}-perturbations.
\end{thm}

\Proof.  Let $\kappa \maps k \TO k' \maps K \To K \maps A' \to A$
in the \cat{Gray}-category $\mathcal{K}$. To check that $F^{\T}$
is a $\cat{Gray}$-natural transformation we must check that:
\[
 \xymatrix{
\mathcal{K}(A,B) \ar[r]^{F^{\T}_A} \ar[d]_{\mathcal{K}(K,B)} &
\cat{$\T$-Alg} \ar[d]^{\hat{K}} \\ \mathcal{K}(A',B)
\ar[r]_{F^{\T}_{A'}} & \cat{$\T$-Alg} }
\]
commutes. Let $\xi \maps h \To h' \maps s \To s'$ in
$\mathcal{K}(A,B)$.  Then
 \ban
 \hat{K}F_A^{\T}(\xi)
 &=& \scriptstyle
T\xi K \maps (ThK,\mu_h^{-1}K)\To
 (Th'K,\mu_{h'}^{-1}K)\maps (TsK,\mu sK, \lambda sK, \alpha sK)\to
 (Ts'K,\mu s'K, \lambda s'K, \alpha s'K) \\
  &=&  \scriptstyle
T\xi K \maps (ThK,\mu_{hK}^{-1})\To
 (Th'K,\mu_{h'K}^{-1})\maps (TsK,\mu sK, \lambda sK, \alpha sK)\to
 (Ts'K,\mu s'K, \lambda s'K, \alpha s'K) \\
  &=&
  F_{A'}^{\T}\mathcal{K}(K,B)(\xi).
 \ean

Next we show that the collection of $e^{\T}_A$ define a
\cat{Gray}-modification $e^{\T} \maps F^{\T} U^{\T} \To 1_{{\rm
T-Alg}}$. To see that $e^{\T}$ is a \cat{Gray}-modification we
must verify the following equality:
\[
 \xy
 (-20,10)*+{ F^{\T}_{A'}U^{\T}_{A'}\hat{K} }="tl";
 (20,10)*+{  F^{\T}_{A'}U^{\T}_{A'}\hat{K}'}="tr";
 (-20,-10)*+{  \hat{K} }="bl";
 (20,-10)*+{   \hat{K}'}="br";
        {\ar@{=>}^{F^{\T}_{A'}U^{\T}_{A'}\hat{k}} "tl";"tr"};
        {\ar@{=>}_{ e^{\T}_{A'}\hat{K}} "tl";"bl"};
        {\ar@{=>}^{ e^{\T}_{A'}\hat{K}'} "tr";"br"};
        {\ar@{=>}_{  \hat{k}} "bl";"br"};
            {\ar@3{->}_-{\hat{k}_{e^{\T}_{A'}}^{-1}} (5,2.5);(-5,-2.5)};
 \endxy
\quad = \quad
 \xy
 (-20,10)*+{  \hat{K}F^{\T}_AU^{\T}_{A}}="tl";
 (20,10)*+{\hat{K}'F^{\T}_AU^{\T}_{A}}="tr";
 (-20,-10)*+{ \hat{K} }="bl";
 (20,-10)*+{  \hat{K}' }="br";
        {\ar@{=>}^{  \hat{k}F^{\T}_AU^{\T}_{A}} "tl";"tr"};
        {\ar@{=>}_{ \hat{K} e^{\T}_{A}} "tl";"bl"};
        {\ar@{=>}^{ \hat{K}'e^{\T}_{A}} "tr";"br"};
        {\ar@{=>}_{ \hat{k}  } "bl";"br"};
            {\ar@3{->}_-{\hat{k}_{e^{\T}_{A}}} (5,2.5);(-5,-2.5)};
 \endxy
\]
Note that
 \ban
 F^{\T}_{A'}U^{\T}_{A'}\hat{K}(s,\nu,\psi,\chi)
    &=&
(TsK,\mu sK, \lambda sK, \alpha sK) =
\hat{K}F^{\T}_{A}U^{\T}_{A}(s,\nu,\psi,\chi) \\
 F^{\T}_{A'}U^{\T}_{A'}\hat{K}'(s,\nu,\psi,\chi)
 &=& (TsK',\mu sK',
\lambda sK', \alpha sK') =
\hat{K}'F^{\T}_{A}U^{\T}_{A}(s,\nu,\psi,\chi).
 \ean
Furthermore,
$
 F^{\T}_{A'}U^{\T}_{A'}\hat{k} (s,\nu,\psi,\chi)
   = (Tsk,\mu_{sk}^{-1})
   = (Tsk,\mu_{s}^{-1}k)
      = \hat{k}F^{\T}_AU^{\T}_A
      (s,\nu,\psi,\chi)$,
and $e^{\T}_{A'}\hat{K}(s,\nu,\psi,\chi)= (\nu K,\chi K) =
\hat{K}e^{\T}_A(s,\nu,\psi,\chi)$. Hence, the desired equality is
satisfied by the $\cat{Gray}$ axioms asserting the following
equality:
\[
 \def\objectstyle{\scriptstyle}
 \def\labelstyle{\scriptscriptstyle}
 \xy
 (-20,16)*+{ F^{\T}_{A'}U^{\T}_{A'}\hat{K} = \hat{K}F^{\T}_{A}U^{\T}_{A}}="tl";
 (20,16)*+{  F^{\T}_{A'}U^{\T}_{A'}\hat{K}' = \hat{K}'F^{\T}_{A}U^{\T}_{A}}="tr";
 (-20,-14)*+{ \hat{K} }="bl";
 (20,-14)*+{ \hat{K}'}="br";
        {\ar@/^1.5pc/@{=>}^{ \hat{k}F^{\T}_{A}U^{\T}_{A}} "tl";"tr"};
        {\ar@/_1.5pc/@{=>}_{F^{\T}_{A'}U^{\T}_{A'}\hat{k}} "tl";"tr"};
        {\ar@/_1.8pc/@{=>}_{\hat{K}e^{\T}_{A}} "tl";"bl"};
        {\ar@/^1.8pc/@{=>}^{e^{\T}_{A'}\hat{K} } "tl";"bl"};
        {\ar@{=>}^{ e^{\T}_{A'}\hat{K}'} "tr";"br"};
        {\ar@{=>}_{ \hat{k}} "bl";"br"};
            {\ar@3{->}^-{(\hat{k}_{e^{\T}_{A'}})^{-1}} (10,2.5);(2,-4.5)};
            {\ar@3{->}_-{  1} (0,18);(0,12)};
            {\ar@3{->}_-{  1} (-15,0);(-25,0)};
 \endxy
\quad = \quad
 \xy
 (-22,16)*+{ F^{\T}_{A'}U^{\T}_{A'}\hat{K} = \hat{K}F^{\T}_{A}U^{\T}_{A}}="tl";
 (20,16)*+{  F^{\T}_{A'}U^{\T}_{A'}\hat{K}' = \hat{K}'F^{\T}_{A}U^{\T}_{A}}="tr";
 (-22,-14)*+{ \hat{K} }="bl";
 (20,-14)*+{ \hat{K}'}="br";
        {\ar@{=>}^{\hat{k}F^{\T}_{A}U^{\T}_{A}} "tl";"tr"};
        {\ar@{=>}_{\hat{K}e^{\T}_{A} } "tl";"bl"};
        {\ar@/^1.8pc/@{=>}^{e^{\T}_{A'}\hat{K}'} "tr";"br"};
        {\ar@/_1.85pc/@{=>}_{ \hat{K}'e^{\T}_{A}} "tr";"br"};
        {\ar@{=>}^{\hat{k}} "bl";"br"};
            {\ar@3{->}_-{\hat{k}_{e^{\T}_{A}}} (-2,4.5);(-10,-2.5)};
            {\ar@3{->}_-{1 } (24,0);(15,0)};
 \endxy
\]

To show that the collection of pseudo natural transformations
$i^{\T}_A$ define a \cat{Gray}-modification $i^{\T}\maps
1_{\mathcal{K}(-,B)} \To U^{\T}F^{\T}$ we must verify the
following equality:
\[
 \xy
 (-25,10)*+{  \mathcal{K}(K,B) }="tl";
 (25,10)*+{   \mathcal{K}(K',B)}="tr";
 (-25,-10)*+{ U^{\T}_{A'}F^{T}_{A'}\mathcal{K}(K,B) }="bl";
 (25,-10)*+{  U^{\T}_{A'}F^{T}_{A'}\mathcal{K}(K',B)}="br";
        {\ar@{=>}^-{ \mathcal{K}(k,B)} "tl";"tr"};
        {\ar@{=>}_-{ i^{\T}_{A'}\mathcal{K}(K,B)} "tl";"bl"};
        {\ar@{=>}^-{ i^{\T}_{A'}\mathcal{K}(K',B)} "tr";"br"};
        {\ar@{=>}_-{ U^{\T}_{A'}F^{T}_{A'}\mathcal{K}(k,B)} "bl";"br"};
            {\ar@3{->}_-{i^{\T}_{A'}   \;_{\mathcal{K}(k,B)}^{-1} \;\;} (5,2.5);(-5,-2.5)};
 \endxy
\ =
\]
\[
 \xy
 (-25,10)*+{  \mathcal{K}(K,B) }="tl";
 (25,10)*+{\mathcal{K}(K',B) }="tr";
 (-25,-10)*+{ \mathcal{K}(K,B)U^{\T}_{A}F^{T}_{A}}="bl";
 (25,-10)*+{  \mathcal{K}(K',B)U^{\T}_{A}F^{T}_{A}}="br";
        {\ar@{=>}^-{  \mathcal{K}(k,B) } "tl";"tr"};
        {\ar@{=>}_-{ \mathcal{K}(K,B) i^{\T}_{A}} "tl";"bl"};
        {\ar@{=>}^-{ \mathcal{K}(K',B)i^{\T}_{A}} "tr";"br"};
        {\ar@{=>}_-{ \mathcal{K}(k,B)U^{\T}_{A}F^{T}_{A} } "bl";"br"};
            {\ar@3{->}_-{\mathcal{K}(k,B)_{i^{\T}_A} \;\;\;\;} (5,2.5);(-5,-2.5)};
 \endxy
\]
Note that $$U^{\T}_{A'}F^{T}_{A'}\mathcal{K}(K,B)(s) = TsK =
\mathcal{K}(K,B)U^{\T}_{A}F^{T}_{A}(s)$$ and
$$U^{\T}_{A'}F^{T}_{A'}\mathcal{K}'(K,B)(s) = TsK' =
\mathcal{K}(K',B)U^{\T}_{A}F^{T}_{A}(s).$$  Furthermore,
$$i^{\T}_{A'}\mathcal{K}(K,B)(s) = (\eta sK) =
\mathcal{K}(K,B)i^{\T}_{A}(s)$$ and
$$U^{\T}_{A'}F^{\T}_{A'}\mathcal{K}(k,B)(s) = (Tsk) =
\mathcal{K}(k,B)U^{\T}_{A}F^{\T}_{A}(s).$$ Hence, the desired
equality follows from the $\cat{Gray}$ axioms asserting the
equality
\[
 \def\objectstyle{\scriptstyle}
 \def\labelstyle{\scriptscriptstyle}
 \xy
 (-36,16)*+{\mathcal{K}(K,B)}="tl";
 (36,16)*+{\mathcal{K}(K',B)}="tr";
 (-36,-14)*+{U^{\T}_{A'}F^{T}_{A'}\mathcal{K}(K,B) =\mathcal{K}(K,B)U^{\T}_{A}F^{T}_{A}}="bl";
 (36,-14)*+{ U^{\T}_{A'}F^{T}_{A'}\mathcal{K}(K',B) =
\mathcal{K}(K',B)U^{\T}_{A}F^{T}_{A}}="br";
        {\ar@{=>}^{ \mathcal{K}(k,B)} "tl";"tr"};
        {\ar@/_1.8pc/@{=>}_{ \mathcal{K}(K,B) i^{\T}_{A}} "tl";"bl"};
        {\ar@/^1.8pc/@{=>}^{i^{\T}_{A'}\mathcal{K}(K,B) } "tl";"bl"};
        {\ar@{=>}^{i^{\T}_{A'}\mathcal{K}(K',B)} "tr";"br"};
        {\ar@/^1.5pc/@{=>}^{U^{\T}_{A'}F^{T}_{A'}\mathcal{K}(k,B)} "bl";"br"};
        {\ar@/_1.5pc/@{=>}_{ \mathcal{K}(k,B)U^{\T}_{A}F^{T}_{A}} "bl";"br"};
            {\ar@3{->}^-{(i^{\T}_{A'}   \;_{\mathcal{K}(k,B)})^{-1}} (10,8.5);(2,-1.5)};
            {\ar@3{->}_-{  1} (-34,0);(-40,0)};
            {\ar@3{->}_-{  1} (0,-12);(0,-18)};
 \endxy
\]
\[ \qquad =
 \def\objectstyle{\scriptstyle}
 \def\labelstyle{\scriptscriptstyle}
 \xy
 (-36,16)*+{\mathcal{K}(K,B)}="tl";
 (36,16)*+{\mathcal{K}(K',B)}="tr";
 (-36,-14)*+{U^{\T}_{A'}F^{T}_{A'}\mathcal{K}(K,B) =\mathcal{K}(K,B)U^{\T}_{A}F^{T}_{A}}="bl";
 (36,-14)*+{ U^{\T}_{A'}F^{T}_{A'}\mathcal{K}(K',B) =
\mathcal{K}(K',B)U^{\T}_{A}F^{T}_{A}}="br";
        {\ar@{=>}^{\mathcal{K}(k,B)} "tl";"tr"};
        {\ar@{=>}_{\mathcal{K}(K,B) i^{\T}_{A} } "tl";"bl"};
        {\ar@/^1.8pc/@{=>}^{i^{\T}_{A'}\mathcal{K}(K',B)} "tr";"br"};
        {\ar@/_1.85pc/@{=>}_{ \mathcal{K}(K',B)i^{\T}_{A}} "tr";"br"};
        {\ar@{=>}_{ \mathcal{K}(k,B)U^{\T}_{A}F^{T}_{A}} "bl";"br"};
            {\ar@3{->}_-{\mathcal{K}(k,B)_{i^{\T}_A}} (-2,4.5);(-10,-2.5)};
            {\ar@3{->}_-{1 } (40,0);(34,0)};
 \endxy
\]

In order to show that the collection of $I^{\T}_A$ define a
\cat{Gray}-perturbation we must establish the following equality:
\[
 \xy
 (-20,14)*+{ sK}="tl";
 (20,14)*+{  SK'}="tr";
 (-20,-14)*+{SK }="bl";
 (20,-14)*+{ SK' }="br";
        {\ar@/^1.5pc/@{=>}^{sk} "tl";"tr"};
        {\ar@/_1.5pc/@{=>}_{sk'} "tl";"tr"};
        {\ar@/_1.8pc/@{=>}_{\nu K.\eta sK } "tl";"bl"};
        {\ar@/^1.8pc/@{=>}^{ sK} "tl";"bl"};
        {\ar@{=>}^{ sK'} "tr";"br"};
        {\ar@{=>}_{sk' } "bl";"br"};
            {\ar@{=} (6,-2);(4,-4)};
            {\ar@3{->}_-{  s\kappa} (3,16);(3,10)};
            {\ar@3{->}_-{  \psi^{-1}K } (-15,0);(-25,0)};
 \endxy
\quad = \quad
 \xy
 (-20,14)*+{ sK}="tl";
 (20,14)*+{  SK'}="tr";
 (-20,-14)*+{SK }="bl";
 (20,-14)*+{ SK' }="br";
        {\ar@/^1.5pc/@{=>}^{sk} "bl";"br"};
        {\ar@/_1.5pc/@{=>}_{sk'} "bl";"br"};
        {\ar@/_1.8pc/@{=>}_{\nu K'.\eta sK' } "tr";"br"};
        {\ar@/^1.8pc/@{=>}^{ sK'} "tr";"br"};
        {\ar@{=>}_{\nu K.\eta sK} "tl";"bl"};
        {\ar@{=>}^{sk } "tl";"tr"};
            {\ar@3{->}_-{k_{\nu.\eta s}^{-1}} (-2,8.5);(-10,1.5)};
            {\ar@3{->}_-{  s\kappa} (3,-10);(3,-16)};
            {\ar@3{->}_-{  \psi^{-1}K' } (25,0);(15,0)};
 \endxy
\]
This clearly follows from the \cat{Gray} axioms. To establish that
the collection of $E^{\T}_A$ define a \cat{Gray}-perturbation we
must verify the following equality of pseudoalgebra
2-homomorphisms:
\[
 \xy
 (-20,14)*+{ TsK}="tl";
 (20,14)*+{ TsK'}="tr";
 (-20,-14)*+{TsK}="bl";
 (20,-14)*+{ TsK'}="br";
        {\ar@/^1.5pc/@{=>}^{Tsk} "tl";"tr"};
        {\ar@/_1.5pc/@{=>}_{Tsk'} "tl";"tr"};
        {\ar@/_1.8pc/@{=>}_{1} "tl";"bl"};
        {\ar@/^1.8pc/@{=>}^{\mu sK.T\eta sK} "tl";"bl"};
        {\ar@{=>}^{ \mu sK'.T\eta sK'} "tr";"br"};
        {\ar@{=>}_{Tsk' } "bl";"br"};
            {\ar@3{->}^-{k_{\nu.\eta s}^{-1}} (10,-1.5);(2,-8.5)};
            {\ar@3{->}_-{  Ts\kappa} (3,16);(3,10)};
            {\ar@3{->}_-{  \rho^{-1}sK } (-15,0);(-25,0)};
 \endxy
\quad = \quad
 \xy
 (-20,14)*+{ TsK}="tl";
 (20,14)*+{ TsK'}="tr";
 (-20,-14)*+{TsK}="bl";
 (20,-14)*+{ TsK'}="br";
        {\ar@/^1.5pc/@{=>}^{Tsk} "bl";"br"};
        {\ar@/_1.5pc/@{=>}_{Tsk'} "bl";"br"};
        {\ar@/_1.8pc/@{=>}_{1} "tr";"br"};
        {\ar@/^1.8pc/@{=>}^{\mu sK'.T\eta sK'} "tr";"br"};
        {\ar@{=>}_{1} "tl";"bl"};
        {\ar@{=>}^{Tsk } "tl";"tr"};
            {\ar@{=} (-6,2);(-4,4)};
            {\ar@3{->}_-{  Ts\kappa} (3,-10);(3,-16)};
            {\ar@3{->}_-{  \rho^{-1}sK' } (25,0);(15,0)};
 \endxy
\]
Again this equality above is follows directly from the \cat{Gray}
axioms. The fact that this  data defines a pseudoadjunction in the
$\cat{Gray}$-category $[\mathcal{K}^{\op},\cat{Gray}]$ follows
from Proposition~\ref{FUApseudoadjoints} since the coherence
conditions are verified pointwise. \qed

Note that the previous theorem can also be adapted to the context
of a pseudocomonad $\G$ on $B$.  In this case, one obtains a
$\cat{Gray}$-functor $\cat{$\G$-CoAlg}\maps \mathcal{K}^{\op} \to
\cat{Gray}$.  As before there exists a forgetful
$\cat{Gray}$-natural transformation $U^{\G}\maps \cat{$\G$-CoAlg}
\to \mathcal{K}(-,B)$.  However, in this case, $U^{\G}$ has a
right pseudoadjoint $F^{\G}$.

\subsection{Pseudoadjoint pseudomonads} \label{pseudoadjsec}

In Section~\ref{secEMcompletion} it was explained how thinking of
an Eilenberg-Moore object as a weighted limit can be used to
construct the free completion of a 2-category under
Eilenberg-Moore objects. In this section we will need to
generalize the notion of an Eilenberg-Moore object for a monad to
an Eilenberg-Moore object for a pseudomonad. It turns out that
thinking of an Eilenberg-Moore object as weighted limit will prove
useful for this task as well.

Denote the `free living monad' or `walking monad' as $\cat{mnd}$,
meaning that a monad in a 2-category $\mathcal{K}$ is a 2-functor
$\cat{mnd} \to \mathcal{K}$.  Street has constructed a 2-functor
$J \maps \cat{mnd} \to \cat{Cat}$ with the property that the
Eilenberg-Moore object of a monad $\T \maps \cat{mnd} \to
\mathcal{K}$ is the $J$-weighted limit of the 2-functor $\T$
\cite{Street3}. This idea was used by Lack to construct a
$\cat{Gray}$-category $\cat{psm}$
--- the `free living pseudomonad' --- such that a pseudomonad $\T$
in the $\cat{Gray}$-category $\mathcal{K}$ is a
$\cat{Gray}$-functor $\T \maps \cat{psm} \to \mathcal{K}$
\cite{Lack}. Lack also constructs a $\cat{Gray}$-functor $P \maps
\cat{psm} \to \cat{Gray}$ with the property that the
Eilenberg-Moore object of the pseudomonad $\T$ is the $P$-weighted
limit of the $\cat{Gray}$-functor $\T \maps \cat{psm} \to
\mathcal{K}$, denoted $\{P,\T\}$.

This limit does not always exist in $\mathcal{K}$, but
$\mathcal{K}$ can always be embedded via the Yoneda embedding $Y
\maps \mathcal{K} \to [\mathcal{K}^{\op},\cat{Gray}]$ where the
$P$-weighted limit of $Y\T$ can be formed.  Then the
Eilenberg-Moore object of $\T$ will exist in $\mathcal{K}$ if and
only if $\{P,Y\T\}$ is representable since the Yoneda embedding
must preserve any limits which exist. The $\cat{Gray}$-functor
$\{P,Y\T\}$ is just the $\cat{Gray}$-functor $\cat{$\T$-Alg}$
constructed in the previous section. Thus, an Eilenberg-Moore
object for the pseudomonad $\T$ is just a choice of representation
for $\cat{$\T$-Alg}$.

If $\cat{$\T$-Alg}$ is representable, then in our previous
notation, it will correspond to the $\cat{Gray}$-functor
$\mathcal{K}(-,B^{\T})$ where $B^{\T}$ is an Eilenberg-Moore
object for the pseudomonad $\T$.  If an Eilenberg-Moore object for
$\T$ does exist then the pseudoadjunction
$I^{\T},E^{\T},i^{\T},e^{\T}\maps F^{\T} \dashv_p U^{\T}\maps
\cat{$\T$-Alg} \to \mathcal{K}(-,B)$ of Theorem~\ref{FextendsGray}
corresponds via the enriched Yoneda lemma to a pseudoadjunction
$I,E,i,e\maps F \dashv_p U \maps B \to B^{\T}$ in $\mathcal{K}$.

The limit description of an Eilenberg-Moore object in a
$\cat{Gray}$-category also facilitates the free completion of an
arbitrary $\cat{Gray}$-category to one that has Eilenberg-Moore
objects. Indeed, because a $\cat{Gray}$-category is just a
$\cat{Gray}$-enriched category, the free-completion is achieved
using the theory of enriched category theory~\cite{kel}.  With all
of the abstract theory in place, we begin by proving the pointwise
version of the categorified adjoint monad theorem.  In
Theorem~\ref{pseudoADJMONTHM} we will prove the full result.

\begin{lem} \label{lemmar}
If $\Sigma,\Upsilon,\iota,\sigma\maps T \dashv_p G$, then the
2-category \cat{$\T$-Alg}$_{\bf A}$ of pseudoalgebras based at $A$
is 2-equivalent to the 2-category \cat{$\G$-CoAlg}$_{\bf A}$ of
pseudocoalgebras based at $A$ for the pseudocomonad $\G$.
Furthermore, this 2-equivalence commutes with the forgetful
2-functors $U_A^{\T} \maps \cat{$\T$-Alg}_{\bf A} \to
\mathcal{K}(A,B)$ and $U^{\G}_A \maps \cat{$\G$-CoAlg}_{\bf A} \to
\mathcal{K}(A,B)$.
\end{lem}

\Proof  This lemma is essentially due to the properties of
pseudomates under pseudoadjunction and the fact that this
association preserves composites up to natural isomorphism. With
$\Theta$ and $\Phi$ as in Proposition~\ref{pseudomates}, define
the 2-functor:
 \ban
\mathcal{M}_A\maps \cat{$\T$-Alg}_{\bf A} &\to&
\cat{$\G$-CoAlg}_{\bf A} \\
    (s, \nu, \psi, \chi)
& \mapsto &
    \big(s, \Phi(\nu),
  \xy
    (-25,0)*+{\Phi(\eta)s.\Phi(\nu)}="1";
    (-0,0)*+{\Phi(\nu.\eta s)}="2";
    (25,0)*+{\Phi(s)=s}="3";
        {\ar@3{->}^-{\Phi(\psi)} "2";"3"};
        {\ar@3{->}^-{Y} "1";"2"};
  \endxy, \\
&  &
  \xy
    (-25,0)*+{G\Phi(\nu).\Phi(\nu)}="1";
    (-0,0)*+{\Phi(\nu.T\nu)}="2";
    (25,0)*+{\Phi(\nu.\mu s)}="3";
    (50,0)*+{\Phi(\nu)s.\Phi(\nu)}="4";
        {\ar@3{->}^-{\Phi(\chi)} "2";"3"};
        {\ar@3{->}^-{Y^{-1}} "3";"4"};
        {\ar@3{->}^-{X} "1";"2"};
  \endxy \big) \\
 (h,\varrho) & \mapsto & (h,
 \xy
    (-25,0)*+{\Phi(\nu').\Phi(h)}="1";
    (3,0)*+{\Phi(\nu'.Th)}="2";
    (28,0)*+{\Phi(h.\nu)}="3";
    (50,0)*+{Gh.\Phi(\nu)}="4";
        {\ar@3{->}^-{\Phi(\varrho)} "2";"3"};
        {\ar@3{->}^-{X^{-1}} "3";"4"};
        {\ar@3{->}^-{X} "1";"2"};
  \endxy) \\
 \xi \maps
(h,\varrho) \To (h',\varrho') & \mapsto & \xi \maps
(h,X^{-1}\circ\Phi(\varrho)\circ X) \To
(h',X^{-1}\circ\Phi(\varrho')\circ X)
 \ean
This data defines a pseudocoalgebra, morphism of pseudocoalgebras,
and 2-morphism of pseudocoalgebras because $\xi \maps (h,\varrho)
\To (h',\varrho')\maps (s,\nu,\psi,\chi) \to
(s',\nu',\psi',\chi')$ is a 2-morphism of pseudoalgebras, and
mateship under pseudoadjunction preserves all composites up to
natural isomorphism.

Since $\mathcal{M}_A \maps h\maps s \To s' \mapsto h \maps s \To
s'$ it is clear that $\mathcal{M}_A$ preserves 1-morphism
identities and to see that $\mathcal{M}_A$ preserves composites of
1-morphisms all we must check is its behavior on $\varrho$.  For
this purpose it will be helpful to have the specific form of
$\mathcal{M}_A(h,\varrho)$. By plugging in the relevant
pseudoadjunctions, one can check that $\mathcal{M}_A(h,\varrho) =
\big(h, \Phi(\varrho) \circ G\nu'.(\iota_h)\big)$. Hence, using
the definition of $\Phi(\varrho)$ and the \cat{Gray} axioms it is
easy to verify the following chain of equalities:
 \ban
    \mathcal{M}_A(h',\varrho').\mathcal{M}_A(h,\varrho)
& = &
    \big( h'.h,Gh'.G\varrho.\iota s \circ Gh'.G\nu'.\iota_{h}\circ
   G \varrho'.\iota s'.h\circ G\nu''.\iota_{h'}.h
             \big)
\\
  & = &  \big(h'.h,G h'.G\varrho.\iota s \circ G \varrho'.GTh.\iota s \circ
G\nu''.\iota_{h'.h}    \big) \\
  & = &
\mathcal{M}_A(h'.h,\varrho'.Th\circ h'.\varrho).
 \ean

Since $\mathcal{M}_A$ maps 2-morphisms to themselves, it is clear
that $\mathcal{M}_A$ preserves composition of 2-morphisms on the
nose as well.  Hence, $\mathcal{M}_A\maps \cat{$\T$-Alg}_A \to
\cat{$\G$-CoAlg}_A$ is a 2-functor. Now we define the other
2-functor taking part in the equivalence:
 \ban
\overline{\mathcal{M}}_A\maps \cat{$\G$-CoAlg}_{\bf A} &\to&
\cat{$\T$-Alg}_{\bf A}
 \\
    (s, \bar{\nu}, \bar{\psi}, \bar{\chi})
 & \mapsto &
    \big(s, \Theta(\bar{\nu}),
 \xy
    (-25,0)*+{\Theta(\bar{\nu}).\Theta(\varepsilon)s}="1";
    (3,0)*+{\Theta(\varepsilon s.\bar{\nu})}="2";
    (28,0)*+{\Theta(s) = s}="3";
        {\ar@3{->}^-{\Theta(\bar{\psi})} "2";"3"};
        {\ar@3{->}^-{W} "1";"2"};
  \endxy,
\\
 & &
\xy
    (-25,0)*+{\Theta(\bar{\nu}).T\Theta(\bar{\nu})}="1";
    (0,0)*+{\Theta(G\bar{\nu}.\bar{\nu})}="2";
    (25,0)*+{\Theta(\delta s.\bar{\nu})}="3";
    (50,0)*+{\Theta(\bar{\nu}).\Theta(\delta)s}="4";
        {\ar@3{->}^-{\Theta(\bar{\chi})} "2";"3"};
        {\ar@3{->}^-{W^{-1}} "3";"4"};
        {\ar@3{->}^-{V} "1";"2"};
  \endxy)
 \\
    (h,\bar{\varrho})
 & \mapsto &
    (h,
 \xy
    (-25,0)*+{\Theta(\bar{\nu}').Th}="1";
    (-3,0)*+{\Theta(\bar{\nu}.h)}="2";
    (19,0)*+{\Theta(Gh.\bar{\nu})}="3";
    (53,0)*+{\Theta(h).\Theta(\bar{\nu}) = h.\Theta(\bar{\nu})}="4";
        {\ar@3{->}^-{\Theta(\bar{\varrho})} "2";"3"};
        {\ar@3{->}^-{V^{-1}} "3";"4"};
        {\ar@3{->}^-{V} "1";"2"};
  \endxy)
 \\
 \xi \maps
(h,\bar{\varrho}) \To (h',\bar{\varrho}') & \mapsto & \xi \maps
(h, V^{-1} \circ\Theta(\bar{\varrho}) \circ V) \To (h', V^{-1}
\circ\Theta(\bar{\varrho'}) \circ V)
 \ean
This will define a 2-functor, again by the functoriality of
mateship under pseudoadjunction and the axioms of
$\cat{Gray}$-category.  It will be helpful to have the explicit
formula for $\overline{\mathcal{M}}_A(h,\bar{\varrho})$.  By
plugging in the relevant pseudoadjunctions one can check that
$$\overline{\mathcal{M}}_A(h,\bar{\varrho})=(h, \sigma_h.T\bar{\nu}
\circ \Theta(\bar{\varrho})).$$

We now show that $\mathcal{M}_A$ and $\overline{\mathcal{M}}_A$
define a 2-equivalence of 2-categories. Define the 2-natural
isomorphism $\Gamma_{A}\maps1_{\cat{$\T$-Alg}_{\bf A}} \To
\overline{\mathcal{M}}_A\mathcal{M}_A$ as follows: Denote
$\overline{\mathcal{M}}_A\mathcal{M}_A\big((s,\nu,\psi,\chi)\big)$
as $(s,\tilde{\nu},\tilde{\psi},\tilde{\chi})$. Define the
morphism of pseudoalgebras $\Gamma_{(s,\nu,\psi,\chi)} \maps
(s,\nu,\psi,\chi)\to (s,\tilde{\nu} ,\tilde{\psi}, \tilde{\chi})$
by letting $h \maps s \To s$ be the identity, so that $\varrho$ is
just a map $\tilde{\nu} \To \nu$.  From the definition of
$\mathcal{M}_A$ and $\overline{\mathcal{M}}_A$ we know that
$\tilde{\nu}=\Theta\Phi(\nu)$.  Hence we can choose $\varrho$ to
be the isomorphism $\bar{\gamma}_{\nu}^{-1} \maps \Theta\Phi(\nu)
\To \nu$ defined in Proposition~\ref{pseudomates}.  The pair
$(1_s,\bar{\gamma}_{\nu}^{-1})$ is a morphism of pseudoalgebras by
the naturality of the isomorphism $\bar{\gamma}$ of Proposition
\ref{pseudomates} applied to the 3-morphisms $\psi$ and $\chi$.
The explicit form of the isomorphism $\bar{\gamma}_{\nu}^{-1}$ is
$\nu.\Upsilon s \circ \sigma_{\nu}^{-1}.T\iota s  $.

To see that  $\Gamma_{A}$ is natural in the one dimensional sense,
suppose that $(h,\varrho) \maps (\psi,\chi) \to (\psi',\chi')$ is
an arbitrary 1-cell in \cat{$\T$-Alg}$_{\bf A}$. Consider the
following diagram:
\[
 \xymatrix{
 (s,\nu,\psi,\chi)
   \ar[rr]^{(h,\varrho)}
   \ar[d]_{(1,\bar{\gamma}_{\nu}^{-1})}
&&(s',\nu',\psi',\chi')
   \ar[d]^{(1,\bar{\gamma}_{\nu'}^{-1})} \\
 (s,\tilde{\nu},\tilde{\psi},\tilde{\chi})
     \ar[rr]_{\overline{\mathcal{M}}_A\mathcal{M}_A(h,\varrho)}
 &&(s',\tilde{\nu}',\tilde{\psi}',\tilde{\chi}')
 }
\]
Note that since $h \maps s \To s'$ for some $s' \maps A \to B$,
the pseudoadjunction determining the mate of $h$ is the identity
adjunction so that $h$ is its own mate under pseudoadjunction.
Thus, this diagram of pseudoalgebra maps commutes if and only if
 \be \label{ffvc}
h.\gamma_{\nu}\circ\tilde{\varrho}  = \varrho\circ\gamma_{\nu'}.Th
.
 \ee

Using the explicit formulas given above we have that
$$
\overline{\mathcal{M}}_A\mathcal{M}_A(h,\varrho) =
(h,\sigma_h^{-1}.TG\nu.T\iota s \circ \sigma s'.TG\varrho.T\iota
s\circ \sigma s'.TG\nu'.T\iota_h^{-1}).
$$
In order to prove the naturality of $\Gamma_A$ we will need the
following equalities that all follow directly from the axioms of a
\cat{Gray}-category:
\begin{eqnarray*}
  h.\sigma_{\nu}^{-1} \circ \sigma_h^{-1}.TG\nu  &=& \sigma_{h.\nu}^{-1} \\
  \sigma_{h.\nu}^{-1}\circ \sigma s'.TG\varrho &=& \varrho.\sigma Ts \circ \sigma_{\nu'.Th}^{-1} \\
  \sigma_{\nu'.Th}^{-1} &=& \nu'.\sigma_{Th}^{-1} \circ \sigma_{\nu'}^{-1}.TGTh \\
  \sigma_{Th}^{-1}.T\iota s \circ \sigma Ts'.T \iota_h^{-1} \ &=& (\sigma.T\iota)_h^{-1} \\
  Th.\Upsilon s \circ (\sigma.T\iota)_h^{-1} &=& Th_s \circ \Upsilon s'.Th=
\Upsilon s'.Th
\end{eqnarray*}
The proof of equation \ref{ffvc} above is as follows:
\begin{eqnarray*}
  h.\gamma_{\nu}\circ\tilde{\varrho}
   &=&
    h.\nu.\Upsilon s \circ h.\sigma_{\nu}^{-1}.T\iota s \circ
\sigma_{h}^{-1}.TG\nu.T\iota s \circ \sigma s'.TG \varrho.T\iota
s\circ \sigma s'.TG\nu'.T\iota_h^{-1}
    \\
   &=&
  h.\nu.\Upsilon s \circ \sigma_{h.\nu}^{-1}.T \iota s \circ \sigma s'.TG \varrho.T\iota
s\circ \sigma s'.TG\nu'.T\iota_h^{-1}
 \\
   &=&
  h.\nu.\Upsilon s \circ \varrho.\sigma Ts.T \iota s \circ
\sigma_{\nu'.Th}^{-1}.T \iota s \circ \sigma s'.TG
\nu'.T\iota_h^{_1}
 \\
   &=&
 h.\nu.\Upsilon s \circ \varrho.\sigma Ts.T \iota s \circ \nu'.\sigma
Ts'.T\iota_h^{-1} \circ \sigma_{\nu'}^{-1}.T\iota s'.Th \quad
(Interchange)
 \\
   &=&
 h.\nu.\Upsilon s \circ \varrho.\sigma Ts.T \iota s\circ
\nu'.(\sigma.T\iota)_h^{-1}\circ \sigma_{\nu'}^{-1}.T\iota s'.Th
 \\
   &=&
\varrho \circ \nu'.Th.\Upsilon s \circ
\nu'.(\sigma.T\iota)_h^{-1}\circ \sigma_{\nu'}^{-1}.T\iota s'.Th
 \quad (Interchange)
 \\
   &=&
 \varrho \circ \nu'.\Upsilon s'.Th \circ \sigma_{\nu'}^{-1}.T\iota s'.Th
\\
   &=&
\varrho\circ\gamma_{\nu'}.Th
\end{eqnarray*}

To see the 2-naturality of $\Gamma_{A}$ let $\xi \maps (h,\rho)
\To (h',\varrho')$, then the equality:
\[
 \xy
   (-25,0)*+{(\psi,\chi)}="l";
   (0,0)*+{(\psi',\chi')}="m";
   (25,0)*+{(\tilde{\psi},\tilde{\chi})}="r";
        {\ar@/^1.3pc/^{(h,\varrho)} "l";"m"};
        {\ar@/_1.3pc/_{(h',\varrho')}  "l";"m"};
        {\ar^{(1,\bar{\gamma}^{-1}_{\nu'})}  "m";"r"};
        {\ar@{=>}^{\xi} (-12.5,3);(-12.5,-3) };
 \endxy
\qquad = \hspace{2in}
\]
\[ \hspace{2in}
 \xy
   (-25,0)*+{(\psi,\chi)}="l";
   (0,0)*+{(\psi',\chi')}="m";
   (30,0)*+{(\tilde{\psi},\tilde{\chi}) .}="r";
        {\ar@/^1.3pc/^{\overline{\mathcal{M}}_A\mathcal{M}_A(h,\varrho)} "m";"r"};
        {\ar@/_1.3pc/_{\overline{\mathcal{M}}_A\mathcal{M}_A(h',\varrho')}  "m";"r"};
        {\ar^{(1,\bar{\gamma}^{-1}_{\nu})}  "l";"m"};
        {\ar@{=>}^{\overline{\mathcal{M}}_A\mathcal{M}_A\xi} (9.5,3);(9.5,-3) };
 \endxy
\]
follows from the fact that
$\overline{\mathcal{M}}_A\mathcal{M}_A(\xi) = \xi$ and the
naturality of $\bar{\gamma}$ applied to the 3-morphism $\xi$ in
$\mathcal{K}$. A 2-natural isomorphism $\overline{\Gamma}_A \maps
\mathcal{M}_A\overline{\mathcal{M}}_A \To 1_{\cat{$\G$-CoAlg}_{\bf
A}}$ can be defined in a similar way.

To see that this 2-equivalence of 2-categories commutes with the
forgetful 2-functors, note that in the above proof we have shown
that the 2-equivalence is the identity on the base map $s$ of the
pseudoalgebra.  Furthermore, in the discussion of naturality we
have shown that for any map $(h,\varrho)$ of pseudoalgebras
$\mathcal{M}_A$ is the identity on $h$ and $\mathcal{M}_A$ also
acts as the identity on every 3-cell defining a 2-morphism of
pseudoalgebras. Thus, by the definition of the forgetful
2-functors $U_A^{\T}$ and $U^{\G}_{A}$ it is clear that the
equivalence $\mathcal{M}_A$ commutes with the forgetful
2-functors. \qed

\begin{thm}[The categorified adjoint monad theorem] \label{pseudoADJMONTHM}
If $\Sigma,\Upsilon,\iota,\sigma\maps T \dashv_p G$ in
$\mathcal{K}$, then the $\cat{Gray}$-functors $\cat{$\T$-Alg}$ and
$\cat{$\G$-Alg}$ are \cat{Gray}-equivalent in the
$\cat{Gray}$-category $[\mathcal{K}^{\op},\cat{Gray}]$.  This
means that there exists \cat{Gray}-natural transformations $
\mathcal{M}\maps \cat{$\T$-Alg} \to \cat{$\G$-CoAlg}$, $
\overline{\mathcal{M}}\maps \cat{$\G$-CoAlg} \to \cat{$\T$-Alg}$
and invertible $\cat{Gray}$-modifications
$\Gamma\maps1_{\cat{$\T$-Alg}} \To
\overline{\mathcal{M}}\mathcal{M}$, $\overline{\Gamma} \maps
\mathcal{M}\overline{\mathcal{M}} \To 1_{\cat{$\G$-CoAlg}}$.
Furthermore, this \cat{Gray}-equivalence commutes with the
forgetful $\cat{Gray}$-natural transformations $U^{\T}$ and
$U^{\G}$.
\end{thm}

\Proof Define a \cat{Gray}-natural transformation
$\mathcal{M}\maps \cat{$\T$-Alg} \to \cat{$\G$-CoAlg}$ which
assigns to each object in $\mathcal{K}^{\op}$ the 2-functor
$\mathcal{M}_A$ defined in the preceding lemma. To see the
naturality of $\mathcal{M}$ let $K \maps A' \to A$ in
$\mathcal{K}$ and note that $\Phi(fK)=\Phi(f)K$ for
$f=\nu,\psi,\chi$ so that:
 \ban
\hat{K}\mathcal{M}_A(s,\nu,\psi,\chi) &=& \big( sK, \Phi(\nu)K,
\Phi(\psi)K\circ
YK,Y^{-1}K\circ \Phi(\chi)K\circ XK \big) \\
&=& \big( sK, \Phi(\nu K), \Phi(\psi K)\circ Y,Y^{-1}\circ
\Phi(\chi K)\circ X \big) \\
&=& \mathcal{M}_{A'}\hat{K}(s,\nu,\psi,\chi).
 \ean
A similar check shows that the collection of 2-functors
$\overline{\mathcal{M}}_A \maps \cat{$\G$-CoAlg}_A \to
\cat{$\T$-Alg}_A$ defines a \cat{Gray}-natural transformation
$\overline{\mathcal{M}} \maps \cat{$\G$-CoAlg} \to
\cat{$\T$-Alg}$.

Define an invertible \cat{Gray}-modification $\Gamma \maps
1_{\cat{$\T$-Alg}} \to \overline{\mathcal{M}}\mathcal{M}$ which
assigns to each object in $\mathcal{K}^{\op}$ the 2-natural
isomorphism $\Gamma_A$ defined in the preceding lemma.  To see
that this data defines a \cat{Gray}-modification first note that
the following diagrams commute:
\[
 \xy
 (-20,10)*+{  \hat{K} }="tl";
 (20,10)*+{   \hat{K}'}="tr";
 (-20,-10)*+{ \overline{\mathcal{M}}_{A'}\mathcal{M}_{A'}\hat{K} }="bl";
 (20,-10)*+{  \overline{\mathcal{M}}_{A'}\mathcal{M}_{A'}\hat{K}'}="br";
        {\ar@{=>}^{ \hat{k}} "tl";"tr"};
        {\ar@{=>}_{ \Gamma_{A'}\hat{K}} "tl";"bl"};
        {\ar@{=>}^{ \Gamma_{A'}\hat{K}'} "tr";"br"};
        {\ar@{=>}_{ \overline{\mathcal{M}}_{A'}\mathcal{M}_{A'}\hat{k}} "bl";"br"};
 \endxy
\quad  \quad
 \xy
 (-20,10)*+{  \hat{K} }="tl";
 (20,10)*+{\hat{K}' }="tr";
 (-20,-10)*+{ \hat{K}\overline{\mathcal{M}}_{A}\mathcal{M}_{A}}="bl";
 (20,-10)*+{  \hat{K}'\overline{\mathcal{M}}_{A}\mathcal{M}_{A}}="br";
        {\ar@{=>}^{  \hat{k} } "tl";"tr"};
        {\ar@{=>}_{ \hat{K} \Gamma_{A}} "tl";"bl"};
        {\ar@{=>}^{ \hat{K}'\Gamma_{A}} "tr";"br"};
        {\ar@{=>}_{ \hat{k}\overline{\mathcal{M}}_{A}\mathcal{M}_{A} } "bl";"br"};
 \endxy
\]
In the first case we must show that the composites of
pseudoalgebra homomorphisms $(sK',\bar{\gamma}_{\nu
K'}^{_1}).(sk,\nu_k^{-1})$ and
$(sk,\Theta\Phi(\nu)_k^{-1}).(sK,\bar{\gamma}_{\nu K}^{-1})$ are
equal.  In the second diagram we must show that
$(sK',\bar{\gamma}^{-1}_{\nu}K').(sK,\nu_k^{-1}) =
(sK',\Theta\Phi(\nu)_k^{-1}).(sK,\bar{\gamma}^{-1}_{\nu}K)$. These
equalities follow from the fact that $\Theta\Phi(\nu
K')=\Theta\Phi(\nu)K'$, $\Theta\Phi(\nu K)=\Theta\Phi(\nu)K$ and
the \cat{Gray} axioms.  Finally, $\Gamma$ is a
\cat{Gray}-modification since
$\bar{\gamma}^{-1}_{\nu}K'=\bar{\gamma}^{-1}_{\nu K'}$ and
$\bar{\gamma}^{-1}_{\nu}K=\bar{\gamma}^{-1}_{\nu K}$.

In a similar fashion one can define an invertible
\cat{Gray}-modification $\overline{\Gamma} \maps
\mathcal{M}\overline{\mathcal{M}} \To 1_{\cat{$\G$-CoAlg}}$.
Hence, we have shown that the $\cat{Gray}$-functors \cat{$\T$-Alg}
and \cat{$\G$-CoAlg} are \cat{Gray}-equivalent in the
$\cat{Gray}$-category $[\mathcal{K}^{\op},\cat{Gray}]$.
Lemma~\ref{lemmar} shows that this \cat{Gray}-equivalence commutes
with the forgetful functors. \qed

\begin{thm} \label{THMab}
Let $\T = (T,\mu,\eta,\lambda,\rho,\alpha,\varepsilon)$ be a
Frobenius pseudomonad on the object $B$ of the
$\cat{Gray}$-category $\mathcal{K}$ with
$\Sigma,\Upsilon,\iota,\varepsilon.\mu\maps T \dashv_p T$.  Denote
the induced pseudocomonad structure on $T$ by $\G$. Then there
exists an invertible $\cat{Gray}$-modification $\Xi \maps
\mathcal{M}F^{\T}\to F^{\G}$.  Alternatively, the morphisms
$\mathcal{M}F^{\T}$ and $F^{\G}$ of the $\cat{Gray}$-category
$[\mathcal{K}^{\op},\cat{Gray}]$ are isomorphic.
\end{thm}

\Proof We begin by defining a 2-natural isomorphism $\Xi_A \maps
\mathcal{M}_AF_A^{\T} \To F_A^{\G} \maps
 \mathcal{K}(A,B) \to \cat{$\G$-CoAlg}_A$.  Recall that
$F_A^{\T}(s)=(Ts,\mu s, \lambda s, \alpha s)$ and that
$F_A^{\G}(s)=(Ts,\delta s, \bar{\rho}^{-1}s, \bar{\alpha}s)$ with
$\delta,\bar{\rho},\bar{\alpha}$ of the form given in
Proposition~\ref{PROPaa}.  Hence,
 \ban
 \mathcal{M}_AF_A^{\T}(s)
& =& \big(Ts,\Phi\left(\left(\mu
 s\right)\right),\Phi\left(\left(\lambda s\right)\right),\Phi\left(\left(\alpha s\right)\right) \big)  \\
 & =& \big(Ts,T\mu s.\iota s, \Phi\left(\left(\lambda s\right)\big),
 \Phi\left(\left(\alpha s\right)\right)\right) \\
 F^{\G}_A(s) & = & \left(Gs,\delta s, \bar{\rho}^{-1}s,\bar{\alpha} s\right) =
    \left(Ts,T^2(\varepsilon.\mu)s.T^2\mu Ts.T\iota T^2s.\iota Ts,\bar{\rho}^{-1}s,\bar{\alpha} s
    \right).
 \ean
The double parenthesis are used in order to distinguish which
pseudoadjunction is intended. For example, $\Phi(\mu s)$ is
determined by the pseudoadjunctions with $F=T,U=T,F'=T^2,U'=G^2$
and morphisms $a=b=s$.  While $\Phi\left(\left(\mu
s\right)\right)$ is given by the pseudoadjunctions with
$F=U=1_{B},F'=T,U'=T$ and morphisms $a=b=Ts$.

We define the isomorphism of pseudocoalgebras
$(h,\bar{\varrho}s)\maps\mathcal{M}_AF_A^{\T}(s) \to F^{\G}_A(s)$
by taking $h=1_{Ts}$, and $\bar{\varrho}s$ given by the following
diagram:
\[
 \xy
   (-60,0)*+{Ts}="1";
   (-40,10)*+{T^3s}="2t";
   (-40,-10)*+{T^3 s }="2b";
   (-20,10)*+{T^2s}="3t";
   (-20,-10)*+{T^5 s }="3b";
    (0,0)*+{T^4s}="4t";
   (0,-20)*+{T^4 s }="4b";
   (40,10)*+{T^2s}="5t";
   (20,-10)*+{T^3 s }="5b";
        {\ar@{=>}^{\iota Ts  } "1";"2t"};
        {\ar@{=>}_{\iota Ts  } "1";"2b"};
        {\ar@{=}_{} "2b";"2t"};
        {\ar@{=>}^{T \mu s  } "2t";"3t"};
        {\ar@{=>}_{T \iota T^2s } "2b";"3b"};
        {\ar@{=>}^{T \iota Ts } "3t";"4t"};
        {\ar@{=>}^{T^3 \mu Ts } "3b";"4t"};
        {\ar@{=>}_{T^2 \mu Ts } "3b";"4b"};
        {\ar@{=>}^{T^2s} "3t";"5t"};
        {\ar@{=>}^{T^2\mu s } "4t";"5b"};
        {\ar@{=>}_{T^2\mu s } "4b";"5b"};
        {\ar@{=>}_{T^2\varepsilon s } "5b";"5t"};
            {\ar@3{->}^{T\Sigma s} (15,6);(15,0)}
            {\ar@3{->}^{T^2\alpha s} (0,-6);(0,-12)}
            {\ar@3{->}^{T \iota_{\mu}^{-1} s} (-30,6);(-30,0)}
 \endxy
\]
One can verify using the definitions of $\mathcal{M}_A$,
$\bar{\psi}$, and $\bar{\chi}$ that this map is indeed a morphism
of pseudocoalgebras.

Let $h \maps s \To s'$ in $\mathcal{K}(A,B)$. To establish the
1-naturality of $\Xi_A$ we must show that:
\[
 \xy
   (-38,-10)*+{(Ts,\delta s,\bar{\rho}^{-1}s,\bar{\alpha}s)}="tr";
   (-38,10)*+{\mathcal{M}_A(Ts,\mu s, \lambda s, \alpha s )}="tl";
   (38,10)*+{\mathcal{M}_A(Ts',\mu s', \lambda s', \alpha s' )}="bl";
   (38,-10)*+{(Ts',\delta s',\bar{\rho}^{-1}s',\bar{\alpha}s')}="br";
        {\ar_-{(Ts,\bar{\varrho}s)} "tl";"tr"};
        {\ar^{\mathcal{M}_A(Th,\mu_h^{-1})} "tl";"bl"};
        {\ar_{(Th,\delta_h^{-1})} "tr";"br"};
        {\ar^-{(Ts',\bar{\varrho}{s'})} "bl";"br"};
 \endxy
\]
commutes where
 \ban
 \mathcal{M}_A(Th,\mu_h^{-1})
 &=&\Phi((\mu_h^{-1}))\circ T\mu s'.(\iota_{Th}^{-1})   \\
 &=&
 T \mu_{Th}^{-1}.\iota Ts\circ T\mu s'.(\iota_{Th}^{-1}).
 \ean
This amounts to the equality of the following diagrams:
\[
 \xy
   (-60,10)*+{Ts'}="1";
   (-40,10)*+{T^3s'}="2t";
   (-40,-10)*+{T^3 s' }="2b";
   (-10,10)*+{T^2s'}="3t";
   (-20,-10)*+{T^5 s' }="3b";
    (0,0)*+{T^4s'}="4t";
   (0,-20)*+{T^4 s' }="4b";
   (40,10)*+{T^2s'}="5t";
   (20,-10)*+{T^3 s' }="5b";
   (-40,25)*+{T^3 s }="5x";
   (40,25)*+{T^2s}="5Y";
   (-60,25)*+{Ts}="2x";
        {\ar@{=>}^{\iota Ts'  } "1";"2t"};
        {\ar@{=>}_{\iota Ts'  } "1";"2b"};
        {\ar@{=}_{} "2b";"2t"};
        {\ar@{=>}^{T \mu s'  } "2t";"3t"};
        {\ar@{=>}_{T \iota T^2s' } "2b";"3b"};
        {\ar@{=>}^{T \iota Ts' } "3t";"4t"};
        {\ar@{=>}^{T^3 \mu Ts' } "3b";"4t"};
        {\ar@{=>}_{T^2 \mu Ts' } "3b";"4b"};
        {\ar@{=} "3t";"5t"};
        {\ar@{=>}^{T^2h} "5Y";"5t"};
        {\ar@{=>}^{T^2\mu s' } "4t";"5b"};
        {\ar@{=>}_{T^2\mu s' } "4b";"5b"};
        {\ar@{=>}_{T^2\varepsilon s' } "5b";"5t"};
        {\ar@{=>}^{T^3 h } "5x";"2t"};
        {\ar@{=>}^{T \mu s } "5x";"5Y"};
        {\ar@{=>}_{T h } "2x";"1"};
        {\ar@{=>}^{\iota Ts } "2x";"5x"};
            {\ar@3{->}^{T\Sigma s'} (15,6);(15,0)}
            {\ar@3{->}^{T^2\alpha s'} (0,-6);(0,-12)}
            {\ar@3{->}^{T \iota_{\mu}^{-1} s'} (-30,6);(-30,0)}
            {\ar@3{->}^{\iota_{Th}} (-50,20);(-50,14)}
            {\ar@3{->}^{T\mu_{Th}^{-1}} (0,20);(0,14)}
 \endxy
\]
\[
 \xy
   (-60,-10)*+{Ts}="1";
   (-40,10)*+{T^3s}="2t";
   (-40,-10)*+{T^3 s }="2b";
   (-20,10)*+{T^2s}="3t";
   (-20,-10)*+{T^5 s }="3b";
    (0,0)*+{T^4s}="4t";
   (0,-20)*+{T^4 s }="4b";
   (40,-10)*+{T^2s}="5t";
   (20,-10)*+{T^3 s }="5b";
   (-60,-30)*+{Ts'}="1'";
   (40,-30)*+{T^2s'}="5'";
   (-40,-30)*+{T^3 s' }="2'";
   (-20,-30)*+{T^5 s' }="3'";
   (0,-40)*+{T^4 s' }="4'";
   (20,-30)*+{T^3 s' }="5b'";
        {\ar@{=>}^{\iota Ts  } "1";"2t"};
        {\ar@{=>}_{\iota Ts  } "1";"2b"};
        {\ar@{=}_{} "2b";"2t"};
        {\ar@{=>}^{T \mu s  } "2t";"3t"};
        {\ar@{=>}_{T \iota T^2s } "2b";"3b"};
        {\ar@{=>}^{T \iota Ts } "3t";"4t"};
        {\ar@{=>}^{T^3 \mu Ts } "3b";"4t"};
        {\ar@{=>}_{T^2 \mu Ts } "3b";"4b"};
        {\ar@{=>}@/^2.5pc/^{T^2s} "3t";"5t"};
        {\ar@{=>}^{T^2\mu s } "4t";"5b"};
        {\ar@{=>}_{T^2\mu s } "4b";"5b"};
        {\ar@{=>}_{T^2\varepsilon s } "5b";"5t"};
        {\ar@{=>}_{Th } "1";"1'"};
        {\ar@{=>}^{T^2h } "5t";"5'"};
        {\ar@{=>}_{\iota Ts'  } "1'";"2'"};
        {\ar@{=>}_{T \iota T^2s' } "2'";"3'"};
        {\ar@{=>}_{T^2 \mu Ts' } "3'";"4'"};
        {\ar@{=>}_{T^2\mu s' } "4'";"5b'"};
        {\ar@{=>}_{T^2\varepsilon s' } "5b'";"5'"};
            {\ar@3{->}^{T\Sigma s} (15,6);(15,0)}
            {\ar@3{->}^{T^2\alpha s} (0,-6);(0,-12)}
            {\ar@3{->}^{T \iota_{\mu}^{-1} s} (-30,6);(-30,0)}
            {\ar@3{->}^{\delta_h^{-1}} (0,-26);(0,-32)}
 \endxy
\]
which are equal by a routine verification using the
$\cat{Gray}$-category axioms. The 2-naturality of $\Xi_A$ follows
from the fact that both $F_A^{\T}$ and $F_A^{\G}$ map the
2-morphism $\xi \maps h \To h' \maps s \to s'$ to $T \xi$, and the
fact that the 2-functor $\mathcal{M}_A$ is the identity on
2-morphisms of pseudocoalgebras.

The collection of $\Xi_A$ define a $\cat{Gray}$-modification by
the commutativity of the following diagrams:
\[
 \xy
 (-28,10)*+{ \mathcal{M}_{A'}F^{\T}_{A'}\mathcal{K}(K,B) }="tl";
 (28,10)*+{  \mathcal{M}_{A'}F^{\T}_{A'}\mathcal{K}(K',B)}="tr";
 (-28,-10)*+{ F^{\G}_{A'}\mathcal{K}(K,B) }="bl";
 (28,-10)*+{ F^{\G}_{A'}\mathcal{K}(K',B)}="br";
        {\ar@{=>}^{\mathcal{M}_{A'}F^{\T}_{A'}\mathcal{K}(k,B)} "tl";"tr"};
        {\ar@{=>}_{ \Xi_{A'}\mathcal{K}(K,B)} "tl";"bl"};
        {\ar@{=>}^{ \Xi_{A'}\mathcal{K}(K',B)} "tr";"br"};
        {\ar@{=>}_{ F^{\G}_{A'}\mathcal{K}(k,B)} "bl";"br"};
 \endxy
\qquad
 \xy
 (-20,10)*+{  \hat{K}\mathcal{M}_AF^{\T}_{A}}="tl";
 (20,10)*+{\hat{K}'\mathcal{M}_AF^{\T}_{A}}="tr";
 (-20,-10)*+{ \hat{K}F^{\G}_{A}}="bl";
 (20,-10)*+{  \hat{K}'F^{\G}_{A}}="br";
        {\ar@{=>}^{  \hat{k}\mathcal{M}_AF^{\T}_{A}} "tl";"tr"};
        {\ar@{=>}_{ \hat{K}\Xi_{A}} "tl";"bl"};
        {\ar@{=>}^{ \hat{K}'\Xi_{A}} "tr";"br"};
        {\ar@{=>}_{ \hat{k}F^{\G}_{A} } "bl";"br"};
 \endxy
\]
which are both equal to
\[
 \xy
   (-30,-10)*+{\left(TsK,\Phi\big((\mu sK)\big)\right)}="bl";
   (-30,10)*+{(TsK,\mu sK)}="tl";
   (30,10)*+{(TsK',\mu sK')}="tr";
   (30,-10)*+{(TsK',\Theta\Phi(\mu)sK')}="br";
        {\ar^-{\left(Tsk, \Phi\big((\mu s_k^{-1})\big)\right)} "tl";"tr"};
        {\ar_{(TsK,\bar{\varrho}sK)} "tl";"bl"};
        {\ar^{(TsK',\bar{\varrho}sK')} "tr";"br"};
        {\ar_-{\big(Tsk,\Phi(\mu)s_k\big)} "bl";"br"};
 \endxy
\]
 \qed

\begin{prop} \label{GT2}
Let
\[
 \padjunction{B}{C}{L_1}{R}
\quad {\rm and} \quad
 \pradjunction{B}{C}{L_2}{R}
\]
(or $L_1 \dashv_p R \dashv_p L_2$)  be pseudoadjunctions in the
$\cat{Gray}$-category $\mathcal{K}$. Also, let $\T_1$ be the
pseudomonad on $B$ induced by the composite $R L_1$, and $T_2$ be
the endomorphism on $B$ induced by the composite $R L_2$. Then
$T_1 \dashv_p T_2$ are pseudoadjoint morphisms, hence $T_2$ is
with the pseudocomonad structure induced via mateship is a right
pseudoadjoint pseudocomonad for the pseudomonad $\T$.
\end{prop}

{\bf Proof .} The composites $RL_1$ and $RL_2$ of pseudoadjoints
are pseudoadjoint by Proposition~\ref{defcomppseudoadj}.  Thus, if
we let $\T_2$ be the pseudocomonad on $B$ determined via mateship
from the pseudomonad $\T_1$ then it is clear that $\T_1 \dashv_p
\T_2$. \qed

\begin{thm} \label{pMainThmII}
If $I,E,J,K,i,e,j,k\maps F \dashv_p U \dashv_p F\maps A \to B$ is
a pseudo ambijunction in the \cat{Gray}-category $\mathcal{K}$,
then the induced pseudomonad $UF$ on $B$ is Frobenius with
$\varepsilon = k$.
\end{thm}

\proof All we must show is that $UF \dashv_p UF$ with counit
$k.UiF$. Define the unit of the pseudo adjunction to be $UjF.i$.
Then $UF \dashv_p UF$ follows by
Proposition~\ref{defcomppseudoadj}.\qed

We now make use of the fact that every $\cat{Gray}$-category
$\mathcal{K}$ can be freely completed to a $\cat{Gray}$-category
\cat{EM}$(\mathcal{K})$ where an Eilenberg-Moore object exists for
every pseudomonad in $\mathcal{K}$.

\begin{thm} \label{pMainThm}
Given a Frobenius pseudomonad $(\T,\varepsilon)$ on an object $B$
in the $\cat{Gray}$-category $\mathcal{K}$, then in
\cat{EM}$(\mathcal{K})$ the left pseudoadjoint $F^{\T}\maps B \to
B^{\T}$ to the forgetful $\cat{Gray}$-functor $U^{\T}\maps
B^{\T}\to B$ is also right pseudoadjoint to $U^{\T} $with counit
$\varepsilon$. Hence, the Frobenius pseudomonad $\T$ is generated
by an ambidextrous pseudo adjunction in \cat{EM}$(\mathcal{K})$.
\end{thm}

\Proof  In \cat{EM}$(\mathcal{K})$ an Eilenberg-Moore object
exists for the pseudomonad $\T$. In particular, this means that
the $\cat{Gray}$-functor $\cat{$\T$-Alg}$ is represented by
$\mathcal{K}(-,B^{\T})$ for some $B^{\T}$ in
$\cat{EM}(\mathcal{K})$. Hence, the pseudoadjunction
$$I^{\T},E^{\T},i^{\T},e^{\T}\maps F^{\T} \dashv_p U^{\T}\maps
\cat{$\T$-Alg} \to \mathcal{K}(-,B)$$ of
Theorem~\ref{FextendsGray} arises via the enriched Yoneda lemma
from a pseudoadjunction
$$I^{\T},E^{\T},i^{\T},e^{\T}\maps F^{\T} \dashv_p U^{\T}\maps
B \to B^{\T}$$ in $\cat{EM}(\mathcal{K})$.  Furthermore, since
$\T$ is a Frobenius pseudomonad we can equip the endomorphism $T$
with the induced pseudocomonad structure of
Proposition~\ref{PROPaa}. We denote this pseudocomonad as $\G$.
Then the pseudoadjunction:
\[
I^{\G},E^{\G},i^{\G},e^{\G}\maps \overline{\mathcal{M}}F^{\G}
\vdash_p U^{\G}\mathcal{M}
 \maps \cat{$\T$-Alg} \to
\mathcal{K}(-,B)
\]
given by the construction of pseudocoalgebras composed with the
\cat{Gray}-equivalence $\cat{$\T$-Alg} \simeq \cat{$\G$-CoAlg}$
must also arise via the enriched Yoneda lemma from a
pseudoadjunction:
\[
I^{\G},E^{\G},i^{\G},e^{\G}\maps U^{\G}\mathcal{M} \dashv_p
\overline{\mathcal{M}}F^{\G} \maps B \to B^{\T}
\]
in $\cat{EM}(\mathcal{K})$.  Since this pseudoadjunction generates
the pseudocomonad $\G$, and $\G$ is defined by mateship under the
self pseudoadjunction determined by $\varepsilon$, we have that
$e^{\G}=\varepsilon$.

By Theorem~\ref{pseudoADJMONTHM} we have that
$U^{\G}\mathcal{M}=U^{\T}$.  Since $\cat{$\T$-Alg}$ is
representable in $\cat{EM}(\mathcal{K})$ the isomorphism
$\mathcal{M}F^{\T} \cong F^{\G}$ of Proposition~\ref{THMab} arises
via the enriched Yoneda lemma from an isomorphism between the
morphisms $\mathcal{M}F^{\T}$ and $F^{\G}$ in
$\cat{EM}(\mathcal{K})$. Hence, $F^{\T}\maps B \to B^{\T}$ is both
a left and right pseudoadjoint to $U^{\T}$, so that the Frobenius
pseudomonad $\T$ is induced from an ambidextrous pseudoadjunction.
\qed

\begin{cor} \label{lastcor}
A Frobenius pseudomonoid in a semistrict monoidal 2-category
$\mathcal{M}$ (or $\cat{Gray}$-monoid) yields a pseudo
ambijunction in \cat{EM}$(\Sigma(\mathcal{M}))$, where
$\Sigma(\mathcal{M})$ is the $\cat{Gray}$-category obtained from
the suspension of $\mathcal{M}$.
\end{cor}

\Proof  Recall that a Frobenius pseudomonoid in the
$\cat{Gray}$-monoid $\mathcal{M}$ is just a Frobenius pseudomonad
in the $\cat{Gray}$-category $\Sigma(\mathcal{M})$. The result
follows by Theorem~\ref{pMainThm}. \qed

\begin{cor}
If $\pseudoambijunction{B}{C}{F}{U}$ is a pseudo ambijunction in
the $\cat{Gray}$-category $\mathcal{K}$, then $UF$ is a Frobenius
pseudomonoid in the semistrict monoidal 2-category
$\mathcal{K}(B,B)$.
\end{cor}

\Proof  By Theorem~\ref{pMainThmII}, $UF$ is a Frobenius
pseudomonad on $B$ in $\mathcal{K}$.  By definition this is a
Frobenius pseudomonoid in $\mathcal{K}(B,B)$. \qed

\subsubsection*{Acknowledgements}

The author would like to thank John Baez, Richard Garner, Martin
Hyland, Steve Lack and Ross Street for their suggestions and
helpful advice. The author is also grateful for the generosity of
the University of Chicago where this work was completed.  Thanks
also to the referee for their careful editing of this paper. All
diagrams in this paper where constructed using Ross Moore and
Kristoffer Rose's \Xy-pic diagram package.

\end{document}